\documentclass[english]{article}
\usepackage[T1]{fontenc}
\usepackage[latin9]{luainputenc}
\usepackage{geometry}
\usepackage{enumitem}
\usepackage{authblk} 
\usepackage{color}
\usepackage{babel}
\usepackage{verbatim}
\usepackage{amsfonts,amsmath,amssymb,amsthm,bm, dsfont,mathtools,thm-restate}
\usepackage[linesnumbered,commentsnumbered,ruled,vlined]{algorithm2e}

\usepackage{multicol}
\usepackage{comment}
\usepackage{xparse}
\usepackage[colorlinks=true, linkcolor=blue, citecolor=red, urlcolor=purple]{hyperref}
\usepackage{cleveref}

\makeatletter

\usepackage[bordercolor=orange,backgroundcolor=orange!20,linecolor=orange,textsize=scriptsize,textwidth=1.5cm]
{todonotes}

\usepackage[most]{tcolorbox}

\newtcolorbox{mybox}[1][]{%
  enhanced,
  breakable,
  colback=blue!5,
  colframe=blue!50!black,
  boxrule=0.8pt,
  arc=2mm,
  left=2mm,
  right=2mm,
  top=1.5mm,
  bottom=1.5mm,
  #1
}

\newtheorem{thm}{Theorem}[section]
\newtheorem*{thm*}{Theorem}
\newtheorem{claim}[thm]{Claim}
\newtheorem{proposition}[thm]{Proposition}
\newtheorem{lem}[thm]{Lemma}
\newtheorem*{lem*}{Lemma}
\newtheorem{cor}[thm]{Corollary}

\theoremstyle{definition}

\newtheorem*{example*}{Example}
\newtheorem{defn}[thm]{Definition}
\newtheorem*{defn*}{Definition}

\allowdisplaybreaks
\newcommand{\ee}{\mathrm{e}}
\newcommand{\N}{\mathbb{N}}
\newcommand{\R}{\mathbb{R}}
\newcommand{\Rpp}{\mathbb{R}_{>\mathbf{0}}}
\newcommand{\Rp}{\mathbb{R}_{\ge\mathbf{0}}}
\newcommand{\Q}{\mathbb{Q}}

\newcommand{\hk}{\hat\kappa}

\newcommand{\cC}{\mathcal{C}}
\newcommand{\cD}{\mathcal{D}}
\newcommand{\cDpp}{\mathcal{D}_{>\mathbf{0}}^n}

\newcommand{\cN}{\mathcal{N}}
\newcommand{\cNc}{\mathcal{N}^{\mathcal{R}}}

\newcommand{\cP}{\mathcal{P}}
\newcommand{\cPs}{\mathcal{P}_\star}

\newcommand{\defeq}{\coloneqq}

\newcommand{\1}{\mathbf 1}
\newcommand{\0}{\mathbf 0}
\newcommand{\targetmu}{\mu_1}

\newcommand{\out}{\mathrm{out}}
\newcommand{\poly}{\mathrm{poly}}
\newcommand{\diag}[1]{\mathrm{diag}\left(#1\right)}

\newcommand{\norm}[1]{\left\lVert#1\right\rVert}

\newcommand{\wt}{\widetilde}

\newcommand{\posin}{\mathcal{I}^+}
\newcommand{\negin}{\mathcal{I}^-}

\newcommand{\length}{\mathrm{length}}
\newcommand{\what}{\widehat}
\newcommand{\pr}[2]{\left\langle #1, #2 \right\rangle}

\newcommand{\LCP}{\mathrm{LCP}}
\newcommand{\rhos}{\rho^\mathcal{S}_0}
\newcommand{\cS}{\mathcal{S}^n}
\newcommand{\cl}{\mathrm{cl}}

\newcommand{\bx}{\bar{x}}

\newcommand{\mc}{\mathrm{maxcut}}

\NewDocumentCommand{\resinvE}{m m m e{^}}{%
\ensuremath{\sigma\IfValueT{#4}{^{#4}}_{#1}(#3,#2)}%
}

\newcommand{\sgn}{\mathrm{sgn}}
\newcommand{\rk}{\mathrm{rk}}
\newcommand{\supp}{\mathrm{supp}}
\NewDocumentCommand{\detw}{m m m e{^}}{%
  \ensuremath{\upsilon\IfValueT{#4}{^{#4}}_{#1}(#3,#2)}%
}
\newcommand{\graphweight}{\rho}

\newcommand{\mmw}{\zeta}

\newcommand{\err}{\mathrm{err}}
\newcommand{\Feas}{\mathcal{F}}
\newcommand{\Feasp}{\mathcal{F}_{>\0}}
\newcommand{\Kp}{K}

\renewcommand{\emptyset}{\varnothing}
\newif\ifcomments\commentsfalse
\ifcomments
\newcommand{\lnote}[2][]{\todo[#1,backgroundcolor=green!20!white,bordercolor=blue,linecolor=green]{Laci: #2}}
\newcommand{\mnote}[2][]{\todo[#1,backgroundcolor=blue!20!white,bordercolor=blue,linecolor=green]{Haoyuan: #2}}
\newcommand{\anote}[2][]{\todo[#1,backgroundcolor=orange!20!white,bordercolor=blue,linecolor=green]{Adrian: #2}}
\else
\newcommand{\lnote}[2][]{}
\newcommand{\mnote}[2][]{}
\newcommand{\anote}[2][]{}
\fi

\makeatletter
\let\original@algocf@latexcaption\algocf@latexcaption
\long\def\algocf@latexcaption#1[#2]{%
  \@ifundefined{NR@gettitle}{%
    \def\@currentlabelname{#2}%
  }{%
    \NR@gettitle{#2}%
  }%
  \original@algocf@latexcaption{#1}[{#2}]%
}

\makeatother
\title{Newton Steps are Faster for Linear Complementarity Problems
}
\author[1]{Haoyuan Ma}
\author[1]{L{\'{a}}szl{\'{o}} A. V{\'{e}}gh}
\author[2]{Adrian Vladu}
\affil[1]{Hertz Chair for Algorithms and Optimization, University of Bonn, Germany}
\affil[2]{Department of Computer Science, Georgetown University, USA}

\date{}

\begin{document}

\maketitle
\begin{abstract}
  Linear Complementarity Problems (LCPs) with sufficient matrices form a well-studied class of problems, dating back from the foundational works of Cottle and Dantzig (1968). LCPs capture convex quadratic programs, and admit reductions from problems such as discounted stochastic games or mean-payoff games. While the original approach for solving these problems was based on pivoting methods, in their 1991 book, Kojima, Megiddo, Noma, and Yoshise introduced an interior point method (IPM) with running time that is polynomial in the input size, but which depends linearly on a certain parameter called the handicap number of the matrix. Despite the ensuing rich literature on sufficient LCP IPMs, the linear dependence on the handicap number has not been improved since then. 

  We present an IPM whose running time depends on the square root of the handicap number. The algorithm can decrease the duality gap by a constant factor  in $O\left(\min\left\{\sqrt{n} (1+\kappa),n^{9/8}\sqrt{1+\kappa}\right\}\right)$ iterations.
  Curiously, this can already be achieved using standard Newton updates. 
  We analyze the IPM in a norm defined by a local weighting, where centralizing Newton steps quadratically converge to the central path, without any dependence on the handicap number. The weights arise from a characterization by Rump (Linear Algebra and Its Applications, 1997, 1998, 2003) of the sign-real spectral radius and its relation to sufficiency; another ingredient of our proof is  a new operator norm characterization of the handicap number.
  
  Moreover, the connection to Rump's work also enables us to show coNP-completeness of testing sufficiency, and approximation hardness of the handicap number.
\end{abstract}

\section{Introduction}

Interior point methods (IPMs) are among the most powerful algorithmic tools for continuous optimization. Their success is closely tied to convexity, which provides the structure needed to control Newton steps and guarantee rapid convergence. A longstanding question is how far these methods can be extended beyond convex optimization, and what properties can replace convexity in their convergence analysis.

An important setting for studying this question is the \emph{linear complementarity problem} (LCP). Given a matrix
$M\in\R^{n\times n}$ and a vector $q\in\R^n$, the goal is to find vectors $x,s\in\R^n$ satisfying
\begin{equation}\label{eq: LCP}\tag{$\LCP(M, q)$}
Mx+q=s\,, \quad x^\top s=0\,, \quad x, s\geq \0_n\,.
\end{equation}
The LCP encompasses a remarkably broad range of problems. In particular, the optimality conditions of linear and convex quadratic programs can be expressed as LCPs, and complementarity formulations arise naturally in economic equilibria, noncooperative games, and variational inequalities.
They were first studied in the 1960s by Lemke~\cite{lemke1965bimatrix} and by Cottle and Dantzig~\cite{cottle1968complementary}, extending the complementary pivot algorithm by Lemke and Howson \cite{lemke1964equilibrium} from bimatrix games to this more general setting. We refer the reader to the books by Cottle, Pang, and Stone \cite{cottle2009linear} and by Murty and Yu \cite{murty1988} for comprehensive treatments of LCP.

The best understood case is that of \emph{monotone} LCPs, where $M+M^\top$ is positive semidefinite, or equivalently, $x^\top Mx\geq 0$ for every $x\in\R^n$. Much like convex optimization problems, monotone LCPs can be solved in polynomial time by standard interior point methods (IPMs) \cite{kojima1989polynomial, ding1991polynomial,kojima1992interior,ye1993quadratic,wright1997primal}. Unfortunately, many important LCPs do not enjoy this structure. For unrestricted matrices, solving an LCP is NP-hard \cite{Chung89}, and the complexity of several natural matrix classes remains a major open problem.

A particularly important example is the class of $\cP$-matrices, which are matrices whose principal minors are all positive. Every $\cP$-matrix $M$ has the property that $\LCP(M,q)$ admits a unique solution for every $q$, a fact established already before the explicit introduction of the LCP model \cite{samelson1958partition}. Cottle and Dantzig \cite{cottle1968complementary} showed that Lemke's algorithm finitely finds a solution for this class, yet there are still no polynomial time algorithms known. 
Megiddo \cite{Meg88} showed that solving $\cP$-matrix LCP cannot be NP-hard unless NP$=$coNP.
The significance of this problem extends well beyond mathematical optimization. $\cP$-matrix LCPs are closely connected to unique sink orientations of hypercubes \cite{stickney1978digraph, szabo2001unique} and to the complexity of total search problems: they are among the first problems shown to lie in PPAD by Papadimitriou
\cite{papadimitriou1994complexity}, and also belong to the subclass
  UEOPL of PPAD \cite{fearnley2020unique}. Prominent problems at the intersection of automata theory and  algorithmic game theory, such as discounted stochastic games, mean-payoff and parity games, admit reductions to $\cP$-matrix LCPs \cite{gartner2005simple,jurdzinski2008simple,hansen2013complexity,fearnley2020unique}. These connections make the search for efficient algorithms for structured LCPs an important direction in theoretical computer science and in optimization.

Analogously to their distinguished role in the context of pivoting algorithms, $\cP$-matrix LCPs are also amenable to interior point methods.
A seminal work by Kojima, Megiddo, Noma, and Yoshise \cite{kojima1991unified} 
 introduced the parametrized matrix class  $\cPs(\kappa)$, and gave an interior point method for $\cPs(\kappa)$-LCP with running time polynomial in the input and linear in $\kappa$.
  The parameter $\kappa$, known as the \emph{handicap number}, measures the extent to which a matrix violates monotonicity; every  $\cP$-matrix belongs to $\cPs(\kappa)$ for some $\kappa\ge 0$. However, this does not yield a polynomial-time algorithm for $\cP$-matrix LCP, because $\kappa$ can be exponentially large in the input encoding length.  

 This was followed by a long line of work on developing IPMs for $\cPs(\kappa)$-matrix LCPs, e.g., \cite{miao1995quadratically,potra1996predictor,potra1997large, peng2000new,potra2005predictor,liu2006corrector,potra2009superlinearly,peng2009self,lesaja2010unified,illes2010polynomial,enagy2024new,darvay2026comprehensive}, see also Section~\ref{sec:previous}. The literature explored a variety of IPMs using different neighborhoods, short and long-steps, as well as different kernel functions and reparametrizations. However, the best known running time bounds for $\cPs(\kappa)$-LCPs have retained a linear dependence on $\kappa$.%

 In this paper, we make the first improvement to this dependence in over three decades, showing that it can be reduced from linear to square root. %
Among extensive efforts to improve IPM iteration complexity across diverse regimes~\cite{madry2013navigating,madry2016computing,cmsv17,liu2019faster,axiotis2020circulation,liu2020faster,vladu2025breaking,liu2026near}, to the extent of our knowledge this is the first result to break a longstanding bound by reanalyzing the classical method in a norm dictated by combinatorial structure.%

\paragraph{Handicap number and LCPs for $\cPs$-matrices.}

To understand the handicap number, recall that monotonicity requires
\[
    x^\top Mx=\sum_{i=1}^n x_i(Mx)_i\geq 0\, ,
    \qquad\text{for all }x\in\R^n.
\]
Although the individual terms in this sum may be negative, their total must be nonnegative. The defining property of $\cPs(\kappa)$-matrices relaxes this requirement by allowing the negative terms to outweigh the positive ones, provided that multiplying the positive contributions by $1+4\kappa$ restores nonnegativity.
Formally, for $x\in\R^n$, define the index sets according to the sign of $x_i(Mx)_i$:
\begin{equation*}
    \posin_M(x)\defeq\left\{i\in [n]: \,x_i(Mx)_i> 0\right\}\,,\quad\mbox{and}\quad \negin_M(x)\defeq\left\{i\in[n]:\, x_i(Mx)_i<0\right\}\,. 
\end{equation*}
The \emph{handicap number} of $M$ is defined as
\begin{equation}
    \hk(M)\coloneqq\inf\left\{\kappa\geq 0:\, (1+4\kappa)\sum\limits_{i\in\posin_M(x)} x_i(Mx)_i+\sum\limits_{i\in\negin_M(x)}x_i(Mx)_i\geq 0\, ,\, \quad\forall x\in\R^n\right\}\,.\label{eq:handicap_def}
\end{equation}

We let $\cPs(\kappa)$ denote the set of matrices with $\hk(M)\leq\kappa$, and write $\cPs:=\bigcup_{\kappa\geq 0}\cPs(\kappa)$ for the class of matrices with finite handicap. Thus, $\cPs(0)$ is precisely the class of monotone matrices.

The class $\cPs$ is of great significance within complementarity theory. A deep theorem by V\"aliaho~\cite{valiaho1996p} (\Cref{thm: Valiaho theorem}) shows that it coincides with the class of \emph{sufficient matrices}, originally introduced by Cottle, Pang, and Venkateswaran \cite{cottle1989sufficient}.    A matrix $M\in\R^{n\times n}$ is \emph{column sufficient} if for all $x\in \R^n$,  $x( Mx)\leq\0_n$ implies $x( Mx)=\0_n$, and \emph{row sufficient} if $M^\top$ is column sufficient; the product between $x$ and $Mx$ is applied entry-wise. A matrix is called \emph{sufficient} if it is both column and row sufficient. 

Moreover, $\cPs$-LCP appears to be the broadest class where one may hope for polynomial-time algorithms: let $\cP_0$ denote the class of matrices whose principal minors are all nonnegative; such matrices arise as limits of sequences of $\cP$-matrices (see \Cref{lem: matrix class relation}).
Solving $\cP_0$-LCP is already NP-hard \cite[Section 3.4]{kojima1991unified}. Sufficient matrices fall strictly between these classes: $\cP\subsetneq\cPs\subsetneq \cP_0$ \cite[Section 3.2]{kojima1991unified}.

While the linear dependence on $\kappa$ has not been improved in the extensive previous literature,
there are compelling reasons to believe that this dependence is suboptimal. One particularly revealing example comes from discounted turn-based stochastic games, which can be encoded as $\cP$-matrix LCPs. Hansen and Ibsen-Jensen~\cite{hansen2013complexity} showed that, for games with $n$ states and discount factor $\gamma\in(0,1)$, the corresponding matrices have handicap numbers being 
$O\left(\frac{n}{(1-\gamma)^2}\right)$,
and gave an example where this bound is attained. Consequently, the classical handicap-based IPM bounds incur a quadratic dependence on $1/(1-\gamma)$. In contrast, specialized strategy-iteration algorithms for these games achieve a dependence that is linear in $1/(1-\gamma)$, up to logarithmic factors~\cite{hansen2013strategy}; see further discussion in \Cref{sec:previous}.

This gap suggests that the handicap-based IPM analysis may be losing substantially more than necessary. The $\sqrt{1+\kappa}$ dependence in \Cref{thm:main} below reduces the dependence on the inverse discount gap from quadratic to linear, matching the exponent achieved by specialized combinatorial algorithms; and thus provides a natural milestone for progress for IPM on $\cPs(\kappa)$-LCPs. Our work may raise hopes towards $\log(1+\kappa)$ dependence; this would imply a polynomial-time algorithm for $\cPs$-LCP, and consequently for all these applications.

\subsection{Our contributions}
\paragraph{New characterization of $\cPs(\kappa)$-matrices.}
Classical analyses of interior point methods for $\cPs(\kappa)$-LCPs build directly on the definition in~\eqref{eq:handicap_def}. Our equivalent characterization of $\hk(M)$ is given in terms of a matrix norm that arises naturally in analyses of  interior point methods. 
This result leads to simpler, more transparent proofs of known results, and is also a key ingredient in our improved analysis. For $A\in\R^{m\times n}$, we let $\|A\|_2=\|A\|_{2\to 2}=\max_{x\in\R^n\, , \|x\|_2=1}\|Ax\|$ denote the spectral norm, i.e., the $\ell_2\to\ell_2$ operator norm. Let $\cDpp\subseteq \R^{n\times n}$ denote the set of positive diagonal matrices.

\begin{restatable}{thm}{equikappathm}\label{thm:equivalent_pstark_char}
Let $M\in\R^{n\times n}$. Then, $M\in\cPs$ if and only if  $I_n+DMD$ is nonsingular for any $D\in\cDpp$ and 
$\sup_{D \in \cDpp} \left\|(I_n+DMD)^{-1}\right\|_2$ is finite. If 
 $M \in \cPs$, then
\begin{equation}\label{eq:idmd}
    \sup_{D \in \cDpp} \left\|(I_n+DMD)^{-1}\right\|_2 = \sqrt{1+\frac{\hk(M)^2}{1/4+\hk(M)}}\,.
\end{equation}
\end{restatable}
 The proof is given in Section~\ref{sec:proof-of-equikappathm}.
 Besides being the key to our improved analysis, it also leads to  simple new proofs of a theorem of V\"aliaho \cite{valiaho1997determining} showing that the handicap number is invariant under transposition, see \Cref{thm: handicap number of matrix transpose}.

\paragraph{Improved IPM analysis via algebraic and combinatorial properties.} 
 As we will see in Section~\ref{sec:basicIPM}, Newton steps for interior point methods involve  matrix inverses of the form as in \eqref{eq:idmd}.
In \Cref{thm:Newton-simple}, we give a simple analysis showing that $O(\sqrt{n}(1+\hk(M))\log(\rho))$ Newton steps suffice to reduce the optimality gap by a factor $\rho$.
Our main result improves the running time dependence to $\sqrt{1+\hk(M)}$. \begin{restatable}{thm}{main}\label{thm:main}
 Let $M\in\cPs\cap\Q^{n\times n}$, $q\in\Q^n$ of total encoding length $L$. There exists an algorithm that either finds a solution to \ref{eq: LCP} or concludes its infeasibility in total running time dominated by 
 \[ O\left(\min\left\{\sqrt{n} (1+\hk(M)),n^{9/8}\sqrt{1+\hk(M)}\right\}\cdot (L+\log(2+ \hk(M)))\right)\] many linear system solves, and in addition the running time of $O(\log(2+\hk(M)))$ linear programs of the form $Mx+q=s$, $x,s\ge\0_n$, $x_I=\0_I$, $s_J=\0_J$ for some $I,J\subseteq [n]$.%
\end{restatable}
In contrast with the more sophisticated IPM algorithms developed over the past decades, our result uses standard  Newton steps, for which we provide a significantly stronger analysis. Analyzing convergence in the Euclidean norm naturally leads to linear $\hk(M)+1$ dependence. We show that one can choose a changing weighted norm in which the relevant matrices have norm bounded independently of $\hk(M)$. Moreover, these weights can be chosen to remain bounded by $O(\sqrt{1+\hk(M)})$ and to change smoothly under Newton steps. This combination allows us to exploit Newton's method much more effectively, yielding the improved square-root dependence. We emphasize that these weights only appear in the analysis.%

The construction relies on a remarkable theorem of Rump \cite{rump1997theorems,rump1998sign}, relating the sign-real spectral radius of a matrix to the maximum geometric mean of products of entries of the matrix along directed cycles. While Rump established a strong connection between the 
 sign-real spectral radius and $\cP$-matrices (see \Cref{thm: Rump's theorem}), as far as we are aware, these tools have not been used in the $\cPs$-LCP IPM literature.
By interpreting the Newton matrix as a weighted directed graph, we obtain the required diagonal scaling from maximum product paths. Proving that these weights remain stable as the iterates change requires additional combinatorial arguments that are central to our analysis.

\paragraph{New hardness results.} 
 It is natural to ask for the complexity status of recognizing $\cPs$-matrices and computing the handicap number. Recognizing $\cP$-matrices is known to be coNP-complete~\cite{coxson1994p}. Recognizing row sufficient (and thus, column sufficient) matrices, as well as recognizing $\cP_0$-matrices, is also shown to be coNP-complete \cite{tseng2000co}. However, these do not imply the complexity status of recognizing sufficient matrices, which remains an open problem despite the wide belief that it is coNP-complete.\footnote{The papers \cite{de2011complexity} and \cite{e2024sufficient} attribute such a statement to \cite{tseng2000co}. This is claimed in particular in the abstract of \cite{e2024sufficient}, while the Introduction mentions the correct claim on column sufficiency. } Our next theorem settles this problem, by proving an even stronger sta{}tement: it is coNP-hard to recognize any class of matrices $\wt\cP$ between   $\cP$ and $\cP_0$.

\begin{thm}\label{thm:hardness-deciding-sufficiency}
    For any matrix class $
\wt \cP$ such that $\cP\subseteq\wt\cP\subseteq\cP_0$, deciding if $M\in\Q^{n\times n}$ is in $\wt\cP$ is coNP-hard. In particular, deciding if $M\in\Q^{n\times n}$ is in $\cPs$ is coNP-complete. 
\end{thm}
Interestingly, the proof also relies heavily on the notion of sign-real spectral radius. For any $M\in\Q^{n\times n}$, we can construct a matrix $\wt M\in\Q^{n\times n}$ such that $\wt M\in \wt\cP$ if and only if $M\in\cP$, through the characterization of $\cP$-matrix in terms of the sign-real spectral radius (\Cref{thm: Rump's theorem}). Then the result follows as a consequence of coNP-completeness of recognizing $\cP$-matrices.
 
Given a promise that a matrix is sufficient, computing its handicap number is still considered to be difficult. As shown in \cite{de2011complexity,e2024sufficient},
 it can grow exponentially with the matrix dimension. It was only recently proven that the logarithm of the handicap number can be polynomially bounded in the encoding length \cite{nagy2026handicap}. The best known algorithm to compute handicap number remains V\"aliaho's recursive algorithm \cite{valiaho1997determining} through enumeration over principal submatrices of principal pivot transforms, which is exponential. 
 
Our next result shows approximation hardness of computing the handicap number of an upper triangular matrix. Recall that in the \emph{Max-Cut} problem, the input is an undirected graph $G=(V,E)$, and the goal is to find  a subset $ S\subseteq V$ maximizing the number of edges between $S$ and $V\setminus S$.

\begin{thm}\label{thm:main-computing-handicap-number}
A polynomial-time multiplicative  $\alpha$-approximation for  $1+4\hk(M)$, restricted to upper triangular matrices  $ M\in\cP\cap \Q^{n\times n}$, implies a polynomial-time multiplicative $\alpha$-approximation for Max-Cut. %
\end{thm}
As shown in \cite{khot2007}, Max-Cut cannot be approximated by a factor better than the reciprocal of the approximation ratio $\approx0.878$ of the Goemans--Williamson algorithm \cite{goemans1995improved},
assuming the Unique Games Conjecture.\footnote{Recently, an LLM-generated manuscript claiming the proof of the Unique Games Conjecture was posted online~\cite{openai2026uniquegames}. If correct, this would imply the same hardness ratio under the assumption $P\neq NP$.} 
It is also known that approximating by a factor better than $17/16$ is NP-hard \cite{Has01,TSSW00}.
We prove our result by giving a construction of an upper triangular $\cP$-matrix $M_G$ that embeds the incidence matrix of a general connected  graph $G$, such that $1+4\hk(M_G)=4\mc(G)$. Our proof is motivated by the connection between the sign-real spectral radius and Max-Cut. Rump showed  \cite{rump1997theorems} that the quantity $\max_{s_1, s_2\in\{-1, 1\}^n}s_1^\top Ms_2$ coincides with the square of the sign-real spectral radius of a matrix that embeds $M$, and it was established by Poljak and Rohn \cite{poljak1993checking} that Max-Cut can be reduced to computing $\max_{s_1, s_2\in\{-1, 1\}^n}s_1^\top Ms_2$ for some integer matrix $M$. Our construction directly relates the Max-Cut problem to the handicap number of $\cPs$-matrices. Finding  efficient approximation algorithms, or establishing stronger inapproximability results is left for future work.

\subsection{Further related work}\label{sec:previous}
 The first polynomial time IPM for solving monotone LCPs was given by Kojima, Mizuno and Yoshise \cite{kojima1989polynomial}, followed by  further algorithms including \cite{ding1991polynomial,kojima1992interior, ye1993quadratic,wright1997primal}.
The  book by  Kojima, Megiddo, Noma, and Yoshise \cite{kojima1991unified} discussed above extended these methods for $\cPs$-LCP. Subsequent work developed faster path-following methods, different proximity measures, and predictor-corrector IPMs. In particular, Miao~\cite{miao1995quadratically} established an $\widetilde{O}((1+\kappa)\sqrt n)$ iteration bound, by an extension of Mizuno-Todd-Ye predictor-corrector algorithm \cite{mizuno1993adaptive}. Potra and Sheng \cite{potra1996predictor,potra1997large} developed variants with infeasible starting point.
Potra and Liu~\cite{potra2005predictor, liu2006corrector} introduced wide neighborhood algorithms which do not require the explicit handicap number as input.
Peng, Roos, and Terlaky \cite{peng2000new,peng2009self} developed IPMs that use more general barrier functions called self-regular kernel functions, and used these to derive long-step IPMs for $\cPs$-LCP nearly matching the above bound; Lesaja and Roos \cite{lesaja2010unified} extended the analysis for a broader class of eligible kernel functions.
Darvay \cite{darvay2003new} introduced the algebraic equivalent transformation technique for LP, a certain reparametrization of the central path.  This has been extended for  $\cPs$-LCP in a series of work, e.g.,~\cite{enagy2024new,darvay2020feasible}. The recent paper \cite{darvay2026comprehensive} presents a unified analysis of such algorithms.  $\cPs$-LCP has been extended to weighted linear complementarity problems by Potra~\cite{potra2012weighted,potra2016sufficient}, see also \cite{nagy2025parabolic}. Despite considerable progress in algorithm design and convergence analysis, the best general handicap dependence in all these works remained linear.

Sufficient matrices were introduced in 
\cite{cottle1989sufficient}. Characterizations were given by Cottle and Guu in \cite{cottle1992two,guu1995subclass}, in particular,  the latter paper conjectured that this class equals $\cPs$.
This was shown by V\"aliaho \cite{valiaho1996p}, see also 
 \cite{valiaho1997determining} and \cite{valiaho1996criteria} for further properties and characterizations. 
 The sign-real spectral radius of a matrix, which plays a key role in our work, was further studied by Rump~\cite{rump1998sign, rump2003p, rump2003perron}.

E.-Nagy and V\'egh \cite{nagy2026handicap} introduced the optimized handicap number $\hk^\star(M)$ that is the best achievable value of $\hk(DM)$ under a diagonal rescaling of the rows, and showed that \ref{eq: LCP} can be solved in $\mathrm{poly}(n,L,\hk^\star(M))$, with linear dependence on $\hk^\star(M)$. The algorithm uses the ellipsoid method with an IPM used as a cutting plane; replacing this by \Cref{alg:rescale-ipm-full} used in the proof of \Cref{thm:main} yields $\sqrt{1+\hk^\star(M)}$ dependence (albeit with higher polynomial terms).

In the context of \Cref{thm:main-computing-handicap-number} on approximation hardness of the handicap number, it is worth noting that for some other matrix condition measures relevant in optimization,  significantly stronger inapproximability results are known. The smallest absolute value subdeterminant \cite{khachiyan1995complexity}, the Dikin--Stewart--Todd condition number $\bar\chi_M$ \cite{tunccel1999approximating}, and the circuit imbalance measure \cite{DHNV24} cannot be approximated by a factor better than $2^{\mathrm{poly}(n)}$, and the maximum subdeterminant by a factor better than $2^{O(n)}$. It remains an interesting open question whether even $\log(2+\hk(M))$ remains similarly hard to approximate.

Among the applications of $\cPs$-LCP, we give a non-comprehensive overview of algorithms for discounted stochastic games, introduced by Shapley in 1953 \cite{Sha53}. %
Zwick and Paterson also reduced mean payoff games to this 
problem \cite{zwick1996complexity}. As mentioned, it was formulated as a $\cP$-LCP by Jurdzinski and Savani \cite{jurdzinski2008simple}. Hansen, Miltersen and Zwick \cite{hansen2013strategy} showed that 
Howard's policy iteration algorithm terminates in $O\left(\frac{m}{1-\gamma}\log\frac{n}{1-\gamma}\right)$ iterations for $n$ states and $m$ actions, each iteration requiring the solution of a Markov Decision Process LP. Independently from the discount factor, Halman \cite{Hal07} gave the first randomized subexponential algorithm, with the current best $\ee^{O(\sqrt{n \log(m/n)} )}$ by Zwick~\cite{zwick2026improved}. Very recently, quasipolynomial algorithms were announced by Soleymani and Farina \cite{soleymani2026} and by an AI company \cite{Ope26smpg}. Using the handicap bound $\kappa=O(n/(1-\gamma)^2)$ from \cite{hansen2013complexity}, our algorithm has running time $O(n^{1.625}L/(1-\gamma))$ many linear system solves. %
\subsection{Roadmap}
The rest of this paper is organized as follows.
\begin{itemize}
\item Section~\ref{sec:preliminaries} -- we cover preliminaries, such as notation, basic properties of $\cPs(\kappa)$-matrices, and their relation to combinatorial properties of directed graphs.
\item Section~\ref{sec:approach-overview} -- we give a comprehensive overview of the IPM analysis and the hardness results.
\item Section~\ref{sec:basicIPM} -- we present a simple analysis of a simple interior point method with Newton steps for $\cPs(\kappa)$-LCPs recovering the standard iteration bound with linear dependence on the handicap number.
\item Section~\ref{sec:improved} -- we introduce Rump weights, and use them to analyze the same algorithm  in a changing weighted norm. We show that Newton's method converges in this norm, improving the dependence on the handicap number from linear to square root.
\item Section~\ref{sec:two-step-analysis} -- we further improve the dependence on the dimension by analyzing two Newton steps at a time.
\item Section~\ref{sec:hardness} -- we prove hardness results for recognizing $\cPs$-matrices and approximating the handicap number, \Cref{thm:hardness-deciding-sufficiency} and \Cref{thm:main-computing-handicap-number}.
\item Section~\ref{sec: initialisation and termination} -- we show how to initialize the IPM, recover an exact solution, and run the algorithm without prior knowledge of $\hk(M)$, and complete the proof of \Cref{thm:main}.
\item Section~\ref{sec:proofs} -- we provide the deferred technical proofs, including \Cref{thm:equivalent_pstark_char}.
\end{itemize}

\section{Preliminaries}\label{sec:preliminaries}

\subsection{Notation}\label{sec:notation}

\paragraph{Vectors and matrices.}  For $n \in \mathbb{N}$, we let $[n]\coloneqq\{1,2,\ldots,n\}$. We write $\0_n$ and $\1_n$ for the all-zero and all-ones vectors of dimension $n$, respectively, and $I_n$ for the $n \times n$ identity matrix. For $i\in [n]$, we use $e_i$ to denote the unit vector at $i$ i.e. $({e_i})_i=1$ and $({e_i})_j=0$ for $j\neq i$. Scalar operations on vectors are interpreted coordinate-wise whenever this is clear from context, and $\langle x, y\rangle$ denotes the inner product $x^\top y$ for $x, y\in\R^n$.
We write $X \defeq \diag{x}$ to be the diagonal matrix with $X_{ii}=x_i$ for any vector $x$. The sets of $n\times n$ positive and nonnegative diagonal matrices are denoted by $\cDpp$ and $\cD_{\geq \0}^n$, respectively.
A symmetric matrix $A\in\R^{n\times n}$ is positive semidefinite, written $A\succeq \0_n$, if $x^\top A x\geq 0$ for all $x\in\R^n$. A matrix satisfying this inequality is called monotone, whether or not it is symmetric. For matrices with rational entries $A\in\Q^{m\times n}$, we let $\length(A)$ denote the total encoding length of $A$. For $i\in [n]$ and $A\in\R^{n\times n}$,  we use $A_{\bullet i}\in\R^n$ to denote the $i$th column of $A$, and $A_{i\bullet}$ to denote the $i$th row of $A$, respectively.

For $x, y\in\R^n$, we denote the Hadamard product as $xy=\diag {x} y=(x_iy_i)_{i\in
[n]}$, and sometimes write $x\circ y$ for clarity. Similarly $x^p$, for $p\in\mathbb{Z}$, $\sqrt{x}$, and $x/y$ are used for vectors $x\in\R^n$ pointwise. For $x\in\R^n$, and $p\in \mathbb{Z}_{\geq 1}\cup\{\infty\}$, we use $\norm{x}_p$ to denote the $\ell_p$-norm.
Whenever clear from the context we denote the $\ell_2$-norm simply as $\|x\|$.

\paragraph{Operator norms and vector norms.}
Any norm $\|\cdot\|$ on $\R^n$ induces the operator norm $\|A\|=\max_{x\neq 0}\frac{\|Ax\|}{\|x\|}$ on $\R^{n\times n}$. More generally, we write $\|A\|_{p\to q}=\max_{x\neq 0}\frac{\|Ax\|_q}{\|x\|_p}$ for the induced $\ell_p\to\ell_q$ norm. Additionally, we use  $\|A\|_{\max}$ to denote the maximum absolute value of the entries of $A$. 

\paragraph{Feasibility and centrality.}
Let $\Feas(M,q)\subseteq \Rp^{2n}$ denote the linear feasible region of \ref{eq: LCP} and $\Feasp(M, q)\subseteq \Rpp^{2n}$ denote the set of positive feasible solutions, that is,
\[
	\Feas(M,q)\coloneqq\{(x,s)\in\Rp^{2n}\,:\, s=Mx+q\}\,,\quad \Feasp(M,q)=\Feas(M,q)\cap\Rpp^{2n}\, .
\]
Whenever clear from the context, we simply write $\Feas=\Feas(M,q)$ and $\Feasp=\Feasp(M,q)$. We call $(x,s)\in\Feas(M,q)$  a \emph{feasible solution}. A solution to \ref{eq: LCP} is called a \emph{complementary feasible solution} or an \emph{LCP solution}. 
Interestingly, for row sufficient matrices, \ref{eq: LCP} is solvable if and only if $\Feas(M,q)\neq\emptyset$.\footnote{However, this cannot be directly turned to an efficient algorithm for solving \ref{eq: LCP}, because one cannot recurse on a problem of the same form after fixing $x_i=0$ or $s_i=0$ for some $i\in [n]$.} This can be derived from the LCP duality theorem \cite{fukuda1992linear}, see \cite{illes2009ep}, as well as \cite[Section 3.5]{cottle2009linear}.
\begin{thm}\label{thm:feasibile}
Let $M\in\R^{n\times n}$ be row sufficient and $q\in\R^n$. 
Then, the problem \ref{eq: LCP} has a complementary solution if and only if $\Feas(M,q)\neq\emptyset$.
\end{thm}

The \emph{central path} is defined as the set%
\[
\mathrm{CP}\coloneqq \left\{ (x,s)\in \Feasp\, ,\, x_is_i=\frac{x^\top s}{n}\,,\quad\forall i\in[n]\right\}.
\]
For $(x,s)\in\Feasp$ and $\mu>0$, the \emph{centrality error vector} of $(x,s,\mu)$ is defined as
\begin{equation}\label{eq:def-centrality-error}
     \err(x,s,\mu)\coloneqq \1_n-\frac{xs}{\mu}\, . 
\end{equation}
For $\beta\ge 0$, we define the $\ell_2$-neighborhood of the central path as
\[
	\cN^2(\beta)\coloneqq\left\{(x,s,\mu)\,:\, (x,s)\in \Feasp\, ,\, \mu\in\Rpp\, , \, \|\err(x,s,\mu)\|_2\le\beta\right\}\, .
\]

\subsection{\texorpdfstring{Background on $\cP$-, $\cPs$- and $\cP_0$-matrices}{Background on P-, Ps- and P0-matrices}}\label{sec:matrix-prelim}
We recall that $M\in\cPs$ if $M\in\cPs(\kappa)$ for some $\kappa\geq 0$, that is,
\[(1+4\kappa)\sum\limits_{i\in\posin_M(x)} x_i(Mx)_i+\sum\limits_{i\in\negin_M(x)}x_i(Mx)_i\geq 0\,. \]
The handicap number $\hk(M)$ is the infimum of such $\kappa$ for $M\in\cPs$, and $\infty$ if and only if $M\notin\cPs$.  
The following lemma is immediate from the definition of handicap number. 
\begin{lem}\label{lem:DMD}
    For any $M\in\cPs\cap\R^{n\times n}$ and $D\in\cDpp$, $\hk(DMD)=\hk(M)$. 
\end{lem}
One can also consider the 
 handicap number of $M$ with respect to a particular vector $x\in\R^n$, defined as 
 \begin{equation}\label{eq:handicap-ratio}
\hk(M,x)\defeq\begin{cases}
\infty\, &\mbox{if }x^\top M x<0 \mbox{ and }\mathcal{I}_M^+(x)=\emptyset \, ,\\
-\frac{\sum_{i\in \mathcal{I}_M^-(x)}x_i (M x)_i}{4 \sum_{i\in \mathcal{I}_M^+(x)}x_i (M x)_i}-\frac{1}{4}\, ,&\mbox{if } x^\top M x<0 \mbox{ and } \mathcal{I}_M^+(x)\neq\emptyset\, ,\\
0\, ,&\mbox{if } x^\top M x\ge 0\, .
\end{cases}
\end{equation}
With this notation, we can equivalently write $
\hk(M)=\sup_{x\in\R^{n}} \hk(M,x)$. 
In general, the supremum of $ \hk(M,x)$ may not be attained at any $x$, but it is the case for $\cP$-matrices, see \cite{valiaho1997determining}.

Recall that $M\in\R^{n\times n}$ is column sufficient if for all $x\in \R^n$, $x(Mx)\leq\0_n$ implies $x(Mx)=\0_n$, and row sufficient if $M^\top$ is column sufficient, and sufficient if both row and column sufficient. A fundamental result by V\"aliaho establishes the equivalence of $\cPs$ and sufficient matrices:
\begin{thm}[V\"aliaho \cite{valiaho1996p}]\label{thm: Valiaho theorem} A matrix $M\in\R^{n \times n}$ is sufficient if and only if $M\in\cPs$. 
\end{thm}
 
The following lemma considers the relations and properties of $\cP$, $\cPs$ and $\cP_0$-matrices; see \cite[Chapter 3]{cottle2009linear}.
\begin{lem}\label{lem: matrix class relation}
We have $\cP\subset \cPs\subset\cP_0$. Further, $M\in\cP_0$ if and only if $M+\varepsilon I_n\in\cP$ for any $\varepsilon>0$. In particular, $\cl(\cP)=\cP_0$. If $M\in\cP$, then $M$ is invertible and $M^{-1}\in\cP$.
\end{lem}
In particular, this justifies the invertibility in \Cref{thm:equivalent_pstark_char}.
\begin{cor}\label{cor:invertiblity}
    For any $M\in\cP_0$ and $D\in\cDpp$, the matrix $I_n+DMD$ is nonsingular.
\end{cor}
The following theorem by V\"aliaho which was originally established through more involved arguments in \cite{valiaho1997determining}, now follows as a corollary of \Cref{thm:equivalent_pstark_char}.
\begin{thm}[V\"aliaho \cite{valiaho1997determining}]\label{thm: handicap number of matrix transpose}%
For any $M\in\cPs\cap\R^{n\times n}$, we have $M^\top \in\cPs$, and 
$\hk(M)=\hk(M^\top)$. %
\end{thm}
\begin{proof}
Immediate from \Cref{thm:equivalent_pstark_char}, noting that 
   $\sup_{D\in\cDpp}\|(I_n+DMD)^{-1}\|_2=\sup_{D\in\cDpp}\|(I_n+DM^\top D)^{-1}\|_2$.
   \end{proof}
\subsection{Rump's characterizations of the sign-real spectral radius}
Here we provide preliminaries necessary for understanding and operating with Rump's characterization of the largest in magnitude real eigenvalue of a sufficient matrix.  Let $\cS\subseteq \R^{n\times n}$ denote the set of signature matrices i.e., diagonal matrices with diagonal entries $\pm 1$.

\begin{defn}
  For $M\in\R^{n\times n}$,  its \emph{real spectral radius} $\rho(M)$ is defined to be the maximum absolute value of its real eigenvalues, or $\rho(M)=0$ if no real eigenvalue exists. The \textit{sign-real spectral radius}  of $M$ is defined to be $\rhos(M)\coloneqq\max_{S\in\cS}\rho(SM)$.
\end{defn}
The following fundamental connection between sign-spectral radius and sufficient matrices is implied by Rump's theorem \cite{rump1997theorems} stated later as \Cref{thm: Rump's theorem}. We include the simple proof for completeness and for intuition.
\begin{lem}\label{lem:rhos-bound}
For $A\in\cP_0\cap\R^{n\times n}$,  $\rhos((I_n+A)^{-1})\le 1$.
\end{lem}
\begin{proof}
For a contradiction, assume $\rhos((I_n+A)^{-1})=\lambda>1$.  Then there
exists a diagonal $S\in\cS$ and $v\in\R^n\setminus\{\0\}$ such that $S(I_n+A)^{-1}v=\lambda v$.
Equivalently, for $v'=Sv$ one has that $(I_n+A)^{-1}Sv'=\lambda v'$.
Thus, $(I_n+A)v'=\frac{1}{\lambda}Sv'$, that is, $\left(D+A\right)v'=0$ for $D=I_n-\frac{1}{\lambda}S$. But $D$ is positive diagonal by the choice of $\lambda$, and therefore $D+A=D^{1/2}(I_n+D^{-1/2}AD^{-1/2})D^{1/2}$ is nonsingular by \Cref{cor:invertiblity}, a contradiction.
\end{proof}

A remarkable theorem by Rump relates the maximum absolute geometric cycle mean of a matrix to its sign-real spectral radius. In the context of matrices, a cycle takes the form $C=(i_1,..., i_\ell, i_{\ell+1}=i_1)$, where $i_1,..., i_\ell$ are distinct integers in $[n]$, and it has length $|C|=\ell\geq 1$. Note that $\ell=1$, in which case $C=(i, i)$ for some $i\in [n]$, i.e., a loop, is allowed.   
\begin{defn} For a matrix $M\in\R^{n\times n}$, we define the maximum absolute geometric cycle mean of $M$ to be 
\[
\mmw(M)\coloneqq\max_{\text{cycles }C}\left(\prod_{(i,j)\in C}|M_{ij}|\right)^{1/|C|}.
\]
\end{defn}

Rump \cite{rump1997theorems} proved that for any $M\in\R^{n\times n}$, $\mmw(M)$ is both upper bounded by a constant factor of $\rhos(M)$ and lower bounded by a $1/n$-factor of $\rhos(M)$, showing the close relation of these two quantities. 
\begin{thm}[{\cite[Theorem 5.5]{rump1997theorems}}]
\label{thm:rump-ineq}One has that 
\[
\frac{\mmw(M)}{3+2\sqrt{2}}\le\rhos(M)\le n\cdot\mmw(M)\,.
\]
\end{thm}
\subsection{Maximum mean cycles}
We use the following well-known result from combinatorial optimization. This is the multiplicative version of Gallai's characterization of graphs with no negative weight cycles \cite{Gallai1958}, see e.g. \cite[Chapter 5]{amo}. Such a multiplicative version is commonly used in the context of generalized flows, see e.g. \cite{truemper1977max} for an early example, and also \cite[Chapter 15]{amo}.
\begin{thm}\label{thm:reduced-costs}
Let $G=(V,E,\graphweight)$ be  a  weighted directed graph, $R\ge 1$, such that the arc weights satisfy  $\graphweight(i,j)\in [0,R]$ for $i\neq j$. Assume
$\prod_{e\in E(C)} \graphweight(e)\le 1$ for every directed cycle $C$.
Let us define
\begin{equation}
	\pi_j\coloneqq\max\left\{\prod_{e\in E(P)}\graphweight(e)\, :\, P\textrm{ is a directed walk ending in j}\right\}\, , \quad \forall j\in V\, , \label{eq:pi-def}
\end{equation}
with the convention that the product is 1 for an empty walk. Then, $\pi_j$ is well-defined, $\pi_j\in [1,R^{|V|}]$, and satisfies
\[
	\pi_i \rho(i,j) \le \pi_j \quad\forall (i,j)\in E\, .
\]
Moreover, $\pi$ can be computed in time $O(|V|^3)$.
\end{thm}

In the  context of \Cref{thm:rump-ineq}, we get the following. We note that such rescalings have been independently introduced in the linear algebra context by Fiedler and Pt\'ak \cite{FP69}, see also \cite{HJ91}.
\begin{cor}\label{prop:Q_H-bound}
Let $H\in\R^{n\times n}$, with  $\mmw(H)\leq \tau$ and $|H_{ii}|\leq \tau$ for all $i\in[n]$ for some $\tau>0$, and let $q_{H}$ be the potential defined by $q_{H}(i)\coloneqq\pi_i$ as in \eqref{eq:pi-def} for the weights $\rho(i,j)\coloneqq|H_{ij}|/\tau$. 
Then, for $Q_H=\diag{q_H}$, it holds that $\|Q_H H Q_{H}^{-1}\|_{\max}\le \tau$.  Moreover, for every $q_{H}(i)$ value, the optimum in  \eqref{eq:pi-def} is taken along a simple path.
\end{cor}

\section{Overview of approach}\label{sec:approach-overview}
Here we provide an overview of our overall approach. We split our approach for solving $\cPs$-LCP into  two parts. In \Cref{subsec:approach ipm 1} we describe a basic IPM which converges in time proportional to the handicap number $\hk(M)$, and we highlight the obstacles that prevent us from achieving faster convergence.
In \Cref{subsec: approach ipm 2} we show how Rump's  theorem (Theorem~\ref{thm:rump-ineq})  asserts the existence of a positive diagonal similarity rescaling $Q$ of a sufficient matrix $M$, which forces the inverse of $I_n+QMQ^{-1}$ to have entries bounded by a constant. This gives a good local step for the IPM, which motivates the introduction of a weighted norm to measure IPM progress.
This is followed by \Cref{subsec: staibility} in which we outline our strategy to prove the stability of $Q$ as the iterates change in the IPM.

We also highlight the key steps for the two hardness results. In \Cref{subsec: hardness overview} we describe how Rump's characterization of $\cP$-matrix through sign-real spectral radius gives a procedure of generating a matrix $\wt M$ for any given matrix $M$ so that $M\in\cP$ if and only if $\wt M\in\cPs$. After that, we introduce our construction that directly relates the handicap number to the max cut problem.

\subsection{Basic IPM overview}\label{subsec:approach ipm 1}
Throughout, let us assume $M\in\cPs(\kappa)$.
The key to any IPM lies in maintaining a coupling between  variables $x$ and  slacks $s=Mx+q$ that allows for a quick reduction in the complementarity gap. In a standard path-following IPM, we attempt to maintain the invariant
\[
x {} (Mx+q) = \mu \1_n
\]
with $x,s>0$, for a quickly decreasing centrality parameter $\mu$. We call $x$ a $\mu$-central solution. To do so, we leverage a recentering procedure which, starting from an iterate $x$ that satisfies $x{} (Mx+q) = w$, for some generic $w\in \R_{>0}^n$, updates $x \rightarrow x + \Delta x$ in an attempt to restore centrality: 
\begin{equation} \label{eq:ipm-update}
(x+\Delta x) {} (M (x +\Delta x) +q) \approx \mu \1_n\,.
\end{equation}
 To find a good update direction, one may simply linearize the left-hand side of ~\eqref{eq:ipm-update} in $\Delta x$. In our method, we use a modified Newton step, replacing $\diag{w}$ by $\mu I_n$ in this linearization. This gives the linear system
\[
\left(\mu I_n +XMX  \right)  \frac{\Delta x}{x} = \mu\1_n - w\,.
\]
To see how this makes progress, we can explicitly calculate the new centrality error  by directly expanding
\begin{align*}
(x +\Delta x){}(M(x+\Delta x)+q) - \mu\1_n &= \left(\1_n+\frac{\Delta x}{x}\right){} (w+XM\Delta x) - \mu \1_n \\
&= \left(\1_n+\frac{\Delta x}{ x}\right) {} \mu \1_n {} \left(\1_n-\frac{\Delta x}{ x}\right) - \mu \1_n\\
&= -\mu \1_n {} \left(\frac{\Delta x}{ x}\right)^{{} 2}\,.
\end{align*}
Thus, recalling \eqref{eq:def-centrality-error}, the  centrality error vector changes from $\err \coloneqq \1_n-\frac{w}{\mu}$ to 
\begin{equation}\label{eq:new-error}
\err' \coloneqq\left( \frac{\Delta x}{x} \right)^{{} 2} =  \left(  \left(I_n+\frac{1}{\mu}XMX\right)^{-1}  \err  \right)^{{} 2}\,.
\end{equation}
The same expansion shows that the update preserves $(x+\Delta x, M(x+\Delta x)+q)\in\Feasp$ whenever $\|\Delta x /x\|_\infty < 1$.
Equation~\eqref{eq:new-error}  captures the quadratic convergence of Newton's method, by showing that error contracts very quickly, provided it was small enough to begin with. For example, if $M$ is monotone, then  $\|\err'\|_2 \leq \| \err \|_2^2$, so error  quickly contracts towards $0$ if we start with $\|\err\|_2 \leq 1/2$. 

While this inequality need not hold for a general  $\cP_\star$ matrix, Theorem~\ref{thm:equivalent_pstark_char} bounds the operator norm of the inverse appearing in~\eqref{eq:new-error} by $O(\sqrt{1+\kappa{}})$. This, therefore, leads to the following error contraction:
\begin{align*}
\| \err' \|_2 = \left\|  \left(I_n+\frac{1}{\mu}XMX\right)^{-1} \err  \right\|_4^2 \leq \left\|  \left(I_n+\frac{1}{\mu}XMX\right)^{-1} \err  \right\|_2^2  
\\
\leq \left\|  \left(I_n+\frac{1}{\mu}XMX\right)^{-1}\right\|^2_2  \|  \err  \|_2^2 \leq O(1+\kappa{}) \|\err\|_2^2\,.
\end{align*}
Thus, if initially $\|\err\|_2 \leq c/(1+\kappa{})$ for a sufficiently small absolute constant $c>0$, repeated Newton steps quickly contract the error toward zero.

With this in mind, we can  see how the handicap number affects the convergence. Imagine that $x$ is a $\mu_0$-central point, and we want to use Newton's method initialized at $x$ to approach a $\mu = \mu_0 / (1+\delta)$-central point, for some preferably large $\delta > 0$. As soon as we lower the target  to $\mu$, the centrality error  becomes $\err = \1_n - \frac{\mu_0}{\mu} \1_n = -\delta \1_n$. To guarantee contraction via the preceding argument, it is enough to require $\|\err\|_2 = \delta n^{1/2} \leq c/(1+\kappa)$. Thus, this analysis permits $\delta=c/(\sqrt{n}(1+\kappa))$, recovering the classical  $\widetilde{O}( \sqrt{n} (1+\kappa{}) )$ iteration bound. A complete analysis is provided in Section~\ref{sec:basicIPM}.

\subsection{Improving the convergence analysis via Rump weights}\label{subsec: approach ipm 2}
The preceding analysis appears to make it difficult to improve the dependence on $\kappa{}$. The quadratic estimate only guarantees convergence within a small neighborhood of the central path, which limits how far we can advance in one step. However, a major limitation lies in the bound we use on $\|(I_n+\frac1\mu XMX)^{-1}\|$.
The bound in Theorem~\ref{thm:equivalent_pstark_char} concerns the \emph{worst-case} scaling $X$ of the input matrix $M$. When we take a Newton step, its permissible length is determined by the current iterate $x$, not the worst case one. This suggests measuring error contraction using the current local scaling of $M$, rather than only the (global) handicap number.

While this does not give a priori reason to yield any faster convergence, we notice that a remarkable theorem of Rump (Theorem~\ref{thm:rump-ineq}, c.f.~\cite{rump1997theorems}) asserts the existence of a similarity rescaling of $H = (I_n+\frac{1}{\mu}XMX)^{-1}$ by a  positive diagonal $Q_H = \diag{q_H}$ which makes all of its entries bounded by a constant in absolute value.

This suggests changing the metric used to measure centrality error so that the contraction depends on $\|Q_HHQ_H^{-1}\|_2$ rather than on $\|H\|_2$.  By defining an appropriate local norm $ \|Q_H x\|_2$ we can obtain a local convergence rate  that is completely independent of the handicap number. We call the entries of $q_H$ \emph{Rump weights}, and the resulting norm the \emph{Rump norm}.

Following the same idea as before, let us see how the error changes after a Newton step when measured in Rump norm. Assuming $q_H\ge\1_n$ entry-wise, we have %
\begin{equation}
\begin{aligned}\label{eq:centering-with-rump-weights}
\|Q_H \err'\|_2 &= \|Q_H (H \err)^{{} 2}\|_2\\
&\leq n^{1/2} \|Q_H H \err \|_\infty^2 \\
&= n^{1/2}  \|Q_H H Q_H^{-1} \cdot Q_H \err \|_\infty^2\\
&\leq n^{1/2} \|Q_H H Q_H^{-1}\|_{2\rightarrow \infty}^2 \|Q_H \err\|_2^2\\
&\leq O(n^{3/2}) \cdot \|Q_H\err\|_2^2\,.
\end{aligned}
\end{equation}
where the last inequality follows from $\|Q_HHQ_H^{-1}\|_{2\rightarrow\infty}=O(n^{1/2})$ because $Q_HHQ_H^{-1}$ has all entries bounded by an absolute constant. This gives quadratic contraction in the current weighted norm with no handicap dependence in the coefficient as long as $\|Q_H\err\|_2=O(1/n^{3/2})$.

There are two remaining challenges.

\begin{itemize}
\item First, the Rump weights depend on the matrix $H$ corresponding to the current iterate $x$ and target $\mu$. As $x$ changes to $x'$, the matrix changes to $H'$, and the relevant weights change to $q_{H'}$. We need to show that this change does not undo the contraction in~\eqref{eq:centering-with-rump-weights}. Indeed, in Lemma~\ref{lem:main-stability}  we show that there exist Rump weights which are sufficiently stable under Newton steps in the sense that $q_{H'} \leq 2 q_H$, so that following \eqref{eq:centering-with-rump-weights} while changing from $Q_H$ to $Q_{H'}$, the inequality $\|Q_{H'} \err'\|_2 \leq O(n^{3/2})  \|Q_H \err\|_2^2$ persists. %

\item Second,  as a result of~\eqref{eq:centering-with-rump-weights}, Newton's method is guaranteed to converge provided that the Rump norm of the error satisfies $\|Q_H \err\|_2 \leq  c/n^{3/2}$ for a sufficiently small absolute constant $c>0$. This determines how much we can advance the target centrality parameter while still being able to correct the remaining error. 

Formally, assuming that we start from exact centrality, after we dial down the target centrality parameter by a factor of $1+\delta$, the new error becomes $ \|Q_H\err\|=\|Q_H\cdot \delta \1_n\|_2 = \delta \|q_H\|_2$, where $H$ is evaluated at the new target. Lemma~\ref{cor:normbounds} bounds Rump weights by
\[
\1_n \leq q_H \leq O(\sqrt{1+\kappa{}}) \1_n\,.
\]

Thus $\|q_H\|_2 \leq O(n^{1/2} \sqrt{1+\kappa{}})$, and we can set $\delta\approx \frac1{n^{2}\sqrt{1+\kappa{}}}$, which  gives $O(n^{2}\sqrt{1+\kappa{}})$ iterations per constant-factor reduction in $\mu$, from a sufficiently central feasible starting point. 
While this provides a proof of concept for our IPM analysis, we can in fact show two strengthenings of this analysis. 
\end{itemize}
For the first one, instead of setting $\delta$ according to the Rump norm of the all-ones vector, we use a slightly larger step size $\delta \approx \frac{1}{n^{5/4} \sqrt{1+\kappa{}}}$, take a Newton step, then verify that the remaining error can be corrected because its Rump norm is $O(n^{-3/2})$.
For the latter, we show that we can set an even larger step size $\delta \approx \frac{1}{n^{9/8} \sqrt{1+\kappa{}}}$, take \emph{two} Newton steps, after which the remaining error is again bounded by $O(n^{-3/2})$ so it can be corrected.

The complete analysis that achieves $\widetilde{O}(n^{5/4} \sqrt{1+\kappa{}})$ iterations is done in Section~\ref{sec:improved}. The further improvement to $\widetilde{O}(n^{9/8} \sqrt{1+\kappa{}})$ is proved in Section~\ref{sec:two-step-analysis}.

\subsection{Bounding Rump weights and proving stability}\label{subsec: staibility}

As we saw, the analysis hinges on proving the existence of Rump weights that are bounded and do not change too aggressively when we update $x$. We never need to compute these weights in the algorithm. For both properties, we use Rump's characterization of the \emph{sign-real spectral radius}, denoted by $\rho_0^S$: the largest absolute real eigenvalue after multiplication by a diagonal sign matrix, maximized over those sign matrices.
The characterization from Theorem~\ref{thm:rump-ineq} enables us to connect matrix sufficiency to combinatorial properties of $H$, by interpreting it as a graph with weights encoded in $|H|$.
This proof is done in two steps. 
\begin{itemize}
\item First, Lemma~\ref{lem:rhos-bound} shows that for any sufficient matrix $A$, the sign-real spectral radius of $H=(I_n+A)^{-1}$ is at most $1$. Rump's theorem then shows that the maximum geometric mean of a cycle in the graph with arc weights $|H_{ij}|/\Kp$, where $\Kp = 8$, is at most
$\frac{3+2\sqrt2}{K} < \frac{3}{4}$.
Equivalently, there is a positive diagonal similarity rescaling $Q_H$ which makes all the entries in $| Q_H H Q_H^{-1}|$ bounded by $\Kp$. The diagonal entries of $Q_H$ are \emph{Rump weights} and can be explicitly computed by letting $q_H(j)$ be the largest product weight path ending at vertex $j$ (Definition~\ref{def:weight-def}).

\item Next, we need to understand how these path products change when we update $x$. Let the relative update in $x$ be  $\Delta \bar x=\Delta x/x$, and let $U=\diag{\1+\Delta\bar  x}$ and $A=\mu^{-1}XMX$, so that $H'=(I+UAU)^{-1}$. For a sufficiently small relative update, comparing the two inverses in the old scaling gives
\[
\left|(Q_HH'Q_H^{-1})_{uv}\right|
\leq \Kp \exp\left(O(\|\Delta\bar x\|_1)\right).
\]
Now look at a simple path $P$ from $i$ to $j$ attaining $q_{H'}(j)$. Multiplying this bound along the path gives
\begin{align*}
q_{H'}(j)
=\prod_{(u,v)\in P}\frac{|H'_{uv}|}{\Kp}
\leq\frac{q_H(j)}{q_H(i)}\exp\left(O(|P| \|\Delta\bar x\|_1)\right)
\leq q_H(j)\exp\left(O(n\|\Delta\bar x\|_1)\right),
\end{align*}
since $q_H(i)\ge1$ and a simple path has fewer than $n$ arcs.  A finer analysis in Lemma~\ref{lem:main-stability} removes the factor of $n$ and gives
\[
\|\log q_{H'}-\log q_H\|_\infty
\leq 16 \Kp \| \log(\1+\Delta\bar x)\|_1
=O(\|\Delta\bar x\|_1)
\]
when $\|\Delta\bar x\|_\infty\le1/2$. We use this improved bound in the IPM analysis.
\end{itemize}
This suffices to show stability under corrector steps. Using $\Delta\bar x=H\err$, $q_H \geq \1_n$ and $\|Q_HHQ_H^{-1}\|_{2\rightarrow\infty}=O(n^{1/2})$, \[\|\Delta\bar x\|_1\leq\|Q_HH\err\|_1
\le O(n^{3/2})\|Q_H\err\|_2\,.\] Thus, if $\|Q_H\err\|\leq c/n^{3/2}$ for a sufficiently small constant $c$, the improved stability bound gives $q_{H'}\leq 2 q_H$.%

Finally, to bound the weights, consider a simple path $P:j\to\cdots\to i$ attaining $q_H(i)>1$. If we replace column $j$ of $H$ by a vector $z$, the path is unchanged: only arcs entering $j$ change, and the path never returns to its starting vertex. The new arc $i\to j$ closes the path into a cycle. If the new matrix $C$ with a swapped column still has maximum cycle geometric mean at most $K$, the normalized cycle product is at most $1$, giving
\[
q_H(i)\frac{|z_i|}{\Kp}\leq 1\,.
\]
Thus, we need to understand how large a vector $z$ we can use while preserving this cycle bound. Let $\tau=\sqrt{1+\kappa{}}$. For any unit vector $v$, Lemma~\ref{lem:column-swap} allows the choice
$z=\frac{2Hv}{\Kp\tau}$. Indeed, this gives $\rho_0^S(C)\le4/3$, so Rump's theorem bounds the maximum cycle geometric mean by $(4/3)(3+2\sqrt2)<K$. Closing the cycle therefore gives
\[
q_H(i)|(Hv)_i|\leq\frac{\Kp^2}{2}\tau.
\]
When $q_H(i)=1$, the same bound follows from $\|H\|_2\le\tau$. Finally, the last arc $k\to i$ of a maximizing path satisfies $\Kp q_H(i) = q_H(k)|H_{ki}|$. Applying the preceding bound at coordinate $k$ with $v=e_i$ gives
\[
\1_n\le q_H\leq\frac{\Kp}{2}\sqrt{1+\kappa}\cdot\1_n\, ,
\]
as claimed in Lemma~\ref{cor:normbounds}.  Throughout the analysis, we use Lemma~\ref{lem:DMD}, which ensures that the  matrix $A=\frac{1}{\mu}XMX$ has the same handicap number as $M$. The Rump weights bounds and stability are proved in Sections~\ref{sec:bounding-Rump-weights} and~\ref{sec:stability-Rump-weights}.

\subsection{Hardness results}\label{subsec: hardness overview}
We recall that the sign-real spectral radius $\rhos(M)$ of a matrix $M\in\R^{n\times n}$ is the maximized real spectral radius attained by $SM$ over all $S$ that are diagonal with entries being $\pm 1$. As another remarkable result by Rump (Theorem~\ref{thm: Rump's theorem}, c.f.~\cite{rump1997theorems}), $\cP$-matrices have an equivalent relation to the sign-real spectral radius. The theorem asserts that the Cayley transform of $M$ is a $\cP$-matrix if and only if $\rhos(M)<1$.
The Cayley transform is defined as $\cC(M)\coloneqq(I_n-M)^{-1}(I_n+M)$ whenever $I_n-M$ is nonsingular; it can be seen as a  M\"obius transformation of matrices, see  \Cref{sec:hardness}.

Rump's theorem can be extended in two ways. For a rational matrix $M$, a rational separation can be found between $\rhos(M)$ and $1$ using the encoding length of $M$. More precisely, we can show that for $M\in\Q^{n\times n}$, the Cayley transform of $M$ is a $\cP$-matrix if and only if  $\rhos(M)\leq 1-\varepsilon$ for some $\varepsilon\in\Q_{>0}$ that depends on $\length(M)$. Secondly, as we have seen in \Cref{lem: matrix class relation}, any $\cP_0$-matrix is the limit of a sequence of $\cP$-matrices, i.e., $\cP_0$ is the closure of $\cP$. Then, it follows from a limiting argument that $\cP_0$-matrices can also be characterized through sign-real spectral radius, namely, the Cayley transform of $M$ is a $\cP_0$-matrix if and only if $\rhos(M)\leq 1$, noting that $\rhos(.)$ is continuous and the limit of a sequence of $\rhos(M_t)<1$ is less than or equal to $1$. 
\paragraph{Overview of the proof of \Cref{thm:hardness-deciding-sufficiency}.} We will exploit the aforementioned guarantee of the rational margin between the sign-real spectral radius  and $1$. For any given rational matrix $M$, we let $N$ be the inverse Cayley transform of $M$; see \Cref{prop: inverse Cayley transform} for its definition and existence. We know that $M\in\cP$ if and only if $\rhos(N)<1$ by Rump's theorem. Next, we pick $\alpha\in\Q$ to be an amplifier according to the rational gap between $\rhos(N)$ and $1$ if $M\in\cP$, so that in the case that $M\in\cP$, we still have  $\rhos(\alpha N)=\alpha\rhos(N)<1$, while if $M\notin\cP$ then we get  $\rhos(\alpha N)=\alpha\rhos(N)\geq \alpha>1$; see \Cref{claim: claim 1 for hardness proof} for the precise statement and its proof.  

Now, setting $\wt M$ to be the Cayley transform of $\alpha N$, we have that  $\rhos(\alpha N)<1$ implies $\wt M\in\cP$, by the standard Rump's theorem, while $\rhos(\alpha N)>1$ implies $\wt M\notin\cP_0$, by the extended Rump's theorem for $\cP_0$-matrices. Hence, for any matrix class $\wt\cP$ between $\cP$ and $\cP_0$, i.e., $\cP\subseteq\wt\cP\subseteq\cP_0$, we have $\wt M\in\wt\cP$ if $\rhos(\alpha N)<1$ and $\wt M\notin \wt\cP$ if $\rhos(\alpha N)>1$. But as shown in the last paragraph, $\rhos(\alpha N)<1$ if $M\in\cP$ and $\rhos(\alpha N)>1$ if $M\notin\cP$, so we can conclude that  $\wt M\in\wt\cP$ if and only if $M\in\cP$. 
As a result,  \Cref{thm:hardness-deciding-sufficiency} follows from the co-NP completeness of testing whether a rational matrix is in $\cP$ or not \cite{coxson1994p}. We provide a pseudo code (\Cref{alg: reduction}) for the detailed construction of $\wt M$.

\paragraph{Overview of the proof of \Cref{thm:main-computing-handicap-number}.} 
The proof of \Cref{thm:hardness-deciding-sufficiency} also indicates that there should be a connection between the sign-real spectral radius and handicap number. Rump proved that $\max_{s_1, s_2\in\{-1, 1\}^n}s_1^\top Ms_2$ is equal to the sign-real spectral radius of a matrix that embeds $M$; see \cite[Theorem 2.8]{rump1997theorems}, and the quantity $\max_{s_1, s_2\in\{-1, 1\}^n}s_1^\top Ms_2$ was shown to be NP-hard to compute via a reduction from the Max-Cut problem, by Poljak and Rohn \cite{poljak1993checking}. Here, we reduce Max-Cut to handicap number computation, by constructing a $\cP$-matrix whose handicap number determines the Max-Cut of the input graph.

In \Cref{lem: main lemma for NP hardness of handicap number}
 for any connected  graph $G=(V, E)$ with $|V|=n\geq 2$, $|E|=m$, we construct a matrix $M_G$ that embeds the node-edge incidence matrix $A_G\in \{0,\pm1\}^{n\times m}$ of $G$, so that $\hk(M_G)=\mc(G)-\frac{1}{4}$. Note that  $A_GA_G^\top\in\R^{n\times n}$ is the Laplacian matrix; we recall that $\mc(G)=\frac{1}{4}\max_{v\in\{-1,1\}^n} \|A_G^\top v\|_2^2$.
 The construction is 
 \[
      M_G\coloneqq \left(\begin{array}{ccc}
        I_m   & 2A_G^\top & \0 \\
        \0   & I_n & 2\1_n \\
        \0 &  \0 & 1\\
      \end{array}\right)\in\mathbb{Z}^{(m+n+1)\times (m+n+1)}\,.
  \]
We use the variable transformation
 $y=(x+M_G x)/2$, and let $N=(M_G-I_{m+n+1})/2$. Using $A_G^\top \1_n=\0$, one can verify that $N$ is nilpotent of degree 2, and that 
\begin{equation}\label{eq:change-of-varible-x^top Mx}
 x^\top M_Gx=\|y\|^2_2-\|N(I_{m+n+1}+N)^{-1}y\|^2_2=\|y_1\|_2^2+\|y_2\|_2^2-\|A_G^\top y_2\|_2^2-(n-1)y_3^2\,.
\end{equation}
Let us choose $y=(y_1,y_2,y_3)$ with $y_1=\0_m$, $y_3=1$ and $y_2\in\{-1, 1\}^n$ being the vector that attains $\mc(G)$, and the corresponding $x=(-A_G^\top y_2, y_2-\1_n, 1)$. Then, we get  that $x^\top M_G x=1-4\mc(G)$, and moreover, $x{} (M_G x)$ has a single positive term that is 1. Consequently, $\hk(M_G)\geq\hk(M_G, x)=\mc(G)-1/4$ by \eqref{eq:handicap-ratio}.

  The upper bound $\hk(M_G)\leq \mc(G)-1/4$ requires a more careful analysis, considering any general $x\in\R^{m+n+1}$ and $\sum_{i\in\mathcal{I}^+_{M_G}(x)}x_i(M_G x)_i$. We still make use of the same change of variable and \eqref{eq:change-of-varible-x^top Mx}, noting that we can bound $\|y_2\|_2^2-\|A_G^\top y_2\|_2^2$ in terms of max cut; see \Cref{prop:identity-for-max-cut}. %

\section[{Path following algorithm for P*(kappa)-LCP}]{Path following algorithm for $\cPs(\kappa)$-LCP}
\label{sec:basicIPM}
\subsection{The algorithm}\label{sec:basic-alg}
\begin{algorithm}[htb]
\DontPrintSemicolon
     \caption{\textsc{Newton's Method for LCP path following}}
    \label{alg:rescale-ipm-simple}
    \SetKwInOut{Input}{Input}
    \SetKwInOut{Output}{Output}
    \SetKw{And}{\textbf{and}}
        \Input{ $M\in \cPs\cap \R^{n\times n}$, $q\in\R^n$, parameters $\beta,\delta\in(0,1/2)$ and  an initial point $(x^0,s^0,\mu_0) \in \cN^2\left(\beta\right)$, and accuracy target $\targetmu> 0$.}
    \Output{ $(x^\out,s^\out)\in\Feasp$ of $\LCP(M, q)$ with $(x^\out,s^\out,\mu^\out)\in \cN^2\left(\beta\right)$ with $\mu^\out\leq \mu_1$. }
    $(x,s,\mu)\gets (x^0,s^0,\mu_0)$ \;
     \While{$\mu>\targetmu$}{ 
    $\mu \gets \frac{\mu}{1+\delta}$ \;
    Compute $(\Delta x,\Delta s)$ as in \eqref{eq: Newton-step for basic IPM} \;
    $(x,s)\gets (x+\Delta x,s+\Delta s)$  
    }
     \Return{$(x^\out,s^\out,\mu^{\out})=(x,s,\mu)$}
     \end{algorithm}

We now describe a generic version of Newton's method for solving LCPs with  $\cPs(\kappa)$-matrices. 
In each iteration, we maintain $(x, s)\in\Feasp$ and a parameter $\mu>0$ that gets decreased along the way. We use the   Newton step computed from the following equations to correct the current iterate $(x, s)$ towards satisfying the central path equation $x s=\mu\1_n$.  

\begin{defn}[Newton matrix and Newton step]\label{def:newton} Let $(x,s)\in\Feasp$ and $\mu>0$. The \emph{Newton matrix} at $(x,s,\mu)$ is defined as 
\[
   H_{x,\mu}\coloneqq \left(I_n+\frac{1}{\mu}XMX\right)^{-1}\, ,
\]
where $X=\diag{x}$, 
and the \emph{Newton step from $(x,s,\mu)$} is $(\Delta x,\Delta s)\in\R^{2n}$ defined by 
\begin{equation}\label{eq: Newton-step for basic IPM}
 \frac{\Delta x}{x}=H_{x,\mu}\err(x,s,\mu)\,,\quad \Delta s=M\Delta x\,.
\end{equation}
\end{defn}
\Cref{cor:invertiblity} guarantees that the Newton step is well defined for $M\in\cPs$.
Algorithm~\ref{alg:rescale-ipm-simple} requires as input a neighborhood parameter $\beta$,  a step size $\delta$, %
and a starting point in the $\cN^2\left(\beta\right)$ neighborhood. In each iteration,  the algorithm decreases $\mu$ by a  $1+\delta$ factor and then performs a Newton step. We discuss the issue of finding an initial solution and guessing the $\kappa$ value for the algorithm in Section~\ref{sec: initialisation and termination}. 

\subsection{Simple analysis with linear dependence on \texorpdfstring{$1+\kappa$}{1+kappa}}
In this section, we give a simple analysis with linear dependence on $1+\kappa$. 
\begin{thm}\label{thm:Newton-simple}
If $M\in \cPs(\kappa)$, for the choice $\beta\in \left[\frac{1}{8(1+\kappa)}, \frac{1}{4(1+\kappa)}\right]$, $\delta\in \left[\frac{1}{32\sqrt{n}(1+\kappa)},\frac{1}{16\sqrt{n}(1+\kappa)}\right]$,  and a suitable starting point $(x^0,s^0,\mu_0) \in \cN^2\left(\beta\right)$, \Cref{alg:rescale-ipm-simple} finds a feasible solution $(x^\out,s^\out,\mu^{\out})\in \cN^2\left(\beta\right)$ with $\mu^\out\le \targetmu$
in time dominated by 
$O\left(\sqrt{n}(1+\kappa)\log\left(\mu_0/\mu_1\right)\right)$  linear system solves.
\end{thm}

The following important identity will also be used in the later analysis:
\begin{lem}[Reducing centrality error] \label{lem:err-square} Let $M\in\cPs$,  $(x,s)\in\Feasp$ and $\mu>0$,  $(\Delta x,\Delta s)$ the solution to \eqref{eq: Newton-step for basic IPM}, and $(x',s')=(x+\Delta x,s+\Delta s)$. Then, 
\[
  \err(x',s',\mu)=\left(H_{x,\mu}\err(x,s,\mu)\right)^2\, .
\]
\end{lem}
\begin{proof}
It follows from \eqref{eq: Newton-step for basic IPM} that 
\begin{equation}\label{eq: centrality error after Newton step}
    \begin{aligned}
    \err(x',s',\mu)&=\1_n-\frac{(x+\Delta x)(s+\Delta s)}{\mu}
   =\1_n-\frac{\left(\1_n+\frac{\Delta x}{x}\right)\left(xs+XMX\frac{\Delta x}{x}\right)}{\mu}\\
   &=\1_n-\left(\1_n+\frac{\Delta x}{x}\right)\left(\left(I_n+\frac{1}{\mu}XMX\right)\frac{\Delta x}{x}+\left(\frac{xs}{\mu}-\frac{\Delta x}{x}\right)\right)\\
   &=\1_n-\left(\1_n+\frac{\Delta x}{x}\right)\left(\1_n-\frac{xs}{\mu}+\frac{xs}{\mu}-\frac{\Delta x}{x}\right)\\
   &=\left(\frac{\Delta x}{x}\right)^2=\left(H_{x,\mu}\err(x,s,\mu)\right)^2\,.
    \end{aligned}
\end{equation}

\end{proof}

The proof of  \Cref{thm:Newton-simple} follows immediately from the following lemma:

\begin{lem}[Progress lemma]\label{lem: basic IPM analysis}
Assume $M\in\cPs(\kappa)$.
  For any $\beta\in \left(0, \frac{1}{4(1+\kappa)}\right]$,  $\delta\in \left(0, \frac{\beta}{2\sqrt{n}}\right]$,   let $(x, s,\mu)\in \cN^2(\beta)$, $\mu'=\mu/(1+\delta)$ and let $(\Delta x, \Delta s)$ be computed from solving \eqref{eq: Newton-step for basic IPM} for $\mu'$, and let $(x',s')=(x+\Delta x,s+\Delta s)$. Then, 
$(x',s',\mu')\in\cN^2(\beta)$.
\end{lem}
\begin{proof}
We start by bounding $\err(x,s,\mu')$. Using $\|xs/\mu-\1_n\|\leq \beta$ from $(x, s, \mu)\in\cN^2(\beta)$ and $\delta\leq \frac{\beta}{2\sqrt{n}}$,
\[
\begin{aligned}
     \|\err(x,s,\mu')\|_2 &=\left\|\frac{xs}{\mu'}-\1_n\right\|_2\le 
     \left\|\frac{xs}{\mu}-\1_n\right\|_2+\delta\left\|\frac{xs}{\mu}\right\|_2\le {\beta}+\frac{\beta}{2\sqrt{n}}\cdot\left(1+\beta\right)\sqrt{n}\leq 2\beta\,. 
\end{aligned}
\]
From \Cref{lem:err-square}, we have $\err(x',s',\mu')=(H_{x,\mu'}\err(x,s,\mu'))^2$.
From \Cref{thm:equivalent_pstark_char}, using that $\hk(M)\le \kappa$, we get $\|H_{x,\mu'}\|_2\le \sqrt{1+\kappa}$. Thus,
\[
   \|\err(x',s',\mu')\|_2\le \left\|H_{x,\mu'}\err(x,s,\mu')\right\|^2_2 \le (1+\kappa)\|\err(x,s,\mu')\|^2\le (1+\kappa) 4\beta^2\le {\beta}\, ,
\]
and by the bound on $\beta$, $(x', s')\in\Feasp$, showing that $(x', s', \mu')\in\cN^2(\beta)$.
\end{proof}
\begin{proof}[Proof of \Cref{thm:Newton-simple}]
       Starting from $(x^0, s^0,\mu_0)\in \cN^2(\beta)$, for any $(x, s, \mu)\in\cN^2(\beta)$ with $ \mu\leq \mu_0$, we compute $(\Delta x, \Delta s)$ at $(x, s)$ for $\mu'=\mu/(1+\delta)$ from \eqref{eq: Newton-step for basic IPM}. With $\beta\in \left[\frac{1}{8(1+\kappa)}, \frac{1}{4(1+\kappa)}\right]$ and $\delta\in\left[\frac{1}{32\sqrt{n}(1+\kappa)}, \frac{1}{16\sqrt{n}(1+\kappa)}\right]$ which satisfies $\delta\leq \frac{\beta}{2\sqrt{n}}$, we are guaranteed to have $(x+\Delta x, s+\Delta s, \mu')\in\cN^2(\beta)$ by Lemma \ref{lem: basic IPM analysis}. Let $(x+\Delta x, s+\Delta s)$ and $\mu'$ be our new iterates. By the choice of $\delta$, \[\mu'\leq \left(1-\frac{1}{64\sqrt{n}(1+\kappa)}\right)\mu\,.\] Thus, it takes $O\left(\sqrt{n}(1+\kappa)\log(\mu_0/\mu_1)\right)$ iterations to obtain an iterate $( x^\out, s^\out, \mu^\out)\in\cN^2(\beta)$ with $\mu^\out\leq \mu_1$,  and each iteration requires one linear system solve.
    \end{proof}

\section{An improved analysis via Rump weights}\label{sec:improved}

In what follows we show that the standard IPM based on Newton steps actually can reduce
$\mu$ by a constant factor in $O(n^{5/4}\sqrt{1+\kappa})$
iterations, thus improving the long standing linear dependence in
the handicap number. The key ingredient is defining a new norm in
which to measure convergence, which heavily relies on  Theorem~\ref{thm:rump-ineq}.
\begin{thm}\label{thm:Newton-improved}
Let $M\in \cPs(\kappa)$, set  $\beta\in \left[\frac{1}{20480n^{3/2}}, \frac{1}{16384 n^{3/2}}\right]$ and $\delta\in \left[\frac{1}{9160n^{5/4}\sqrt{1+\kappa}}\,, \frac{1}{4580 n^{5/4}\sqrt{1+\kappa}}\right]$ in \Cref{alg:rescale-ipm-simple}. Given  an initial point $(x^0,s^0,\mu_0)\in\cNc(\beta)$, \Cref{alg:rescale-ipm-simple} finds a feasible solution $(x^\out,s^\out,\mu^\out)\in\cNc(\beta)$ with $\mu^\out\leq\targetmu$ in time dominated by 
$O\left(n^{5/4}\sqrt{1+\kappa}\cdot\log\left(\mu_0/\mu_1\right)\right)$  linear system solves.
\end{thm}
To achieve this result, we need to show that Newton steps converge in a local metric, defined by \emph{Rump weights}, which we formally define  in what follows.

\subsection{The Rump weights neighborhood}
For the proof of \Cref{thm:Newton-improved}, we introduce the following weights based on Rump's theorem (\Cref{thm:rump-ineq}). Throughout, we set the parameter $\Kp\coloneqq 8$,
but one could use any value such that $K\geq 5+2\sqrt{2}$
for the constant in \Cref{thm:rump-ineq} to get the same asymptotic bounds.
\begin{defn}\label{def:weight-def}
Let $H=(I_n+A)^{-1}$ for some sufficient matrix $A\in\cPs$. The \emph{Rump weights} are defined as $q_{H}(i)\coloneqq\pi_i$ as in \eqref{eq:pi-def} for the weights $\rho(i,j)\coloneqq|H_{ij}|/\Kp$  
and $Q_H\coloneqq \diag{q_H}$. That is, for $j\in [n]$,
\[  q_H(j):=\max\left\{\prod_{k=1}^{\ell}\frac{|H_{i_ki_{k+1}}|}{K}\,:\, {i_1,...,i_{\ell+1}=j}\textrm{ is a simple path }\right\}\,,\] 
with the convention that the product is $1$ for an empty path.
\end{defn}
By applying \Cref{prop:Q_H-bound} to the matrix $H=(I_n+A)^{-1}$, we get that the diagonal similarity transform of $H$ via Rump weights gives a matrix whose maximum entry is $O(1)$. We first show that all the diagonal entries of $H$ are between $0$ and $1$. \begin{proposition}\label{prop:diagonal-entry-H}
    Let $A\in\cPs\cap \R^{n\times n}$ and $H=(I+A)^{-1}$. Then for all $i\in [n]$, $H_{ii}\in (0, 1]$.
\end{proposition}
\begin{proof}
  Since $A\in\cPs$, $I+A$ is non-singular by \Cref{cor:invertiblity}. Let $i\in [n]$ and $h=H_{\bullet i}=He_i$, we have \[Ah=A(I+A)^{-1}e_i=e_i-(I+A)^{-1}e_i=e_i-h\,.\]  
  Hence, $h_i(Ah)_i=h_i(1-h_i)$ and $h_j(Ah)_j=-h_j^2\leq 0$ for all $j\neq i$. Since $A$ is column sufficient, $h\circ Ah\leq \0_n$ implies $h\circ Ah=\0_n$. Thus, if $h_i\notin (0, 1]$ then column sufficiency of $A$ implies $h=H_{\bullet i}=\0_n$ which is a contradiction to $H$ being non-singular.
\end{proof}

\begin{lem}\label{lem:weights-well-defined}
For $A\in\cPs\cap\R^{n\times n}$, let $H=(I_n+A)^{-1}$, then $\mmw(H)<6$, and thus $q_H$ is well-defined. Further, $\|Q_H H Q_H^{-1}\|_{\max}\le \Kp$.
\end{lem}
\begin{proof}
By \Cref{lem:rhos-bound}, $\rhos(H)\le 1$, and therefore by \Cref{thm:rump-ineq}, $\mmw(H)\le(3+2\sqrt{2})\rhos(H)\le3+2\sqrt{2}<6$. For each $i\in[n]$, the value $q_H(i)$ can be obtained along a simple path according to \Cref{prop:Q_H-bound}. For all $i\in [n]$, $|Q_HHQ_H^{-1}|_{ii}=|H_{ii}|\leq 1$ by \Cref{prop:diagonal-entry-H}. Thus as a result of \Cref{prop:Q_H-bound}, we get $\|Q_HHQ_H^{-1}\|_{\max}\leq K$.
\end{proof}
 Note that for $M\in\cPs$, $x\in\Rpp^n$, $\mu>0$, the Newton matrix $H_{x,\mu}$ defined in \Cref{def:newton} is of the form $H_{x,\mu}=(I_n+\frac{1}{\mu}XMX)^{-1}$ where $X=\diag{x}$. The matrix $\frac{1}{\mu}XMX$ is sufficient and has the same handicap number as $M$ by \Cref{lem:DMD}, and hence the above construction is applicable for $H=H_{x,\mu}$. The following weighted norm will play a key part in our analysis, and will be referred to as the \emph{Rump-norm}:
\[
     \|Q_{H_{x,\mu}} v\|_2=\sqrt{\pr{q_{H_{x,\mu}}^2}{v^2}} \, .
\]
\begin{defn}
Consider an instance of \ref{eq: LCP} with $M\in\cPs$.
For $\beta\ge0$, the \emph{Rump-norm neighborhood of the central path} is defined as
\[
	\cNc(\beta)\coloneqq\left\{(x,s,\mu)\,:\, (x,s)\in \Feasp\, ,\, \mu\in\Rpp\, , \, \|Q_{H_{x,\mu}}\err(x,s,\mu)\|_2\le\beta\right\}\, ,
\]
where $H_{x,\mu}$ is the Newton matrix  defined in \Cref{def:newton} and $\err(x,s,\mu)$ the centrality error defined in~\eqref{eq:def-centrality-error}.
\end{defn}
In \Cref{sec:bounding-Rump-weights}, we prove the following bound on the $q_H$ weights:
\begin{restatable}[Bounding Rump weights]{lem}{normbounds}
\label{cor:normbounds}
For $A\in\cPs(\kappa)$ and $H=(I_n+A)^{-1}$, let $Q_H=\diag{q_H}$ be defined as in \Cref{def:weight-def}. 
Then, 
$\|Q_H H v\|_\infty\le \frac{\Kp^2}{2}\sqrt{1+\kappa}\|v\|_2$ for all $v\in\R^n$
and $\1_n\le q_H \le \frac{\Kp}{2}\sqrt{1+\kappa} \cdot \1_n$.
\end{restatable}
Note that because of $q_H\ge \1_n$, $\cNc(\beta)\subseteq \cN^2(\beta)$ holds. From this lemma, we also get  the reverse containment:
\begin{equation}
\cNc(\beta)\subseteq \cN^2(\beta)\subseteq \cNc\left(\frac{\Kp}{2}\sqrt{1+\kappa}\cdot\beta\right)\,,\quad \forall\,\kappa\geq\hk(M)\,.\label{eq:rump-vs-narrow}
\end{equation}

We need one more bound on the stability of Rump weights, proven in Section~\ref{sec:stability-Rump-weights}.
\begin{restatable}[Stability of Rump weights]{lem}{mainstab}\label{lem:main-stability}
   For $A\in\cPs\cap\R^{n\times n}$, $u\in\Rpp^n$ and $U=\diag{u}$, let $H\coloneqq (I_n+A)^{-1}$ and $H'\coloneqq (I_n+UAU)^{-1}$.
Then, for $q_H, q_{H'}$ defined in \Cref{def:weight-def}, we have 
\begin{equation}\label{eq:log-d-bound}
       \left\|\log q_{H'}-\log q_H\right\|_\infty\leq 16\Kp\left\|\log(u)\right\|_1\,.
    \end{equation}
\end{restatable}

\subsubsection{Analyzing Newton steps in the Rump neighborhood }
Using the properties provided in this section, we can show that Newton's method converges in the Rump neighborhood at a rate independent of the  handicap number of the input matrix. Therefore, assuming we start with a point that is in a small Rump neighborhood of the central path, running Newton's method for $O(\log\log \kappa + \log \log n)$ steps restores centrality up to inverse polynomial error. The following proposition proves quadratic convergence when centrality error is $O(n^{-3/2})$ in Rump norm. While this is not directly used in the proof of \Cref{thm:Newton-improved}, it demonstrates the significance and role of Rump weights, and will be also needed for the improved analysis in \Cref{sec:two-step-analysis}. %

\begin{proposition}[Reducing centrality error]\label{lem:centralize}
Let $M\in \cPs(\kappa)$ and  $(x,s, \mu)\in \cNc(\beta)$ for 
 $\beta\in \left(0, \frac{1}{64\Kp^2 n^{3/2}}\right]$. Let $(\Delta x,\Delta s)$ be the Newton step as in \eqref{eq: Newton-step for basic IPM}, and $(x',s')=(x+\Delta x,s+\Delta s)$.
Then 
$(x',s',\mu)\in \cNc\left(2K^2 n^{3/2} \beta^2\right)$.
\end{proposition}

The proof requires the following simple technical lemma that we will also use later:
\begin{lem}\label{lem:H-err}
For a matrix $H=(I_n+A)^{-1}$ and Rump weights $q_H$, for any $v\in \R^n$ we have 
\[
   \|Q_HHv\|_\infty \leq   \Kp n^{1/2} \cdot \|Q_H v\|_2 \, ,\quad\mbox{and}\quad  \|Q_H (H v)^{{} 2}\|_2 \le \Kp^2 n^{3/2}\|Q_H v\|_2^2\, .
\]
\end{lem}
\begin{proof}
Recall that $q_H \geq \1_n$, and therefore $\|Hv\|_\infty\le \|Q_H Hv\|_\infty$.
By Lemma~\ref{lem:weights-well-defined},  $\|Q_H H Q_H^{-1}\|_{\max} \leq \Kp$. Thus, $\|Q_H H Q_H^{-1}\|_{2\rightarrow \infty} \leq \Kp n^{1/2}$, and we can bound
\[
\|Q_H Hv\|_\infty \leq  \|Q_H H Q_H^{-1} \|_{2\rightarrow \infty} \cdot \|Q_H v\|_2 \leq \Kp n^{1/2} \cdot \|Q_H v\|_2\,.
\]
The second claim then follows by
\[
     \|Q_H (H v)^{{} 2}\|_2 \leq \|H v\|_\infty \|Q_H Hv \|_2  \leq \|Q_H H v\|_\infty \|Q_H Hv\|_2 \leq n^{1/2} \|Q_H Hv\|_\infty^2\, .
\]
\end{proof}
\begin{proof}[Proof of \Cref{lem:centralize}]

Let us simplify the notation by setting $H\coloneqq H_{x,\mu}$,  $H'=H_{x',\mu}$, $\err\coloneqq\err(x,s,\mu)$,  and $\err'\coloneqq\err(x',s',\mu)$. Thus, per \eqref{eq: Newton-step for basic IPM}, $x'=x {} (\1_n+H\err)$.  By assumption,  $\|Q_{H}\err\|_2\le \beta$. Thus, our goal is to show that $\|Q_{H'}\err'\|_2\le {2\Kp^2 n^{3/2}\beta^2}$.

First, let us bound the new centrality error $\err'$ in the old $Q_H$ norm.  From \Cref{lem:err-square}, we have $\err'=(H\err)^2$. Then, \Cref{lem:H-err} gives
\begin{align}
\|Q_H \err'\|_2 &= \|Q_H (H \err)^{{} 2}\|_2 \leq K^2n^{3/2}\|Q_H\err\|_2^2\leq K^2n^{3/2}\beta^2 \,.\label{eq:qhherrp_bd3}
\end{align}
We finish the proof by showing that when performing a Newton step in this neighborhood, Rump weights increase only by at most a constant factor. Let $u \coloneqq \1_n+H\err$ denote the multiplicative change when moving from $x$ to $x'$. 
Using $\left|\log(1+\alpha)\right|\leq 2|\alpha|$ for $\alpha\geq -1/2$ and \Cref{lem:H-err},  we bound
\begin{align}
\left\|\log u\right\|_\infty \leq 2\|u-\1_n\|_\infty = 2\|H\err\|_\infty\le 2\|Q_HH\err\|_\infty \leq 2Kn^{1/2}\|Q_H\err\|_2 \leq 2 Kn^{1/2}\beta\,.
\end{align}
Then, applying Lemma~\ref{lem:main-stability} to $A\coloneqq 
\frac{1}{\mu}XMX$, and noting that $H'=(I_n+UAU)^{-1}$ for $U=\diag{u}$, we get 
\[
     \left\|\log q_{H'}-\log q_H\right\|_\infty\le 16\Kp  \left\|\log u\right\|_1 \leq 16 \Kp \cdot n \left\|\log u\right\|_\infty \leq 32 \Kp^2 n^{3/2}\beta \leq \frac{1}{2}\,, 
\]
by the range of $\beta$. Consequently, $\|q_{H'}/q_H\|_\infty\le \ee^{1/2}<2$. Following~\eqref{eq:qhherrp_bd3}, $\|Q_{H'}\err'\|\leq 2\|Q_H\err'\|\leq 2K^2n^{3/2}\beta^2$.
 \end{proof}
Therefore, starting from a point in the $O(n^{-3/2})$ Rump neighborhood, Newton's method yields a contraction in Rump norm. This alone already gives an $\widetilde{O}(n^{2} \sqrt{1+\kappa})$ bound on the number of iterations required to reduce complementarity gap by a constant factor, since reducing $\mu$ by a factor of $1+\delta$ from exact centrality introduces an error whose Rump norm is $\delta\|Q_H\1_n\|=\delta \|q_H\| = \delta \cdot O(n^{1/2} \sqrt{1+\kappa})$, per \Cref{cor:normbounds}, so we can use \Cref{lem:centralize} to restore centrality after setting $\delta=O(1/(n^2 \sqrt{1+\kappa}))$.
In what follows, we provide a tighter analysis which shows an improved bound on the iteration count.

\subsection{Proof of Theorem~\ref{thm:Newton-improved}}

The proof of \Cref{thm:Newton-improved} hinges on \Cref{lem:move-down2}. Starting from a sufficiently small Rump norm neighborhood of radius $O(n^{-3/2})$, we can reduce the centrality parameter $\mu$ by a factor of $1+\delta$, where $\delta = 1/\Omega(n^{5/4}\sqrt{1+\kappa})$, and take a single Newton step that returns to the same neighborhood. The improvement over the analysis sketched in the overview comes from bounding the centrality error \emph{after} this step, which permits a larger decrease in $\mu$.

\begin{lem}[Progress lemma]\label{lem:move-down2}
Let $M\in \cPs(\kappa)$, and $(x,s,\mu)\in \cNc\left( \beta\right)$ for $\beta\in \left(0, \frac{1}{16384 n^{3/2}}\right]$.
Let $\mu'=\frac{\mu}{1+\delta}$ for
$0<\delta\leq \sqrt{\frac{\beta}{1024 n(1+\kappa)}}$  and $(\Delta x,\Delta s)$ be the Newton step towards $\mu'$
\[
\frac{\Delta x}{x} =  H_{x,\mu'}\err(x,s,\mu')\,,\quad \Delta s = M \Delta x\, ,
\]
and let $(x',s')=(x+\Delta x,s+\Delta s)$. Then, $(x',s',\mu')\in \cNc(\beta)$.
\end{lem}
\begin{proof}
Since $(x, s, \mu)\in\cNc(\beta)$, 
$\|Q_{H_{x,\mu}} \err(x,s,\mu)\|\leq \beta$.
After scaling down the centrality parameter to $\mu'=\frac{\mu}{1+\delta}$, we apply Lemma~\ref{lem:main-stability} to $A:=\frac{1}{\mu} XM X$, $u=\sqrt{1+\delta}\1_n$ and $U = \diag{u}=\sqrt{1+\delta}\cdot  I_n$, and obtain that the new Rump weights satisfy 
\[
\left\|\log q_{H_{x,\mu'}} - \log q_{H_{x,\mu}}\right\|_\infty \leq 16\Kp n \log \sqrt{1+\delta}\leq 8\Kp n \delta\,.
\]
Thus, by the range of $\delta$, 
\begin{align}\label{eq:new-rump-weights-bound-l57}
q_{H_{x,\mu'}}\leq  \exp(8 \Kp n\delta) \cdot q_{H_{x,\mu}}\leq 2 q_{H_{x,\mu}}\,.
\end{align}
From here on to simplify notation we will denote $H \coloneqq H_{x,\mu'}$ and $H' \coloneqq H_{x',\mu'}$. With this notation and the bound from~\eqref{eq:new-rump-weights-bound-l57}, we have  that 
\begin{align}
    \|Q_{H} \err(x,s,\mu)\|_2 \leq 2\beta\,.\label{eq:eq-new-norm-l57}
\end{align}
Furthermore, one has that
\begin{equation}\label{eq:x-s-mu-prime}
\err(x,s,\mu') = \1_n - \frac{xs}{\mu'} = \1_n - (1+\delta)\frac{xs}{\mu}  = (1+\delta)\left(\1_n - \frac{xs}{\mu}\right) - \delta \1_n = (1+\delta) \cdot \err(x,s,\mu) - \delta \1_n\,.
\end{equation}
After taking a Newton step from this point, per Lemma~\ref{lem:err-square}, we reach a new point with centrality error
\begin{align}
\err(x',s',\mu')=\left(H  \err(x,s,\mu')\right)^{2} = \left(H \left( (1+\delta)\err(x,s,\mu)-\delta \1_n \right) \right)^{2}\,.\label{eq:new-err-l57}
\end{align}
We will show that this new error has small Rump norm for our   choice of $\delta$, so we could continue reducing error toward zero via repeated applications of \Cref{lem:centralize}. While it turns out that further centering steps are not necessary, this shows that controlling the Rump norm of the error after one Newton step, which is essentially determined by $\delta^2 \| Q_{H'} (H \1_n)^{{} 2}\|$, allows for a larger step size.
Formally,
using~\eqref{eq:new-err-l57} and AM-GM inequality,
we have
\begin{align}
\left\|Q_{H} \err(x',s',\mu')\right\|_2
&= \left\|Q_H \left( (1+\delta) H\err(x,s,\mu) -\delta H\1_n \right)^{ 2}\right\|_2 \nonumber\\
&\leq 2\delta^2 \left\|Q_{H} \left( H\1_n \right)^{2}\right\|_2 
+ 2(1+\delta)^2 \left\|Q_{H} \left(H \err(x,s,\mu)\right)^{ 2}\right\|_2\,.\label{eq:lem-move-down-e1}
\end{align}
First, noting that  $\|H\|_{2\rightarrow 2}\leq \sqrt{1+\kappa}$ by Theorem~\ref{thm:equivalent_pstark_char}, we bound
\begin{equation}\label{eq:H1-bound}
 \| H \1_n\|_2 \le \sqrt{(1+\kappa)n}
\end{equation}
Combining this with  $\|Q_H H \1_n\|_\infty \leq \frac{\Kp^2}{2}\cdot\sqrt{(1+\kappa)n}$ from Lemma~\ref{cor:normbounds}, we
 can bound the first term as
\begin{align}
    \left\|Q_{H} \left( H\1_n \right)^{{} 2}\right\|_2
    \leq \left\|Q_H H\1_n\right\|_\infty \left\| H \1_n\right\|_2  \leq \frac{\Kp^2}{2} (1+\kappa)n \,.\label{eq:lem-move-down-e2}
\end{align}
We can bound the error term by using \Cref{lem:H-err} and \eqref{eq:eq-new-norm-l57}: 
\begin{align}
\left\|Q_H (H\err(x,s,\mu))^{{} 2}\right\|_2\leq K^2n^{3/2}\|Q_H\err(x, s,\mu)\|^2\leq 4\Kp^2  n^{3/2}\beta^2\,.\label{eq:lem-move-down-e4}
\end{align}
Following \eqref{eq:lem-move-down-e1}, we combine \eqref{eq:lem-move-down-e2} and \eqref{eq:lem-move-down-e4} to obtain
\begin{align*}
    \|Q_H \err(x',s',\mu')\|_2 \leq \delta^2 \cdot {\Kp^2}(1+\kappa)n 
    + 8(1+\delta)^2\cdot \Kp^2 n^{3/2}\beta^2\,.
\end{align*}
Since $0<\delta \leq \sqrt{\frac{\beta}{1024 n(1+\kappa)}}$, and $0<\beta \leq \frac{1}{16384n^{3/2}}$,
we can further upper bound 
\begin{align}
 \|Q_H \err(x',s',\mu')\|_2 \leq 
 \frac{\Kp^2 \beta}{1024}  
 + 9 \Kp^2 \beta^2 n^{3/2}
 \leq \frac{\beta}{16}  + \frac{ \beta}{4} \leq \frac{\beta}{2}\,.\label{eq:lem-move-down-e5}
\end{align}
Finally, to bound the final Rump weights $q_{H'}$, we apply Lemma~\ref{lem:main-stability} to $u = \1_n+H\err(x,s,\mu')$ and $A=\frac{1}{\mu'} XMX$. Since  $\left|\log(1+a)\right|\le2|a|$ for $a\geq -1/2$, we have $\left\|\log(u)\right\|_1\le 2\|H\err(x,s,\mu')\|_1$.
From \eqref{eq:x-s-mu-prime},
\[
\|H\err(x,s,\mu')\|_1\le (1+\delta)\|H\err(x,s,\mu)\|_1+\delta\|H\1_n\|_1
\]
We bound the first term using  \Cref{lem:H-err} and the second term from \eqref{eq:H1-bound}, respectively. From \Cref{lem:H-err}, we get \[\|H\err(x, s, \mu)\|_1\leq \|Q_HH\err(x, s, \mu)\|_1\leq n\|Q_HH\err(x, s,\mu)\|_\infty\leq Kn^{3/2}\|Q_H \err(x, s, \mu)\|_2\leq 2Kn^{3/2}\beta\,,\] and from \eqref{eq:H1-bound}, we get $\|H\1_n\|_1\leq n\sqrt{1+\kappa}$. Hence, using the assumed bounds on $\delta$ and $\beta$, we get
\[
\|H\err(x,s,\mu')\|_1\le 2(1+\delta)Kn^{3/2}\beta+ \delta n\sqrt{1+\kappa}\le 2.5 Kn^{3/2}\beta + \sqrt{\frac{\beta n}{1024}}< \frac{1}{512}\, .
\]
Substituting this to \Cref{lem:main-stability}, we get 
\[
     \|\log q_{H'} - \log q_H\|_\infty \leq
     16\Kp \|\log(u)\|_1\le 32\Kp\|H\err(x,s,\mu')\|_1 <\frac{1}{2}\, ,
\]
so $q_{H'} \leq \exp(1/2)\cdot q_H \leq 2 q_H$. Together with~\eqref{eq:lem-move-down-e5}, we get $\|Q_{H'}\err(x', s', \mu')\|_2\leq 2\|Q_H\err(x', s', \mu')\|_2\leq \beta$.

\end{proof}

We are now ready to prove \Cref{thm:Newton-improved}.

\begin{proof}[Proof of \Cref{thm:Newton-improved}]

We show inductively that for the specified choices of $\beta$ and $\delta$, every iterate of \Cref{alg:rescale-ipm-simple} belongs to $\cNc(\beta)$.
This holds initially by the condition on $(x^0, s^0, \mu_0)$. Given
$(x,s,\mu)\in\cNc(\beta)$, the algorithm sets
$\mu'=\mu/(1+\delta)$ and takes one Newton step towards $\mu'$. Since $\beta$ and $\delta$ satisfy the conditions in  \Cref{lem:move-down2}, we get
$(x',s',\mu')\in\cNc(\beta)$ by \Cref{lem:move-down2}, completing the induction. By the choice of $\delta$, the number of iterations is $O(\delta^{-1} \log(\mu_0/\mu_1))=O(n^{5/4}\sqrt{1+\kappa}\log(\mu_0/\mu_1))$, giving the iteration bound.
Note that the algorithm never needs to compute the weights $q_H$, they are used only for the analysis. The runtime of every step is dominated by the Newton step, which is computed  by solving a structured linear system.
\end{proof}

\subsection{Bounding the Rump weights}\label{sec:bounding-Rump-weights}
 In what follows, we derive the technical lemmas needed to prove \Cref{cor:normbounds}. We start with a simple corollary of \Cref{thm:equivalent_pstark_char}.
\begin{cor}
\label{cor:inverse-bd} 
Let $A\in\cPs(\kappa)$, and $D\in\cDpp$. Then, 
\[
\left\Vert (D+A)^{-1}\right\Vert_2 \le \frac{\sqrt{1+\kappa}}{\min_{i\in[n]}d_i}\, .
\]
\end{cor}
\begin{proof}
We write $D+A=D^{1/2}(I_n+D^{-1/2}AD^{-1/2})D^{1/2}$. Noting that the expression in \Cref{thm:equivalent_pstark_char} can be upper bounded by $\sqrt{1+\kappa}$, we get 
\[
\|(D+A)^{-1}\|_2\le\|D^{-1/2}\|^{2}_2\sqrt{1+\kappa}\le \frac{\sqrt{1+\kappa}}{\min_{i\in[n]}d_i}\, .	
\]\end{proof}

We use this to derive the next technical lemma that shows stability of $\rhos((I+A)^{-1})$ under column changes.
\begin{lem}\label{lem:column-swap}
Let $A\in\cPs(\kappa)$, $H\coloneqq (I_n+A)^{-1}$ and $u\in\R^n$. For $k\in [n]$, let $M\in\R^{n\times n}$ be defined by replacing the $k$-th column of $H$ by $Hu$; that is, $M_{\bullet k}= Hu$ and $M_{\bullet i}=H_{\bullet i}$ for $i\neq k$. Then, 
\[
\rhos(M)\le 1+\sqrt{1+\kappa}\cdot \|u\|_2\,.	
\]
\end{lem}
\begin{proof}
The claim is trivial if $\rhos(M)\le 1$, so we assume $\rhos(M)>1$.
Let $\lambda\coloneqq\rhos(M)$. Thus, there exists a diagonal $S\in\text{diag}\{\pm1\}^{n}$ and $v\in\R^n$ such that $SMv=\lambda v$. Letting $v'=Sv$, this is equivalent to $MSv'=\lambda v'$. Multiplying both sides by $H^{-1}=(I+A)$ from the left and rearranging, we get
\[
\left(I_n+A-\frac{1}{\lambda} H^{-1}MS\right)v'=\0_n\, .
\]
Note that $(H^{-1}M)_{\bullet i}=e_i$ for $i\neq k$ and $(H^{-1}M)_{\bullet k}=H^{-1}Hu= u$. Defining the diagonal matrix $D\in\R^{n\times n}$ by $D_{ii}\coloneqq 1-\frac{S_{ii}}{\lambda}>0$ for $i\neq k$ and $D_{kk}\coloneqq 1$, we get 
\[
\left(D+A\right)v'=\frac{S_{kk}v'_k}{\lambda}\cdot u \, ,
\]
that is,
\[
	v' = \frac{S_{kk} v'_k}{\lambda} \cdot \left(D+A\right)^{-1} u\, .
\]
By 
\Cref{cor:inverse-bd}, we have $\|\left(D+A\right)^{-1}\|_2\le \frac{\lambda}{\lambda-1}\cdot \sqrt{1+\kappa}$, and thus
\[
	\|v'\|_2\leq \frac{|S_{kk}v_k'|}{\lambda}\cdot \|(D+A)^{-1}\|_2\cdot \|u\|_2\le  \frac{|v_k'|}{\lambda-1}\cdot \sqrt{1+\kappa} \cdot \|u\|_2\, , 
\]
implying $1\le \frac{1}{\lambda-1} \cdot \sqrt{1+\kappa} \cdot\|u\|_2$, and hence, $\lambda\le 1+\sqrt{1+\kappa}\|u\|_2$, as claimed.
\end{proof}
We can now prove \Cref{cor:normbounds}.
\normbounds*
\begin{proof}
We can w.l.o.g. assume that $\|v\|_2=1$ and let $\tau\coloneqq\sqrt{1+\kappa}$.
We first prove that
\begin{equation}\label{eq:two-to-infty}
     \|Q_H H v\|_\infty\le \frac{\Kp^2}{2}\cdot\tau\cdot\|v\|_2\,, \quad\forall v\in\R^n.
\end{equation}
Let us fix some index $\ell\in [n]$. Thus, we need to prove $q_H(\ell) |(Hv)_\ell|\le \Kp^2\tau/2$. Recall that $(q_H)_\ell\ge 1$. If $q_H(\ell)=1$, then this follows since $\|Hv\|_2\le \tau$ according to \Cref{thm:equivalent_pstark_char}. Otherwise, $q_H(\ell)$ is obtained as the product of the weights $w_{ij}=|H_{ij}|/\Kp$ along a simple path $P$ ending in $\ell$; let $k$ denote the starting point of this path.

Let us now construct the matrix $M$ by replacing the $k$-th column of $H$ by $2Hv/(\Kp\tau)$. Then by \Cref{lem:column-swap}, with $u=2v/(\Kp\tau)$ where $\|v\|_2=1$, we have \[\rhos(M)\le 1+\tau\|u\|_2\leq 1+\frac{2}\Kp\,.\] Then using \Cref{thm:rump-ineq}, $\mmw(M)\leq(\Kp-2)\rhos(M)\leq (\Kp-2)(1+2/\Kp)< \Kp$.
Let us define $w'_{ij}=|M_{ij}|/\Kp$. For all $i\in [n]$, we have $w'_{ij}=w_{ij}$ for $j\neq k$, and $w'_{ik}\coloneqq 2|(Hv)_i|/(\Kp^2\tau)$. Consider now the cycle $C$ obtained by concatenating the $k-\ell$ path $P$ with the arc $(\ell, k)$. We have
\[
	1>\frac{\mmw(M)}{\Kp}\ge \prod_{e\in C} w'_e=\frac{2\left|(Hv)_\ell\right|}{\Kp^2\tau}\cdot\prod_{e\in P} w_e = \frac{2\left|(Hv)_\ell\right|}{\Kp^2\tau}\cdot q_H(\ell)\, .
\]
Thus, $q_H(\ell) |(Hv)_\ell|\le \Kp^2\tau/2$ as claimed.

The first part of the lemma follows since $\|Q_H Hv\|_\infty=\max_{\ell\in [n]}q_H(\ell)|(Hv)_\ell|\le \frac{\Kp^2}{2}\tau$.
For the second part, let $\ell\in [n]$. We have $q_H(\ell)\ge 1$ by definition. If $q_H(\ell)>1$, then $q_H(\ell)=q_H(j) |H_{j\ell}|/\Kp$ for some $j\neq \ell$. Setting $v=e_\ell$,  \eqref{eq:two-to-infty} yields 
\[
	\frac{\Kp^2}{2}\tau\ge q_H(j) \left|(Hv)_j\right| = q_H(j) \left|H_{j\ell}\right|=\Kp q_H(\ell)\, ,
\]
which implies $q_H(\ell)\le \frac{\Kp}{2}\tau$. Since this holds for all $\ell\in [n]$, $q_H\leq \frac{K}{2}\tau \1_n$.
\end{proof}

\subsection{Stability of Rump weights}\label{sec:stability-Rump-weights}

\mainstab*
To prove the lemma, we will rely on a few helper statements. Let $z\coloneqq\log u$, and $Z\coloneqq\diag{z}$. For $t\in[0,1]$, define $A(t)\coloneqq\ee^{tZ}A\ee^{tZ}$, $H(t)\coloneqq(I_n+A(t))^{-1}$, and $d(t)\coloneqq q_{H(t)}$.

Since $A(t)$ is sufficient, the matrix $H(t)$ is well defined. Note that $H(0)=H$ and $H(1)=H'$. We will bound the entries of $H(1)$ by studying its derivatives along the path $H(t)$. To do so, we first differentiate $A(t)$, which gives
\begin{equation}\label{eq:derivative-of-A(t)}
    \dot A(t)=ZA(t)+A(t)Z\,.
\end{equation}
Using this fact we can study the change of $H(t)$.
\begin{claim}\label{cl:H(t)-diff} The matrix $H(t)$ is continuously differentiable on $t\in [0, 1]$ and
\[
    \dot H(t) =-H(t)Z-ZH(t)+2H(t)ZH(t)\,,\quad \textnormal{for all } t\in [0,1]\, .
\]
\end{claim}
\begin{proof}
Recalling that $H(t)=(I_n+A(t))^{-1}$, we differentiate both sides of $(I_n+A(t))H(t)=I_n$, giving $\dot A(t)H(t)+(I_n+A(t))\dot H(t)=\0$. 
Then, using \eqref{eq:derivative-of-A(t)}, we have  \[\dot H(t)=-(I_n+A(t))^{-1}\dot A(t) H(t)=-H(t)\dot A(t)H(t)=-H(t)(ZA(t)+A(t)Z)H(t)\,.\]
Since $A(t)H(t)=H(t)A(t)=I_n-H(t)$, \begin{equation}\label{eq:derivative-H(t)}\begin{aligned}
  \dot H(t)&=-H(t)ZA(t)H(t)-H(t)A(t)ZH(t) \\
  &=-H(t)Z(I_n-H(t))-(I_n-H(t))ZH(t)\\
  &=-H(t)Z-ZH(t)+2H(t)ZH(t)\,,
\end{aligned} 
\end{equation}
as needed.
\end{proof}
We write $w_P(t)$ for the weight of a simple path $P$ in \Cref{def:weight-def}:
\[
    w_P(t)\coloneqq\prod_{r=1}^{\ell}\frac{|H_{i_ri_{r+1}}(t)|}{\Kp}\,,\qquad P=(i_1,\ldots,i_{\ell+1}=j)\,.
\]
By that definition, $d_j(t)=\max_P w_P(t)$ over simple paths ending at $j$, including the empty path of weight $1$.

\begin{proposition}\label{prop:differetiable-of-d(t)}
For every $j\in[n]$, both $d_j$ and $\log d_j$ are Lipschitz on $[0,1]$. For any $t\in(0,1)$ where $d_j$ is differentiable, every maximizing path $P$ satisfies
\[ 
    \dot d_j(t)=\dot w_P(t)\,.
\]
In particular, $\dot d_j(t)=0$ if $d_j(t)=1$.
\end{proposition}
\begin{proof}
Each matrix entry is continuously differentiable on $[0,1]$, so each path weight $w_P$ is Lipschitz. The maximum is over finitely many paths, hence $d_j$ is Lipschitz, and so is $\log d_j$, since $d_j\ge1$.
Fix a differentiable point $t\in(0,1)$ and a maximizing path $P$. Since $w_P(t)=d_j(t)\ge1$, none of its factors vanish at $t$, so $w_P$ is differentiable there. The function $d_j-w_P$ is nonnegative and vanishes at $t$, so its derivative there is zero. If $d_j(t)=1$, apply this to the empty path.
\end{proof}

The following bound controls the derivative along a maximizing path.
\begin{lem}\label{lem:constant-bound-on-entryproduct}
Let $t\in[0,1]$ and let $P=(i_1,\ldots,i_{\ell+1}=j)$ be a nonempty simple path attaining $d_j(t)=w_P(t)$. Then, for every $r\in[\ell]$,
\[
    \frac{d_{i_r}(t)}{d_{i_{r+1}}(t)}\left|H_{i_ri_{r+1}}(t)\right|=\Kp\,.
\]
Furthermore, for every $k\in[n]$,
\[
    \sum_{r=1}^{\ell}\frac{d_{i_r}(t)}{d_{i_{r+1}}(t)}
    \left|H_{i_rk}(t)H_{ki_{r+1}}(t)\right|\le7\Kp^2\,.
\]
\end{lem}

We first prove the combinatorial bound needed for this lemma.
\begin{lem}\label{lem:graph-weight}
Let $G=(V,E,\rho)$ be a directed graph, possibly with self-loops, with weights $\rho\in[0,1]^E$. Suppose every directed cycle $C$ satisfies
\[
    \prod_{e\in C}\rho(e)\le(3/4)^{|C|}\,.
\]
Let $P=(i_1,\ldots,i_\ell)$ be a simple path whose arcs all have weight $1$. Then, for every $r\in V$,
\[
    \sum_{j=1}^{\ell-1}\rho(i_j,r)\rho(r,i_{j+1})\le7\,.
\]
\end{lem}
\begin{proof}
Let $a_j\coloneqq\rho(i_j,r)$ and $b_j\coloneqq\rho(r,i_j)$. The cycle bound also holds for closed walks, by decomposing them into cycles. Applying it to $(r,i_k,\ldots,i_j,r)$ gives
\begin{equation}\label{eq:a-b-bound}
    a_jb_k\leq(3/4)^{j-k+2}\qquad(1\le k\leq j\leq\ell)\,.
\end{equation}
This remains valid when $r$ lies on $P$.

Let $t\in\{1,\ldots,\ell\}$ be the smallest such that $a_j\le(3/4)^{j-t+1}$ for all $t\le j\le\ell-1$. Such a $t$ exists, since $t=\ell$ gives an empty range. Then
\[
    \sum_{j=t}^{\ell-1}a_jb_{j+1}\leq\sum_{h=1}^{\infty}(3/4)^h=3\,.
\]
If $t=1$, this proves the claim. Otherwise, minimality gives some $j^\star\in\{t-1,\ldots,\ell-1\}$ with $a_{j^\star}>(3/4)^{j^\star-t+2}$. By \eqref{eq:a-b-bound}, $b_k\le(3/4)^{t-k}$ for $k\le t-1$, so
\[
    \sum_{j=1}^{t-2}a_jb_{j+1}\le\sum_{h=1}^{\infty}(3/4)^h=3\,.
\]
These two bounds cover all terms except $j=t-1$. Since $a_{t-1},b_t\in[0,1]$, we conclude that
\[
    \sum_{j=1}^{\ell-1}a_jb_{j+1}
    =\sum_{j=1}^{t-2}a_jb_{j+1}+a_{t-1}b_t+\sum_{j=t}^{\ell-1}a_jb_{j+1}
    \le3+1+3=7\,.
\]
\end{proof}
Now we are ready to prove \Cref{lem:constant-bound-on-entryproduct}, which we will then use to conclude the main stability lemma.
\begin{proof}[Proof of \Cref{lem:constant-bound-on-entryproduct}]
To simplify notation, for the rest of the proof we will drop the time parameter $t$.  Let $\rho(a,b)\coloneqq (d_a/d_b) \cdot |H_{ab}| / \Kp$ for all $a, b\in [n]$.
By \Cref{lem:weights-well-defined}, these weights lie in $[0,1]$. Since $P$ attains $d_j$,
\[
    \prod_{r=1}^{\ell}\rho(i_r,i_{r+1})
    =\frac{d_{i_1}}{d_j}w_P=d_{i_1}\ge1\,.
\]
Every factor is at most $1$, so every arc of $P$ has weight $1$, which proves the first claim.
For every directed cycle $C$, the factors $d_a/d_b$ cancel, and \Cref{lem:weights-well-defined}, with $\Kp=8$, gives
\[
    \prod_{(a,b)\in C}\rho(a,b)
    =\prod_{(a,b)\in C}\frac{|H_{ab}|}{\Kp}
    \le\left(\frac{\mmw(H)}{\Kp}\right)^{|C|}
    \le(3/4)^{|C|}\,.
\]
Applying \Cref{lem:graph-weight} with $r=k$, we obtain
\[
    \sum_{r=1}^{\ell}\frac{d_{i_r}}{d_{i_{r+1}}}\left|H_{i_rk}H_{ki_{r+1}}\right|
    =\Kp^2\sum_{r=1}^{\ell}\rho(i_r,k)\rho(k,i_{r+1})\le7\Kp^2\,.
\]
\end{proof}

We are now ready to prove the main statement of this section,  \Cref{lem:main-stability}.
\begin{proof}[Proof of \Cref{lem:main-stability}]
Let $j \in [n]$ and $t \in (0,1)$ such that $d_j$ is differentiable at $t$. If $d_j(t)=1$, then $(\log d_j)'(t)=0$ by \Cref{prop:differetiable-of-d(t)}. Otherwise, there is a non-empty maximizing path $P=(i_1,\ldots,i_{\ell+1}=j)$, in particular, all the entries of $H$ along $P$ are nonzero. 
By \Cref{cl:H(t)-diff}, $\dot H(t)=-H(t)Z-ZH(t)+2H(t)ZH(t)$. Since $Z$ is diagonal, after suppressing the time parameter $t$, for all $r\in [\ell]$ we have $(ZH)_{i_ri_{r+1}}=z_{i_r}H_{i_ri_{r+1}}$, $(HZ)_{i_ri_{r+1}}=z_{i_{r+1}}H_{i_ri_{r+1}}$ and $(HZH)_{i_ri_{r+1}}=\sum_{k=1}^n z_k H_{i_rk}H_{ki_{r+1}}$. Hence, since $H_{i_ri_{r+1}}\neq 0$ for all $r\in [\ell]$,  \begin{equation*}
    \frac{\dot H_{i_ri_{r+1}}}{H_{i_ri_{r+1}}}=-z_{i_r}-z_{i_{r+1}}+2\sum\limits_{k=1}^n z_k\cdot \frac{H_{i_rk}H_{ki_{r+1}}}{H_{i_ri_{r+1}}}\,.
\end{equation*}
Using $(\log|H_{i_ri_{r+1}}|)^{\prime}=\frac{\dot H_{i_ri_{r+1}}}{H_{i_ri_{r+1}}}$, as $H_{i_ri_{r+1}}\neq 0$, and $\log w_P=\sum_{r=1}^\ell\log |H_{i_ri_{r+1}}|-\ell\log\Kp$, we get
\[
   (\log w_P)'=\sum\limits_{r=1}^\ell \frac{\dot H_{i_ri_{r+1}}}{H_{i_ri_{r+1}}}
    =-\sum_{r=1}^{\ell}(z_{i_r}+z_{i_{r+1}})
    +2\sum_{k=1}^n z_k\sum_{r=1}^{\ell}
    \frac{H_{i_rk}H_{ki_{r+1}}}{H_{i_ri_{r+1}}}\,.
\]
Each vertex occurs at most twice in the first summation, so $\left|\sum_{r=1}^\ell z_{i_r}+z_{i_{r+1}}\right|\leq 2\|z\|_1$. For the second summation, we apply both parts of \Cref{lem:constant-bound-on-entryproduct}, obtaining
\[
    \sum_{r=1}^{\ell}\left|\frac{H_{i_rk}H_{ki_{r+1}}}{H_{i_ri_{r+1}}}\right|
    =\frac{1}{\Kp}\sum_{r=1}^{\ell}\frac{d_{i_r}}{d_{i_{r+1}}}\left|H_{i_rk}H_{ki_{r+1}}\right|
    \le7\Kp\,.
\]
Consequently, \Cref{prop:differetiable-of-d(t)} yields, almost everywhere,
\[
    \left|(\log d_j)'\right|\le(2+14\Kp)\|z\|_1\le16\Kp\|z\|_1\,.
\]
Since $\log d_j$ is Lipschitz, integration over $[0,1]$ gives
\[
    \left|\log d_j(1)-\log d_j(0)\right|\leq 16 \Kp\|z\|_1\,.
\]
Using $d(0)=q_H$, $d(1)=q_{H'}$, and $z=\log u$ proves \eqref{eq:log-d-bound}.
\end{proof}

\section{\texorpdfstring{Improving the analysis to $O(n^{9/8} \sqrt{1+\kappa })$ iterations}{Improving the analysis to O(n\^{9/8} sqrt(kappa + 1)) iterations}}\label{sec:two-step-analysis}

In this section we provide an improved analysis of the IPM. While the conceptual contribution of \Cref{thm:Newton-improved} was to show that after taking a Newton step from a $\mu$ central point $x$ for the target centrality parameter $\mu'=\mu/(1+\delta)$, the resulting centrality error of the new point $x'$  is $\|Q_{H_{x',\mu'}} (H_{x,\mu'}\cdot (\delta  \1))^{ 2}\|_2 = O(\delta^2 n (1+\kappa))$, which allows setting $\delta = 1/\Omega(n^{5/4}\sqrt{1+\kappa})$ to stay in the $O(n^{-3/2})$ Rump neighborhood, here we show that two Newton steps after taking an even larger target centrality update $\delta = 1/\Omega(n^{9/8}\sqrt{1+\kappa})$ produce an iterate in the $O(n^{-3/2})$ Rump neighborhood. To do so we essentially need to bound the centrality error after two Newton steps from  error $\delta \1$. The key of this section is to show that it is $O(\delta^4 n^3 (1+\kappa)^2)$.

\begin{lem}[Centrality error after Newton steps]\label{lem:two-step-progress}
Let $M\in\cPs(\kappa)$, and let scalars $\beta=\frac{1}{4096n^{3/2}}$, $\delta=\frac{1}{1024n^{9/8}\sqrt{1+\kappa}}$, and $\varepsilon=\frac{\beta}{2^{17}n(1+\kappa)^2}$.
Suppose $(x,s,\mu)\in\cNc(\varepsilon)$, and set $\mu'=\mu/(1+\delta)$. Two successive Newton steps toward $\mu'$, as in \eqref{eq: Newton-step for basic IPM}, produce iterates $(x',s'),(x'',s'')\in\Feasp$ with $(x'',s'',\mu')\in\cNc(\beta)$.
\end{lem}

To achieve this result, we need to first prove a few helper lemmas.

\begin{lem}\label{lem:upper-bound-column-row-product}
    Let $M\in\cPs$ and $H=(I+M)^{-1}$. Then for all $j\in [n]$, \[\|H_{\bullet j}\|_2\|H_{j\bullet}\|_2\leq 2H_{jj}\sqrt{1+\hk(M)}\,.\] 
\end{lem}

\begin{proof}
    Let $t\geq 0$. Increasing the diagonal entries of a sufficient matrix  by a nonnegative amount can only decrease its handicap number; see \cite[Theorem 3.6]{valiaho1997determining}. Together with \Cref{thm:equivalent_pstark_char}, this shows that 
      \begin{equation}
        \|(I+M+te_je_j^\top)^{-1}\|_2\leq 
    \sup\limits_{D\in\cDpp}\|(I+D(M+te_je_j^\top)D)^{-1}\|_2\leq \sqrt{1+\hk(M+te_je_j^\top)}\leq \sqrt{1+\hk(M)}\,.
    \end{equation}
    By Sherman-Morrison formula, \[(I+M+te_je_j^\top)^{-1}=(I+M)^{-1}-t\frac{(I+M)^{-1}e_je_j^\top(I+M)^{-1}}{1+te_j^\top(I+M)^{-1}e_j}=H-\frac{tHe_je_j^\top H}{1+tH_{jj}}\,.\] 
    Then, since $\|H_{\bullet j}\|_2\|H_{j\bullet}\|_2=\|H_{\bullet j}H_{j\bullet}\|_2=\|He_je_j^\top H\|_2$,  we have \[\|H_{\bullet j}\|_2\|H_{j\bullet}\|_2\leq \frac{1+tH_{jj}}{t}\left(\|H\|_2+\|(I+M+te_je_j^\top)^{-1}\|_2\right)\leq \frac{2(1+tH_{jj})}{t}\sqrt{1+\hk(M)}\,,\]
 which holds for all $t>0$.   By taking $t\rightarrow\infty$, we have
 \[\|H_{\bullet j}\|_2\|H_{j\bullet}\|_2\leq \inf\limits_{t>0}\frac{2\sqrt{1+\hk(M)}(1+tH_{jj})}{t}= 2H_{jj}\sqrt{1+\hk(M)}\,.
\]

\end{proof}
\begin{lem}\label{lem:double-H-Hadamard-product}
    Let $M\in\cPs$ and $H=(I+M)^{-1}$, then for any $v\in\R^n$,
     \[
        \|H\diag{ Hv} \|_2
    \leq 2\sqrt{n(1+\hk(M))}\|v\|_2\,.
    \] 
\end{lem}
\begin{proof}
    Since for all $i, j\in [n]$, $ [H\diag{Hv}]_{ij}=(Hv)_j H_{ij}$, using \Cref{lem:upper-bound-column-row-product} we have \begin{equation*}
        \|H\diag{Hv}\|_2^2\leq \sum\limits_{j=1}^n \|H_{\bullet j}\|_2^2|(Hv)_j|^2\leq \|v\|_2^2\sum\limits_{j=1}^n\|H_{\bullet j}\|_2^2\|H_{j\bullet }\|_2^2\leq 4(1+\hk(M))\|v\|_2^2\sum\limits_{j=1}^n H_{jj}^2\,.
    \end{equation*}
As a result, the bound for $\|H\diag{Hv}\|_2$ follows from $H_{jj}\in (0, 1]$ for all $j\in [n]$ by \Cref{prop:diagonal-entry-H} and then taking the square root. 
\end{proof}

\begin{proof}[Proof of \Cref{lem:two-step-progress}]

Let $M\in\cPs(\kappa)$, $(x,s,\mu)\in\cNc(\varepsilon)$, $\beta=\frac{1}{4096n^{3/2}}$ and $\delta=\frac{1}{1024n^{9/8}\sqrt{1+\kappa}}$. We decrease $\mu$ to $\mu'=\mu/(1+\delta)$ and take two Newton steps, obtaining $(x',s')$ and $(x'',s'')$. We will show that $(x'',s'',\mu')\in\cNc(\beta)$. The initial centrality error of $(x, s, \mu)\in\cNc(\varepsilon)$ is
$   \|Q_{H_{x,\mu}}\err(x,s,\mu)\|_2 \le \varepsilon$, for $\varepsilon=\frac{\beta}{2^{17}n(1+\kappa)^2}$. We keep the initial error in the calculation and bound its contribution at the end.

To simplify notation, let $H\coloneqq H_{x,\mu'}$, $H'\coloneqq H_{x',\mu'}$, and $H''\coloneqq H_{x'',\mu'}$. By \eqref{eq:def-centrality-error},
\begin{equation}\label{eq:two-step-initial-shift}
    \err(x,s,\mu')=(1+\delta)\err(x,s,\mu)-\delta\1\,.
\end{equation}
Since $q_{H_{x,\mu}}\ge\1$, we have $\|\err(x,s,\mu)\|_2\le\varepsilon$. Since $\delta\le1$ and $2\varepsilon\le\delta\sqrt n$, the triangle inequality gives:
\begin{equation}\label{eq:two-step-smallness}
    \|\err(x,s,\mu')\|_2
    \le\delta\sqrt n+2\|\err(x,s,\mu)\|_2
    \le2\delta\sqrt n\,.
\end{equation}

\paragraph{The first Newton step.}
After taking a single Newton step, the relative change in the iterate is $\Delta\bx\coloneqq(x'-x)/x=H\err(x,s,\mu')$. By \Cref{thm:equivalent_pstark_char} and \Cref{cor:normbounds}, respectively, we can bound
\begin{equation}\label{eq:two-step-compact-first}
\begin{aligned}
    \|\Delta\bx\|_2
    &\le\sqrt{1+\kappa}\|\err(x,s,\mu')\|_2,\\
    \|Q_H\Delta\bx\|_\infty
    &\le32\sqrt{1+\kappa}\|\err(x,s,\mu')\|_2.
\end{aligned}
\end{equation}
Additionally, using \Cref{lem:double-H-Hadamard-product} for $v=\err(x,s,\mu')$:
\begin{equation}\label{eq:two-step-first-matrix}
    \|H\diag{\Delta\bx}\|_2
    \le2\sqrt{n(1+\kappa)}\|\err(x,s,\mu')\|_2\,.
\end{equation}
The Newton step \eqref{eq: Newton-step for basic IPM} gives
\begin{equation}\label{eq:two-step-product-identity}
    \frac{x's'}{\mu'}=(\1_n+\Delta\bx)(\1_n-\Delta\bx)
    =\1_n-(\Delta\bx)^{2}\,.
\end{equation}
By \eqref{eq:two-step-compact-first} and then \eqref{eq:two-step-smallness}, $\|\Delta\bx\|_2\le2\delta\sqrt{n(1+\kappa)}\le1/512$, so $x'>0$. Together with~\eqref{eq:two-step-product-identity}, this gives $s'>0$. The relation $s'-s=M(x'-x)$ preserves linear feasibility, and the error after the first step is $\err(x',s',\mu')=(\Delta\bx)^{2}$.

\paragraph{The second Newton step.}
The second relative change is $\Delta\bx'\coloneqq\frac{x''-x'}{x'}=H'(\Delta\bx)^{2}$ and the resulting centrality error is $\err(x'',s'',\mu')=(\Delta\bx')^{2}$, respectively.
As in the proof of \Cref{lem:move-down2}, we bound the Rump norm of centrality error by a product of two norms:
\begin{equation}\label{eq:two-step-error-product}
    \|Q_{H'}\err(x'',s'',\mu')\|_2=\|Q_{H'} (\Delta\bx')^{2}\|_2
    \le\|Q_{H'}\Delta\bx'\|_\infty\|\Delta\bx'\|_2\,.
\end{equation}
We will bound these norms using \eqref{eq:two-step-compact-first} and \eqref{eq:two-step-first-matrix}. To do so, we first account for the changes from $H$ to $H'$ and from $q_H$ to $q_{H'}$.
Set $U\coloneqq\diag{\1+\Delta\bx}$. Using the definition of the Newton matrix and Woodbury's formula, we have
\begin{equation}\label{eq:two-step-matrix-update}
    H'=U^{-1}\bigl(U^{-2}-I_n+H^{-1}\bigr)^{-1}U^{-1}
    =U^{-1}\bigl(I_n+H(U^{-2}-I_n)\bigr)^{-1}HU^{-1}\,.
\end{equation}
Since $|\Delta\bx_i|\le1/4$, the identity
\[
    (U^{-2}-I_n)_{ii}
    =-\Delta\bx_i\frac{2+\Delta\bx_i}{(1+\Delta\bx_i)^2}
\]
and \eqref{eq:two-step-first-matrix} imply
\[
    \|H(U^{-2}-I_n)\|_2
    \le4\|H\diag{\Delta\bx}\|_2
    \le8\sqrt{n(1+\kappa)}\|\err(x,s,\mu')\|_2\le16\delta n\sqrt{1+\kappa}\leq \frac12\,,
\]
where we bound $\|\err(x, s, \mu')\|_2$ using \eqref{eq:two-step-smallness}. Hence $\|(I_n+H(U^{-2}-I_n))^{-1}\|_2\le2$. Also, $\|U^{-1}\|_2\le4/3$, and $U^{-1}$ commutes with $\diag{\Delta\bx}$. Substituting in \eqref{eq:two-step-matrix-update} gives
\begin{equation}\label{eq:two-step-matrix-transfer}
    \|H'\diag{\Delta\bx}\|_2
    \le4\|H\diag{\Delta\bx}\|_2
    \le8\sqrt{n(1+\kappa)}\|\err(x,s,\mu')\|_2\,.
\end{equation}
For the weights, using \Cref{lem:main-stability}, $|\log(1+a)|\le2|a|$ for $|a|\le1/2$, $\|\Delta\bar x\|_2\leq 2\delta\sqrt{n(1+\kappa)}$ and \eqref{eq:two-step-compact-first} give
\begin{equation}\label{eq:two-step-first-weights}
    \|\log q_{H'}-\log q_H\|_\infty
    \le32\Kp\|\Delta\bx\|_1
    \le 32\Kp n^{1/2} \|\Delta\bx\|_2\leq 64\delta Kn\sqrt{1+\kappa}\le\frac12\,,
    \end{equation}
and consequently, $q_{H'}/q_H\leq e^{1/2}\leq 2$.

We can now bound the two norms. From \eqref{eq:two-step-matrix-transfer} and \eqref{eq:two-step-compact-first},
\begin{equation}\label{eq:two-step-compact-second}
    \|\Delta\bx'\|_2=\|H'(\Delta\bar x)^2\|_2
    \le\|H'\diag{\Delta\bx}\|_2\|\Delta\bx\|_2
    \le8\sqrt n(1+\kappa)\|\err(x,s,\mu')\|_2^2\,.
\end{equation}
For the other norm, \Cref{lem:weights-well-defined} gives $\|Q_{H'}H'Q_{H'}^{-1}\|_{\max}\le\Kp$. Together with \eqref{eq:two-step-first-weights}, this yields
\begin{equation}\label{eq:two-step-second-infty}
\begin{aligned}
    \|Q_{H'}\Delta\bx'\|_\infty=\|Q_{H'}H'(\Delta\bar x)^2\|_\infty
    \le\Kp\|Q_{H'}(\Delta\bx)^{2}\|_1
    \le2\Kp\|Q_H\Delta\bx\|_\infty\|\Delta\bx\|_1
    \le512\sqrt n(1+\kappa)\|\err(x,s,\mu')\|_2^2,
\end{aligned}
\end{equation}
where the last inequality uses \eqref{eq:two-step-compact-first} and $\|\Delta\bx\|_1\le\sqrt n\|\Delta\bx\|_2$.

\paragraph{Bounding the centrality error.}
Combining \eqref{eq:two-step-compact-second} and \eqref{eq:two-step-second-infty} in \eqref{eq:two-step-error-product} gives
\begin{equation}\label{eq:two-step-compact-error}
    \|Q_{H'}\err(x'',s'',\mu')\|_2
    \le4096n(1+\kappa)^2\|\err(x,s,\mu')\|_2^4\,.
\end{equation}
To replace $Q_{H'}$ by the final weights $Q_{H''}$, note that \eqref{eq:two-step-compact-second} and \eqref{eq:two-step-smallness} give
\[
    \|\Delta\bx'\|_1\le\sqrt n\|\Delta\bx'\|_2
    \le8n(1+\kappa)\|\err(x,s,\mu')\|_2^2
    \le\frac{8}{512^2}\le\frac{1}{512}\,.
\]
Thus $x'',s''>0$ by applying the same identity as in \eqref{eq:two-step-product-identity} to the second step. The same application of \Cref{lem:main-stability} as in \eqref{eq:two-step-first-weights} gives $q_{H''}\le2q_{H'}$. Consequently,
\begin{equation}\label{eq:two-step-final-error}
\begin{aligned}
    \|Q_{H''}\err(x'',s'',\mu')\|_2
    &\le8192n(1+\kappa)^2\|\err(x,s,\mu')\|_2^4\\
    &\le8192\delta^4n^3(1+\kappa)^2
       +2^{16}n(1+\kappa)^2\|\err(x,s,\mu)\|_2\\
    &\le\frac{\beta}{2^{15}}+\frac{\beta}{2}\le\beta\,.
\end{aligned}
\end{equation}
For the second inequality, we apply to \eqref{eq:two-step-smallness} the bound $(a+b)^4\le a^4+4b$ for $0\le a\le a+b\le1$, with $a=\delta\sqrt n$ and $b=2\|\err(x,s,\mu)\|_2$; the condition $a+b\le 1$ follows from \eqref{eq:two-step-smallness}. The last line uses the choice of $\delta$ and $\|\err(x,s,\mu)\|_2\le\varepsilon$. This proves $(x'',s'',\mu')\in\cNc(\beta)$.
\end{proof}
 \Cref{lem:centralize} reduces the Rump norm error from $\beta=O(n^{-3/2})$ to $\varepsilon$ in $O(\log\log(n+\kappa))$ Newton steps. Thus, starting from an iterate in $\cNc(\beta)$, we can first repeatedly correct it to near centrality, and then apply two successive Newton steps towards $\mu'$. This leads to a constant-factor decrease in $\mu$ using $O(n^{9/8}\sqrt{1+\kappa})$ such Newton steps by \Cref{lem:two-step-progress}, and in total $ O(n^{9/8}\sqrt{1+\kappa}\log\log(n+\kappa))$ Newton steps including centering.

The extra $\log\log(n+\kappa)$ factor can be removed by showing that in fact it suffices to start from an $O(n^{-3/2})$ Rump neighborhood, which shows that the additional centering steps are not necessary. Thus we can provide a strengthening of \Cref{lem:two-step-progress}. Since the required changes in the proof are technical and do not contribute to the main ideas, we defer the proof of this strengthened version to Section~\ref{sec:two-step-no-recentering} in the appendix. \footnote{In fact, we believe that the proof can be further strengthened to show that after taking $O(\log \log n)$ Newton steps, the iterate is in the $O(n^{-3/2})$ Rump neighborhood, which would lead to an iteration bound of $\widetilde{O}(n \sqrt{ 1+ \kappa})$. However, we do not pursue this direction here, since it appears highly technical, and beyond the scope of this paper.}

\begin{restatable}{thm}{thmimprovedanalysis}\label{lem:two-step-progress-rump}
Let $M\in\cPs(\kappa)$, and set $\beta\in \left[\frac{1}{8192n^{3/2}}\,, \frac{1}{4096n^{3/2}}\right]$ and $\delta\in\left[\frac{1}{2436n^{9/8}\sqrt{1+\kappa}}\,,\frac{1}{1218n^{9/8}\sqrt{1+\kappa}}\right]$.
Suppose $(x,s,\mu)\in\cNc(\beta)$, and set $\mu'=\mu/(1+\delta)$. Two successive Newton steps toward $\mu'$, as in \eqref{eq: Newton-step for basic IPM}, produce two iterates $(x',s'),(x'',s'')\in\Feasp$ with $(x'',s'',\mu')\in\cNc(\beta)$.
\end{restatable}

\section{Hardness results}\label{sec:hardness}
\subsection{Hardness of testing membership between \texorpdfstring{$\cP$}{P} and \texorpdfstring{$\cP_0$}{P0}}
In this section, we prove \Cref{thm:hardness-deciding-sufficiency}. Our main tool is based on the characterization of $\cP$-matrix using sign-real spectral radius, first developed by Rump \cite{rump1997theorems}. For convenience, we restate the definition of the sign-real spectral radius of a matrix:   
\begin{defn}[Sign-real spectral radius restated]
 Let $\cS$ denote the set of signature matrices, i.e., diagonal matrices with entries $\pm 1$. For $M\in\R^{n\times n}$, its \textit{sign-real spectral radius} is defined as \[\rhos(M):=\max\limits_{S\in\cS}\left\{|\lambda|: SMx=\lambda x\,, \ \textnormal{for some} \ x\in\R^n_{\neq\0}\right\}\,,\]
with the convention that $\rhos(M):=0$ when $SM$ has no real eigenvalue for all $S\in\mathcal{S}^n$.
\end{defn}
 It follows directly from the definition that $\rhos(\alpha M)=|\alpha|\rhos(M)$ for all $\alpha\in\R$, which we shall exploit in the proof of \Cref{thm:hardness-deciding-sufficiency}. We refer to \cite{rump1997theorems} for further properties of $\rhos(M)$. To establish the connection between sign-real spectral radius and LCPs, we first define the Cayley transform below. Note that there are also different conventions and more general ways to define Cayley transform in the literature. 
\begin{defn}[Cayley transform] For $M\in\R^{n\times n}$ such that $I_n-M$ is nonsingular, the Cayley transform of $M$ is defined to be $\cC(M)\coloneqq(I_n-M)^{-1}(I_n+M)$. 
\end{defn}
Whenever a Cayley transform of $M$ exists, the inverse Cayley transform also exists if $I_n+M$ is nonsingular. The following proposition can be proven by direct verification. 
\begin{proposition}[Inverse Cayley transform]\label{prop: inverse Cayley transform}
    For $M\in\R^{n\times n}$ such that $I_n+M$ is nonsingular, the inverse Cayley transform defined by $\cC^{-1}(M)\coloneqq(M-I_n)(M+I_n)^{-1}$ satisfies $\cC(\cC^{-1}(M))=M$. If additionally $I_n-M$ is nonsingular so that $\cC(M)$ is defined, we have $\cC^{-1}(\cC(M))=M$. 
\end{proposition}

As a main result in \cite{rump1997theorems} that relates the sign-real spectral radius to the matrix classes considered in LCPs, Rump showed that $\rhos(M)<1$ if and only if its Cayley transform is a $\cP$-matrix.

\begin{thm}[Rump \cite{rump1997theorems}] \label{thm: Rump's theorem} For $M\in\R^{n\times n}$ such that $I_n-M$ is nonsingular,  $\rhos(M)<1$ if and only if $\cC(M)\in\cP$. 
\end{thm}
We refer to \cite{rump2003p} for a self-contained proof of this theorem. For a rational matrix $M\in\Q^{n\times n}$ whose Cayley transform is a $\cP$-matrix, a rational gap exists between $\rhos(M)$ and the upper bound $1$ asserted by \Cref{thm: Rump's theorem}, that is, there exists $R>0$ such that  $\rhos(M)<1$ if and only if $\rhos(M)\leq 1-1/R$. Moreover, such $R$ can be computed in terms of the entries of $M$.
\begin{restatable}{lem}{rationalgap}
\label{lem: rational gap lemma}
    For any $M\in\Q^{n\times n}$, let $R\coloneqq C^n(1+n\|M\|_{\max})^{n-1}$ where $C$ is the least positive integral common denominator of entries of $M$. Then $\rhos(M)<1$ if and only if $\rhos(M)\leq 1-\frac{1}{R}$, with $\log R=O(\poly(n, \length(M)))$.
\end{restatable}
The proof is similar to bounding the maximum absolute value of eigenvalues by the encoding length of the matrix. We defer it to \Cref{sec:deferred_hardness}. By \Cref{lem: matrix class relation}, any $\cP_0$-matrix is the limit of a sequence of $\cP$-matrices, i.e., $\cP_0=\cl(\cP)$. We can use a limiting argument and the continuity of $\rhos(.)$ to extend \Cref{thm: Rump's theorem} to $\cP_0$-matrices.

\begin{restatable}[Closed Rump's theorem]{thm}{closedRump}\label{thm: closed Rump's theorem}
  For $M\in\R^{n\times n}$ such that $I_n-M$ is nonsingular, $\rhos(M)\leq 1$ if and only if $\cC(M)\in\cP_0$.
\end{restatable}
The proof can be found in \Cref{sec:deferred_hardness}. Following the road-map described in \Cref{subsec: hardness overview}, we provide \Cref{alg: reduction} for the detailed construction of a matrix $\wt M$ for any given $M$ such that $M\in\cP\cap\Q^{n\times n}$ if and only if $\wt M\in\wt \cP\cap\Q^{n\times n}$ for any matrix class $\cP\subseteq\wt\cP\subseteq\cP_0$. We will prove its correctness, and therefore \Cref{thm:hardness-deciding-sufficiency}, using \Cref{thm: Rump's theorem}, \Cref{lem: rational gap lemma} and \Cref{thm: closed Rump's theorem}.  
\begin{algorithm}[ht]
     \caption{\textsc{Hardness reduction}}
    \label{alg: reduction}
    \SetKwInOut{Input}{Input}
    \SetKwInOut{Output}{Output}
    \SetKw{And}{\textbf{and}}
        \Input{ $M\in \Q^{n\times n}$, a matrix class $\wt \cP$ such that $\cP\subseteq\wt\cP\subseteq\cP_0$.}
    \Output{$\wt M\in\Q^{n\times n}$ such that $\wt M\in\wt \cP$ if and only if $ M\in\cP$.}
\If{$\det(I_n+M)\leq 0$}{$\wt M\gets -I_n$\;}\Else{$N\gets (M-I_n)(M+I_n)^{-1}$\;
Compute $R$ for $N$ as in \Cref{lem: rational gap lemma}\; 
$\alpha\gets 1+\frac{1}{2R}$\;
\If{$\det(I_n-\alpha N)= 0$}{$\wt M\gets -I_n$\;}\Else{$\wt M\gets (I_n-\alpha N)^{-1}(I_n+\alpha N)$\;}
}
     \Return{$\wt M$}
     \end{algorithm}

\begin{lem}\label{lem: hardness main theorem}
    Given $M\in\Q^{n\times n}$ and a matrix class $\wt\cP$ such that $\cP\subseteq \wt\cP\subseteq \cP_0$, Algorithm \ref{alg: reduction} outputs $\wt M\in\Q^{n\times n}$ such that $M\in\cP$ if and only if  $\wt M\in\wt\cP$, in polynomial time.
\end{lem}
\begin{proof}
   Using $\det(I_n+M)=\sum_{I\subseteq [n]}\det(M_{I, I})$,  if $\det(I_n+M)\leq 0$ then $M\notin\cP$. In this case, the algorithm returns $\wt M=-I_n\notin \wt\cP$. Otherwise, $N=(M-I_n)(M+I_n)^{-1}$ satisfies $\cC(N)=M$ by \Cref{prop: inverse Cayley transform}. $R$ is computed as in \Cref{lem: rational gap lemma} so that $\rhos(N)<1$ if and only if $\rhos(N)\leq 1-\frac{1}{R}$. Then as a result of \Cref{thm: Rump's theorem} and \Cref{lem: rational gap lemma}, $M\in \cP$ if and only if $\rhos(N)\leq 1-\frac{1}{R}$. Now we set $\alpha=1+\frac{1}{2R}$. \begin{claim}\label{claim: claim 1 for hardness proof}
   If $M\in\cP$ then $\rhos(\alpha N)<1$. If $M\notin \cP$ then $\rhos(\alpha N)>1$.    
   \end{claim} 
   \begin{proof}
       If $M\in\cP$, then $\rhos(N)\leq 1-\frac{1}{R}$, so $\rhos(\alpha N)=\alpha\rhos(N)\leq (1+\frac{1}{2R})(1-\frac{1}{R})<1$. If $M\notin\cP$, then $\rhos(N)\geq 1$, so $\rhos(\alpha N)=\alpha\rhos(N)\geq \alpha>1$. 
   \end{proof}
   \begin{claim}\label{claim: claim 2 for hardness proof}
       If $\det(I_n-\alpha N)= 0$ then $M\notin \cP$. 
   \end{claim}
   \begin{proof}
       Suppose that $M\in \cP$.  Since $\det(I_n-\alpha N)=0$, there exists some $x\neq \0_n$ such that $\alpha Nx=x$. By the definition of $\rhos(\alpha N)$, we get $\rhos(\alpha N)\geq 1$. This is a contradiction to \Cref{claim: claim 1 for hardness proof}.
   \end{proof}
In case  $\det(I_n-\alpha N)=0$ then  the algorithm returns $\wt M=-I_n\notin\wt\cP$. Otherwise, the algorithm returns $\wt M=(I_n-\alpha N)^{-1}(I_n+\alpha N)=\cC(\alpha N)$. We now show that $M\in\cP$ if and only if $\wt M\in\wt\cP$. Indeed, if $M\in\cP$, $\rhos(\alpha N)<1$ by \Cref{claim: claim 1 for hardness proof}, then $\wt M\in \cP$ by \Cref{thm: Rump's theorem}. Since $\cP\subseteq\wt\cP$, we get $\wt M\in\wt\cP$. If $M\not\in\cP$, $\rhos(\alpha N)>1$ by \Cref{claim: claim 1 for hardness proof}. Using \Cref{thm: closed Rump's theorem}, $\wt M\notin \cP_0$. Since $\wt\cP\subseteq\cP_0$, we get $\wt M\notin \wt\cP$.

Since $\log R=O(\poly(n, \length(M)))$ by \Cref{lem: rational gap lemma} and also $\length(N)=O(\poly(n, \length(M)))$, the sizes of $\alpha$ and hence $\wt M$ are polynomial in $n$ and $\length(M)$. 
\end{proof}
\begin{proof}[Proof of \Cref{thm:hardness-deciding-sufficiency}]
The first statement follows from \Cref{lem: hardness main theorem} and that deciding if $M\in\Q^{n\times n}$ is a $\cP$-matrix is coNP-complete \cite{coxson1994p}. By \Cref{thm: Valiaho theorem}, a $\cPs$-matrix is both column sufficient and row sufficient, and a certificate of being non-column sufficient or non-row sufficient is in NP \cite{tseng2000co}, thus the second statement follows.
\end{proof}

\subsection{Hardness of handicap number computation}
In this section, we prove \Cref{thm:main-computing-handicap-number}: for any graph $G=(V,E)$ we construct  a $\cP$-matrix whose handicap number determines the max cut of $G$.
By the \emph{incidence matrix} $A_G$ of $G$, we mean the directed incidence matrix for an arbitrary orientation of $G$. That is, a matrix $A_G\in \{0,\pm1\}^{|V|\times |E|}$, where for each edge $e=(i,j)\in E$, the $e$-th column of $A_G$ contains $+1$ for $i$ and $-1$ for $j$, or vice versa, and 0 in all other entries. This matrix is not unique, but the argument works for an arbitrary such choice.
Further, $L_G\coloneqq A_G A_G^\top \in \R^{n\times n}$ is the \emph{Laplacian matrix} of $G$; this is independent of the choice of the orientation. We recall the well-known connection between the Laplacian and Max-Cut, see e.g. \cite{goemans1995improved}. %
\begin{proposition}\label{prop:maxcut-laplacian}
For a graph $G=(V,E)$, 
\[
    \mc(G)=\frac{1}{4}\max_{v\in\{-1,1\}^n} v^\top L_G v\, .
\]
\end{proposition}

 When $M\in\cPs$ and we want to lower bound $x^\top Mx$ for $x\in\R^n$, the following change of variable via the midpoint of $x$ and $Mx$ would be useful, as it expresses $x^\top Mx$ as a difference of squares; the proof is deferred to \Cref{sec:deferred_hardness}.
\begin{restatable}{proposition}{changeofvariable}\label{prop:change-of-variable-handicap-number}
    Let $M\in\cPs$ and $x\in\R^n$, then the change of variable $y=\frac{1}{2}(x+Mx)$ is bijective and for $N=(M-I_n)/2$, we have \begin{equation*}
        x^\top Mx=\|y\|^2-\|N(I_n+N)^{-1} y\|^2\,. 
    \end{equation*}
   If $N^2=\mathbf{0}_{n\times n}$, then $x^\top Mx=y^\top(I_n-N^\top N)y$.
\end{restatable}
We derive the following corollary from \Cref{prop:maxcut-laplacian}; the proof is also given in \Cref{sec:deferred_hardness}.
\begin{restatable}{proposition}{maxcutidentity}\label{prop:identity-for-max-cut}
    Let $G=(V, E)$ be a connected  graph with $|V|=n\geq 2$ and $|E|=m$.  Then for all $v\in\R^n$, \[v^\top (L_G -I_n)v\leq \|v\|_\infty^2(4\mc(G)-n)\,.\]
\end{restatable}

We now introduce our construction of $M_G$.
\begin{lem}\label{lem: main lemma for NP hardness of handicap number}
  Let $G=(V, E)$ be a connected graph with $|V|=n\geq 2$ and $|E|=m$, and $A_G\in\R^{n\times m}$ be its incidence matrix. We define \begin{equation}
      M_G\coloneqq \left(\begin{array}{ccc}
        I_m   & 2A_G^\top & \0 \\
        \0   & I_n & 2\1_n \\
        \0 &  \0 & 1\\
      \end{array}\right)\in\mathbb{Z}^{(m+n+1)\times (m+n+1)}\,,
  \end{equation}
  then $M_G\in\cP$ and $\hk(M_G)=\mc(G)-1/4$.
\end{lem}
\begin{proof}
    $M_G\in\cP$ because all its principal minors are $1$. By \cite[Theorem 2.3]{valiaho1997determining}, $\hk( M_G)=\max_{z\in\R^{m+n+1}}\hk(M_G, z)$. Letting $y\coloneqq\frac{1}{2}(x+M_Gx)\in\R^{m+n+1}$, we get \begin{equation*}
        x=2(I_{m+n+1}+M_G)^{-1}y=\left(\begin{array}{ccc}
            I_m & -A^\top_G & \0 \\
           \0  & I_n & -\1_n \\
           \0 & \0 & 1
\end{array}\right)\left(\begin{array}{ccc}
             y_1  \\
             y_2 \\
             y_3
\end{array}\right)=\left(\begin{array}{cc}
             y_1-A_G^\top y_2  \\
             y_2-y_3\1_n \\
             y_3
        \end{array}\right)\,,\quad  y_1\in\R^m\,,\  y_2\in \R^{n}\,, \  y_3\in\R\,.
    \end{equation*}
    and \begin{equation}\label{eq: x circ Mx in hardness proof}
        x\circ M_Gx=\left(\begin{array}{cc}
             y_1-A_G^\top y_2  \\
             y_2-y_3\1_n \\
             y_3
        \end{array}\right)\circ \left(\begin{array}{ccc}
                 y_1+A_G^\top y_2  \\
                 y_2+y_3\1_n\\
                 y_3
            \end{array}\right)=\left(\begin{array}{ccc}
              y_1^{ 2}-(A_G^\top y_2)^{ 2}\\
                y_2^{ 2}-y_3^2\1_n\\
                  y_3^2
            \end{array}\right)\,.
    \end{equation}
 Since  $A_G$ is the incidence matrix of $G$, $A_G^\top\1_n=\0$. For $N=(M_G-I_{m+n+1})/2$, we have  \[N^2=\left(\begin{array}{ccc}
        \0   & A_G^\top & \0 \\
        \0   & \0 & \1_n \\
        \0 &  \0 & 0\\
      \end{array}\right)^2= \left(\begin{array}{ccc}
        \0   & \0 & A_G^\top\1_n \\
        \0   & \0 & \0 \\
        \0 &  \0 & 0\\
      \end{array}\right)=\0\,, \quad I_{m+n+1}-N^\top N=\left(\begin{array}{ccc}
        I_m   & \0 & \0 \\
        \0   & I_n-L_G & \0 \\
        \0 &  \0 & 1-n\\
      \end{array}\right)\,. \]
Using \Cref{prop:change-of-variable-handicap-number}, we have \begin{equation}\label{eq: x^top Mx in y}
\begin{aligned}
    x^\top M_G x&=y^\top(I_{m+n+1}-N^\top N)y\\
    &=\|y\|^2-y_2^\top L_G^\top y_2-ny_3^2\\
    &=\|y_1\|^2+y_2^\top(I_n- L_G)y_2-(n-1)y_3^2\,.
\end{aligned}
\end{equation}
\begin{claim}
    $\hk(M_G)\geq \mc(G)-1/4$.
\end{claim}
\begin{proof}
Let $y_2\in\{-1, 1\}^{n}$ be the vector attaining $\mc(G)$, that is, $y_2^\top L_G y_2=4\mc(G)$. Also, set $y_1=\0_m$ and $y_3=1$. Then in \eqref{eq: x circ Mx in hardness proof}, we have \begin{equation*}
    x\circ  M_G x=\left(\begin{array}{ccc}
         -(A_G^\top y_2)^2 \\
         \0_n \\
         1
    \end{array}\right)\,,
\end{equation*} giving $x^\top M_Gx=-y_2^\top L_G y_2+1=-4\mc(G)+1$ and $\sum_{i\in\mathcal{I}^+_{ M_G}(x)}x_i( M_G x)_i=1$. Thus, \[x^\top M_Gx=-(4\mc(G)-1) \sum_{i\in\mathcal{I}^+_{ M_G}(x)}x_i( M_G x)_i\,,\] so for $x=(-A_G^\top y_2, y_2-\1_n, 1)\in\R^{m+n+1}$, $\hk( M_G, x)=\mc(G)-1/4$. Hence, $\hk( M_G)\geq \mc(G)-1/4$.  
\end{proof}
\begin{claim}
    $\hk(M_G)\leq \mc(G)-1/4$.
\end{claim} 
\begin{proof}
From \eqref{eq: x^top Mx in y} and \Cref{prop:identity-for-max-cut}, we get
\begin{equation}\label{eq: x^top Mx in hardness proof}
  \begin{aligned}
        x^\top M_G x \geq -(4\mc(G)-n)\|y_2\|_\infty^2-(n-1)y_3^2\,.
    \end{aligned}
    \end{equation}
We now consider the positive contribution in $x^\top M_G x$. By looking at the possible positive contribution from only the middle and last components of the vector $x\circ M_G x$ in \eqref{eq: x circ Mx in hardness proof}, \begin{equation}\label{eq: upper bound on handicap number (2)}
\sum\limits_{i\in\mathcal{I}^+_{M_G}(x)}x_i(M_Gx)_i\geq  y_3^2 +\sum\limits_{i=1}^{n}\max \{{y_{2}}_i^2-y_3^2\,, 0\} \geq \max\{ \|y_2\|_\infty^2\,, y_3^2\}
\end{equation}
Since $\mc(G)\geq m/2$ (see e.g. \cite{goemans1995improved}) and $G$ is connected, we have $\mc(G)\geq m/2\geq (n-1)/2\geq n/4$. Hence, we can combine \eqref{eq: x^top Mx in hardness proof} and \eqref{eq: upper bound on handicap number (2)} to obtain \begin{equation}\label{eq: upper bound for handicap number of x}
    x^\top M_G x\geq -(4\mc(G)-1) \sum\limits_{i\in\mathcal{I}^+_{ M_G}(x)}x_i( M_Gx)_i\,.
\end{equation}
This shows that $\hk(M_G, x)\leq \mc(G)-1/4$ for all $x\in\R^{m+n+1}$, and hence $\hk(M_G)\leq \mc(G)-1/4$.
\end{proof}
These claims together show $\hk(M_G)=\mc(G)-1/4$. 
\end{proof}
\begin{proof}[Proof of \Cref{thm:main-computing-handicap-number}] As in \Cref{lem: main lemma for NP hardness of handicap number}, we can construct $M_G\in\mathbb{Z}^{(m+n+1)\times (m+n+1)}$ for any  connected graph $G$ with $|V|=n\geq 2$ and $|E|=m$.   We have $M_G\in\cP$ and $\hk(M_G)=\mc(G)-1/4$. Thus, a polynomial algorithm that can approximate $1/4+\hk(M)$ for a sufficient matrix $M$ within a multiplicative factor of $\alpha$ yields the same approximation guarantee for the max cut of $G$ in polynomial time, proving the theorem. 
\end{proof}
\subsection{Rump's characterization of sufficient matrix}
Having seen Rump's original characterization of $\cP$-matrices using sign-real spectral radius \Cref{thm: Rump's theorem} and also its extension for $\cP_0$-matrices \Cref{thm: closed Rump's theorem}, we conclude this section by giving a similar characterization for $\cPs$-matrices.

Throughout this section, we assume $M\in\cP_0$ so that $I_n+M$ is nonsingular (\Cref{cor:invertiblity}).
We again make use of the change of variable by
midpoint in \Cref{prop:change-of-variable-handicap-number} to express $v\circ (Mv)$ as a difference of squares. Let $x\coloneqq\frac{1}{2}(v+Mv)$ and $N\coloneqq\frac{1}{2}(M-I_n)$. As in the proof of \Cref{prop:change-of-variable-handicap-number}, 
\begin{equation}\label{eq:change-of-variable-v-u}
Mv=\left(I_n+N(I_n+N)^{-1}\right)x\,,\quad v=\left(I_n-N(I_n+N)^{-1}\right)x\,,
\end{equation}
where $N(I_n+N)^{-1}=(M-I_n)(M+I_n)^{-1}=\cC^{-1}(M)$ is the inverse Cayley transform of $M$ (see \Cref{prop: inverse Cayley transform}). Hence, we have \begin{equation}\label{eq:difference-of-squre-F_M}
    v\circ (Mv)=\left(I-\cC^{-1}(M)\right)x\circ \left(I+\cC^{-1}(M)\right)x=x^2-(\cC^{-1}(M)x)^2\,.
\end{equation}
We recall $M$ is column sufficient if $v\circ (Mv)\leq \0_n$ implies $v\circ (Mv)=\0_n$ for all $v\in\R^n$, and is row sufficient if $M^\top$ is column sufficient, and further, $M\in\cPs$ if and only if $M$ is both column and row sufficient by V\"aliaho's theorem (\Cref{thm: Valiaho theorem}). The recent result by E.-Nagy and V\'egh \cite[Theorem 1.4]{nagy2026handicap} also gives a strengthening of V\"aliaho's theorem, characterizing sufficient matrices using both $v\circ Mv$ and $u\circ Mv+v\circ Mu$.

\begin{thm}[\cite{nagy2026handicap}]\label{thm: sufficient theorem} For any $M\in \R^{n\times n}$, $M$ is sufficient if and only if for all $u, v\in\R^n$,  $v\circ Mv\leq \0_n$ and $u\circ Mv+v\circ Mu\leq \0_n$ imply $v\circ Mv=\0_n$ and $u\circ Mv+v\circ Mu=\0_n$.  
\end{thm}

  Similar to \eqref{eq:difference-of-squre-F_M}, we can express $u\circ Mv+v\circ Mu$ as a difference of squares. We let $x=\frac{1}{2}(v+Mv)$, $N=\frac{1}{2}(M-I_n)$ as before and also $y=\frac{1}{2}(u+Mu)$. Then recalling \eqref{eq:change-of-variable-v-u} and $\cC^{-1}(M)=N(I_n+N)^{-1}$, 
\begin{equation}\label{eq:difference-of-square-derivative}
    \begin{aligned}
        u\circ Mv+v\circ Mu&=\left(I_n-\cC^{-1}(M)\right)y\circ \left(I_n+\cC^{-1}(M)\right)x+\left(I_n-\cC^{-1}(M)\right)x\circ \left(I_n+\cC^{-1}(M)\right)y\\
        &=2\left(x\circ y-\cC^{-1}(M)x\circ \cC^{-1}(M)y\right)\,.
    \end{aligned}
\end{equation}
We can now give a full characterization of $\cP$-, $\cPs$- and $\cP_0$-matrices using sign-real spectral radius.  
\begin{thm}[Rump's theorem for $\cPs$-matrix]\label{thm:sufficient-Rump-theorem}
    For $M\in\R^{n\times n}$, if $I_n+M$ is singular, then $M\notin\cP_0$. Suppose that $I_n+M$ is nonsingular, and let $\wt M\coloneqq\cC^{-1}(M)=(M-I_n)(M+I_n)^{-1}$. Then, \begin{enumerate}[label=(\roman*)]
        \item $M\in\cP$ if and only if $\rho_0^S(\wt M)<1$.
        \item $M\in\cPs\setminus \cP$ if and only if it satisfies the following: $\rho_0^S(\wt M)=1$; for all $x\in\R^n$, if $|\wt Mx|\geq |x|$, then $|\wt Mx|=|x|$; for all $x, y\in\R^n$ if $|\wt Mx|=|x|$ and $\wt Mx\circ \wt My\geq x\circ y$ then $\wt Mx\circ \wt M y=x\circ y$. 
        \item $M\in\cP_0$ if and only if $\rho_0^S(\wt M)\leq 1$. 
    \end{enumerate}
\end{thm}
\begin{proof}
 By \Cref{cor:invertiblity}, $M\in\cP_0$ implies $I_n+M$ is nonsingular. Part (i) and (iii) follow from \Cref{thm: Rump's theorem} and \Cref{thm: closed Rump's theorem}. Using \Cref{thm: sufficient theorem}, \eqref{eq:difference-of-squre-F_M}, and \eqref{eq:difference-of-square-derivative}, $M\in\cPs\setminus\cP$ if and only if $\rhos(\wt M)=1$, and for all $x, y\in\R^n$, if $|\wt Mx|\geq |x|$ and $\wt M x\circ \wt My\geq x\circ y$, then $|\wt Mx|=|x|$ and $\wt M x\circ \wt My=x\circ y$. This is equivalent to Part (ii).  
\end{proof}
Criteria for a general matrix to be sufficient were extensively studied by V\"aliaho \cite{valiaho1996criteria}.
 Apart from providing a spectral viewpoint of $\cPs$-matrices which could be of theoretical interest, \Cref{thm:sufficient-Rump-theorem} may also be used to develop new properties of sufficient matrices, which then give practical screening tests to detect non-sufficiency.

For example, \Cref{thm:sufficient-Rump-theorem} yields the following corollary on the structure of the zero eigenspace of a sufficient matrix.
An eigenvalue is \textit{semi-simple} if its algebraic multiplicity (possibly zero) is equal to its geometric multiplicity, that is, the associated Jordan block of that eigenvalue is diagonal. For $x\in\R^n$, we define $\sgn(x)\in\{\pm1, 0\}^n$ to be $(\sgn(x))_i=\sgn(x_i)$ for all $i\in [n]$, where $\sgn(0):=0$. For $U\subseteq \R^n$, we let $\sgn(U):=\{\sgn(x): x\in U\}$.
\begin{restatable}{cor}{zeroeigencor}
    Let $M\in\cPs\cap \R^{n\times n}$. Then $0$ is a semi-simple eigenvalue of $M$, in particular, $\ker(M^p)=\ker(M)$ and $\rk(M^p)=\rk(M)$ for all $p\in\mathbb{N}$. Further, $\sgn(\ker(M))=\sgn(\ker(M^\top))$. 
\end{restatable}
The proof can be found in \Cref{sec:deferred_hardness}.

\section{Initialization, termination and guessing of the handicap number}\label{sec: initialisation and termination}

We start by two simple statements on (column) sufficient matrices that will be needed for the initialization scheme. The first one is true for any sufficient matrix, and follows easily by the convexity of the solution set of \ref{eq: LCP}, see e.g. \cite[Theorem 3.5.8]{cottle2009linear}.
\begin{lem}\label{lem:col-suff}
    Let $M\in\cPs\cap\R^{n\times n}$, $q\in\R^n$, and let $(x,s)$ and $(x',s')$ be two solutions to \ref{eq: LCP}. Then, for any $i\in [n]$, either $x_i=x'_i$ or $s_i=s'_i$.
\end{lem}
For sufficient matrices, the above statement is also implied by the next bound, that can be seen as a generalization of the near monotonicity of the central path of linear programs, see e.g., \cite[Lemma 16]{vavasis1996primal}. We will use this to identify the optimal support from a sufficiently good approximate solution.
\begin{lem}\label{lem:near-mon}
 Let $M\in\cPs(\kappa)$, $q\in\R^n$, and let $(x,s,\mu)\in\cN^2(1/2)$, and let $(\bar x,\bar s)$ be any solution to \ref{eq: LCP}. Then,
 \[
 (\bar x,\bar s)\le 32(\kappa+1)n\cdot (x,s)\, .
 \]
\end{lem}
\begin{proof}
 Let $v\coloneqq x-\bar x$.
Then,
\[
v (Mv)= 
(x-\bar x)(s-\bar s)=x s+\bar x \bar s-\bar x s-x \bar s\le \frac{3}{2}\mu\cdot \1_n -\bar x s-\bar s x\le \frac{3}{2}\mu\cdot \1_n \, , 
\]
using the bound on $xs\le \frac{3}{2}\mu\cdot \1_n$ from $(x,s,\mu)\in\cN^2(1/2)$, and  $\bar x \bar s=0$. By the definition of $\kappa$, 
\[
0\le (4\kappa+1) \sum_{i\in \posin_M(v)}v_i (M v)_i +  \sum_{i\in \negin_M(v)}v_i (M v)_i\le \frac{3}{2}\cdot (4\kappa+1)|\posin_M(v)|\mu + \sum_{i\in \negin_M(v)}v_i (M v)_i\, ,
\]
implying that for any $i\in [n]$
\[
-\frac{3}{2}\cdot (4\kappa+1)(n-1)\mu  \le v_i(Mv)_i \le \frac{3}{2}\mu\cdot \1_n -\bar x_i s_i-\bar s_i x_i\, ,
\]
giving
\[
\bar x_i s_i+\bar s_i x_i\le 3(4\kappa+1)n \mu\, .
\]
By $(x,s,\mu)\in\cN^2(1/2)$, $x_is_i\geq \mu/2$, giving
\[
\frac{\bar x_i}{x_i}+\frac{\bar s_i}{s_i}\le 6(4\kappa+1) n\, ,
\]
implying the claim.
\end{proof}

\subsection{The extended system}
To find an initial point $(x^0, s^0)$ for \Cref{alg:rescale-ipm-simple}, we follow the same procedure as in Section 5 of \cite{kojima1991unified}, by considering an  $\LCP(\what M, \hat q)$ of size $2n$ that embeds \ref{eq: LCP}. We construct an embedding along with an initial point that lies \emph{exactly} on the central path. Solving the embedded problem optimally (or to sufficiently high accuracy) either gives a solution to \ref{eq: LCP}, or we may conclude that no feasible solution to the original problem exists.

First, we recall standard bit-complexity bounds, see e.g. \cite{SchrijverLIP}. %
\begin{lem}\label{lem:bit-complexity}
    $M\in\Q^{n\times n}$, $q\in\Q^n$ and $L$ be the total encoding length of $(M, q)$. Then, there exists $\Delta\le 2^{O(L)}$, $\Delta\in\mathbb{Z}$ such that  $\|M\|_{\max}+\|q\|_{\infty}\le \Delta$, and 
    for every basic feasible solution $(x,s)$ to the polyhedron $\Feas(M,q)$, and every $i\in [n]$, $x_i,s_i\in\{0\}\cup [1/\Delta,\Delta]$.
\end{lem}
We say that $(x,s)$ is a complementary basic feasible solution  if it is a complementary solution to \ref{eq: LCP} as well as a basic feasible solution to $\Feas(M,q)$. Clearly, for any complementary solution $(x,s)$ to \ref{eq: LCP}, there exists a complementary basic feasible solution $(\bar x,\bar s)$ with $\supp(\bar x)\subseteq \supp(x)$ and $\supp(\bar s)\subseteq \supp(s)$. 

The embedding is given as
 \begin{equation}\label{eq: extended LCP}\tag{$\LCP(\what M, \hat q)$}
   \what M \hat x+\hat q=\hat s, \quad \hat x^\top \hat s=0,\quad \hat x, \hat s\geq\0_{2n}\, , 
\end{equation}
where \begin{equation*}
   \what M=\left(\begin{array}{cc}
       M  & I_n \\
      -I_n   & \0
    \end{array}\right)\,,\quad \what q=\left(\begin{array}{cc}
         q  \\
         \tilde q 
    \end{array}\right)\,,\quad \hat x=\left(\begin{array}{cc}
         x  \\
         \tilde x 
    \end{array}\right)\,,\quad \hat s=\left(\begin{array}{cc}
         s  \\
         \tilde s 
    \end{array}\right)\,.
\end{equation*}
\begin{lem}[Initialization]\label{lem: initial point}
    Let $M\in\cPs\cap\Q^{n\times n}$, $q\in\Q^n$ and $L$ be the total encoding length of $(M, q)$, and let $\Delta\le 2^{O(L)}$ as in \Cref{lem:bit-complexity}.  Let us define $\hat x^0=(x^0, \tilde x^0)$, $\hat s^0=(s^0, \tilde s^0)$, $\tilde q$ and $\hat \mu_0$ by   \begin{equation}\label{eq: initial point formula}
    \begin{aligned}
        x^0=2\Delta\1_n\,,\quad s^0=4\Delta^2 n \1_n\,, \quad \tilde x^0=s^0-Mx^0-q\,,\quad \tilde s^0=\frac{8n\Delta^3}{\tilde x^0}\,,\quad \tilde q=x^0+\tilde s^0\,, \quad \hat\mu_0=8n\Delta^3\,.
  \end{aligned}
    \end{equation} 
  $M\in\cPs(\kappa)$ if and only if $\what M\in\cPs(\kappa)$ and $\hk(M)=\hk(\what M)$, and $(\hat x^0, \hat s^0)$ is a central path point to \ref{eq: extended LCP} with $\hat\mu_0$.  Let $(\hat x^*, \hat s^*)=((x^*, \tilde x^*), (s^*, \tilde s^*))$ be any solution to \ref{eq: extended LCP}. If $\tilde x^*=\0_n$, then $(x^*, s^*)$ is a solution to \ref{eq: LCP}. If $\tilde x^*\neq\0_n$, then \ref{eq: LCP} has no solution. 
\end{lem}
\begin{proof}
    The proof of $M\in\cPs(\kappa)$ if and only if $\what M\in\cPs(\kappa)$ can be found in \cite[Lemma 5.3]{kojima1991unified}, implying  $\hk(M)=\hk(\what M)$. It can be directly verified that $\hat s^0=\what M\hat x^0+\what q$ and that $\hat x^0\hat s^0=8n\Delta^3\1_{2n}$.  Further,
    \[\tilde x^0\geq s^0-\|Mx^0+q\|_\infty\1_n\geq (4\Delta^2 n-2\Delta^2 n -\Delta) \1_n>\0_n\,, 
    \]
   Hence, $(\hat x^0,\hat s^0)$ lies on the central path with $\hat\mu_0=8n\Delta^3$. If $\tilde x^*=\0_n$ then $(x^*, s^*)$ is a solution to \ref{eq: LCP} by the construction of \ref{eq: extended LCP}.  
   
   If $\tilde x^*\neq\0_n$, then let $i\in [n]$ with $\tilde x^*_i>0$. For a contradiction, assume \ref{eq: LCP} is feasible, and select $(\bar x,\bar s)$ as any complementary basic feasible solution of \ref{eq: LCP}. This extends to a solution $((\bar x,\0_n),(\bar s,\tilde s))$ of \ref{eq: extended LCP} by setting  $\tilde s\coloneqq \tilde q-\bar x$.  Since $\tilde x_i=0<\tilde x^*_i$, \Cref{lem:col-suff} implies that $\tilde s_i=\tilde s_i^*=0$, where the latter equality follows by $\tilde x_i^*\tilde s_i^*=0$.
But 
   \begin{equation}\label{eq:s-big}
   \tilde s_i=\tilde q_i-\bar x_i\ge x^0_i -\bar x_i=2\Delta-\bar x_i>0\, ,
   \end{equation}
   by \Cref{lem:bit-complexity}, a contradiction.
\end{proof}
We now use the near-monotonicity property \Cref{lem:near-mon} to identify the optimal support of \ref{eq: LCP} from an approximate solution to the embedded system, assuming that $M\in\cPs(\kappa)$.
\begin{lem}\label{lem:terminate-kappa}
  Let $M\in\Q^{n\times n}$, $q\in\Q^n$, $\kappa>0$, and $L$ be the total encoding length of $(M, q)$, and let $\Delta\le 2^{O(L)}$ as in \Cref{lem:bit-complexity}.  Consider the embedded system \ref{eq: extended LCP} with $\tilde q$ defined as in \eqref{eq: initial point formula}, and let $((x,\tilde x), (s,\tilde s), \mu)\in\cN^2(1/2)$  and $\mu\le 1/(16384 n^3\Delta^2(\kappa+1)^2)$. Define
\[I\coloneqq\{i\in [n]: x_i\leq\sqrt{2n\mu}\}\,,\quad J\coloneqq\{j\in[n]: s_j\leq\sqrt{2 n\mu}\}\,.\] 
Then $I\cup J=[n]$, and  at least one of the following hold:
\begin{enumerate}[label=(\roman*)]
\item There exists a solution $(\bar x,\bar s)\in\Feas(M,q)$ with  $\bar x_I=\0_I$ and $\bar s_J=\0_J$, which is a complementary solution to \ref{eq: LCP};
\item \ref{eq: LCP} is infeasible ;
\item $\hk(M)>\kappa$.
\end{enumerate}
\end{lem}
\begin{proof}
Since $((x, \tilde x), (s, \tilde s), \mu)\in\cN^2(1/2)$, by looking at the first $n$ coordinates,
 $|x^\top s/\mu -n|\le  \|\err(x,s,\mu)\|_1\le \sqrt{n}/2$, and consequently $x^\top s\le 2n\mu$. Hence, for each $i\in [n]$, $0\le x_is_i\le 2n\mu$ implies that $I\cup J=[n]$.
 Assume that $\hk(M)\le \kappa$, that is, $M\in\cPs(\kappa)$, and that \ref{eq: LCP} is feasible. We show that {\em (i)} must hold.

  Let $(\bar x,\bar s)$ be a complementary basic feasible solution to the original system  \ref{eq: LCP}. Thus, $((\bar x,\0_n),(\bar s,\tilde q-\bar x))$ is complementary feasible to \ref{eq: extended LCP}; for this, $\tilde q-\bar x>\0_n$ follows as in \eqref{eq:s-big}. 
  We apply \Cref{lem:near-mon} to $(( x,\tilde x),( s,\tilde s))$  and to  $((\bar x,\0_n),(\bar s,\tilde q-\bar x))$; we use that  $\hk(\what M)=\hk(M)\le \kappa$ by \Cref{lem: initial point}. Hence,
  \[
  \bar x_i \le 64(\kappa+1)nx_i\leq \sqrt{2n\mu}\cdot  64(\kappa+1)n<\frac{1}{\Delta}\, ,\quad \forall i\in I\, , \quad\mbox{and}\quad  \bar s_i \leq 64(\kappa+1)ns_i\le \sqrt{2n\mu}\cdot  64(\kappa+1)n<\frac{1}{\Delta}\, ,\quad \forall i\in J\, .
  \]
  By \Cref{lem:bit-complexity}, $\bar x_I=\0_I$ and $\bar s_J=\0_J$.
  Since  \ref{eq: LCP} is assumed to be feasible and this holds for the chosen complementary basic feasible solution, we can conclude {\em (i)}.
\end{proof}

\subsection{The overall algorithm}
\Cref{alg:rescale-ipm-full} presents the full version of our algorithm needed to prove \Cref{thm:main}, where we do not assume knowledge of $\hk(M)$.
As shown in \Cref{sec:hardness}, given an instance of \ref{eq: LCP}, it is coNP-hard to decide if $M$ is sufficient or not, and even if there is a guarantee that $M\in\cPs$, it is still NP-hard to compute its handicap number. By using similar strategies developed by Potra and Liu \cite{potra2005predictor} to handle unknown handicap number and Ill\'es, E.-Nagy and  Terlaky \cite{illes2010polynomial} to adapt the algorithm to treat general matrices, we can modify our algorithm to require no knowledge of $\hk(M)$ or whether $M\in\cPs$ or not. We make use of the following upper bound on the handicap number of any rational  $\cPs$-matrix in terms of its encoding length.

 \begin{thm}[\cite{nagy2026handicap}]\label{thm:encoding-length}
    For any $M\in\cPs\cap\Q^{n\times n}$ with encoding length $\length(M)$, $\hk(M)< 2^{O(n\length(M))}$.
\end{thm}

In order to get the better of the running times in \Cref{thm:Newton-simple} and \Cref{lem:two-step-progress-rump}, for a guess $\kappa\ge 1$ on $\hk(M)$, we define the parameters $\beta_n(\kappa)$ and $\delta_n(\kappa)$ as 
\begin{equation}\label{eq:greek}
\beta_n(\kappa)\in\begin{cases}
 \left[\frac{1}{8(1+\kappa)}, \frac{1}{4(1+\kappa)}\right]\, ,\,\quad &\mbox{if }\kappa^4 \le n^5\, \\
 \left[\frac{1}{8192n^{3/2}}\,, \frac{1}{4096n^{3/2}}\right]\, ,\,\quad &\mbox{if }\kappa^4 >n^5\, 
\end{cases}\, \quad 
\delta_n(\kappa)\in\begin{cases}
\left[\frac{1}{32\sqrt{n}(1+\kappa)},\frac{1}{16\sqrt{n}(1+\kappa)}\right]\, ,\,\quad &\mbox{if }\kappa^4 \le n^5\, \\
\left[\frac{1}{2436n^{9/8}\sqrt{1+\kappa}}\,,\frac{1}{1218n^{9/8}\sqrt{1+\kappa}}\right]\,\quad &\mbox{if }\kappa^4 >n^5\, 
\end{cases}
\end{equation}

We separate the iterations for a fixed promise $\kappa$ from the guessing procedure. The subroutine \textsc{FixedKappaNewton} in \Cref{alg:rescale-ipm-full} takes an input $(M,q)$ and an initial point $(x^0,s^0,\mu_0)$ on the central path, along with a fixed guess $\kappa\geq1$ on $\hk(M)$. It stays in the neighborhood $\beta_n(\kappa)$ and uses step size $\delta_n(\kappa)$. There are three possible outputs:
\begin{enumerate}[label=(\roman*)]
\item $(\hat x,\hat s,\mu)\in\cN^2(\beta_n(\kappa))$ with $\mu\leq\targetmu$, or
\item $\hk(M)>\kappa$, or
\item $M\notin\cPs$.
\end{enumerate}

\Cref{alg:rescale-ipm-full} takes a rational input $(M,q)$. According to \Cref{thm:feasibile}, if $M\in\cPs$ and
$\Feas(M,q)\neq \emptyset$, then \ref{eq: LCP} must be feasible. Hence, we can first run an LP solver to solve the linear relaxation $Mx+q=s$, $x,s\ge\0_n$. If this is infeasible, we terminate. Otherwise, either the LCP is feasible or $M\notin\cPs$ must hold.

 We let $L=\length(M,q)$, and set $\Delta=2^{O(L)}$ as in \Cref{lem:bit-complexity}.
We start with the guess $\kappa=1$; this will be doubled after each unsuccessful attempt. For each value of $\kappa$, we call \textsc{FixedKappaNewton} on the 
 embedded system 
 $(\what M,\what q)$ as in \ref{eq: extended LCP} with the  initial point $(\hat x^0,\hat s^0,\hat\mu_0)$ on the central path as in \Cref{lem: initial point}, and $\targetmu$ as in \Cref{lem:terminate-kappa}.

 If it returns $\what M\notin\cPs$, we terminate. If the output is $\hk(\what M)>\kappa$, then we double $\kappa$ and restart. If it returns a solution $((x,\tilde x), (s,\tilde s), \mu)\in\cN^2(\beta_{2n}(\kappa))$, we identify $(I,J)$ as in  \Cref{lem:terminate-kappa}, and solve the LP $M\bar x+q=\bar s$, $\bar x,\bar s\ge \0_n$, $\bar x_I=\0_I$, $\bar s_J=\0_J$. This is obtained from an LP solver denoted as \textsc{LPSolve}. 
 If this is feasible, we return $(\bar x,\bar s)$. Otherwise, we conclude again $\hk(M)>\kappa$ and restart after doubling the guess. If $\kappa$ exceeds the bound $T$ as in \Cref{thm:encoding-length}, then we terminate by concluding again $M\notin\cPs$.

\begin{algorithm}[htb!]
    \caption{\textsc{Full version Newton's Method for LCP path following}}
    \label{alg:rescale-ipm-full}
    \SetKwInOut{Input}{Input}
    \SetKwInOut{Output}{Output}
      \SetKwFunction{LPSolve}{LPSolve}
    \SetKwFunction{FixedKappaNewton}{FixedKappaNewton}
    \SetKwProg{Fn}{Function}{:}{}
    \Input{$M\in\Q^{n\times n}$ and $q\in\Q^n$ }
    \Output{Complementary solution $(x,s)$ to \ref{eq: LCP},  or the conclusion \ref{eq: LCP} is infeasible, or that $M\notin\cPs$}
    \Fn{\FixedKappaNewton{$M,q,x^0,s^0,\mu_0,\targetmu,\kappa$}}{
        $(x,s,\mu)\gets (x^0,s^0,\mu_0)$\;
     $\beta\gets\beta_n(\kappa)$,     $\delta\gets\delta_n(\kappa)$  as in \eqref{eq:greek}\;
        \While{$\mu>\targetmu$}{
            Compute Newton step $(\Delta x,\Delta s)$ for $(x,s,\frac{\mu}{1+\delta})$ as in \eqref{eq: Newton-step for basic IPM}\label{line:compute-Newton-1}\;
            \lIf{\textnormal{\eqref{eq: Newton-step for basic IPM} has no solution or non-unique solution}}{\Return{$M\notin\cPs$}\label{line:fail-not-in-P_0-1}}
            $(x',s')\gets (x+\Delta x,s+\Delta s)$\;
            \lIf{$(x',s')\notin\Feasp$}{\Return{$\hk(M)>\kappa$}\label{line:feasibility-test}}
            Compute Newton step $(\Delta x,\Delta s)$ for $(x',s',\frac{\mu}{1+\delta})$ as in \eqref{eq: Newton-step for basic IPM}\label{line:compute-Newton-2}\;
            \lIf{\textnormal{\eqref{eq: Newton-step for basic IPM} has no solution or non-unique solution}}{\Return{\textnormal{$M\notin\cPs$}}\label{line:fail-not-in-P_0-2}}
            $(x_g,s_g,\mu_g)\gets (x'+\Delta x,s'+\Delta s,\frac{\mu}{1+\delta})$\label{line:test-iterate-update}\;
            \lIf{$(x_g,s_g,\mu_g)\notin\cN^2(\beta)$}{\Return{$\hk(M)>\kappa$}\label{line:neighbourhood-test}}
            $(x,s,\mu)\gets (x_g,s_g,\mu_g)$\;
        }
        \Return{$(x,s,\mu)$}\;
    }
    \BlankLine
   $U\gets$ \LPSolve{$Mx+q=s\,,\, x,s\ge \0_n$}\;
        \lIf{$U=\textrm{Infeasible}$}{\Return{Infeasible}}
  $\Delta\gets 2^{C_1\length(M,q)}$\, \textnormal{as in \Cref{lem:bit-complexity} for an appropriately chosen constant $C_1$}\;
  $T\gets 2^{C_2 n\length(M)}$\, \textnormal{as in \Cref{thm:encoding-length} for an appropriately chosen constant $C_2$}\;
    $(\what M,\what q)$ as in \ref{eq: extended LCP} with $\tilde q$ as in \eqref{eq: initial point formula}\;
    $(\hat x^0,\hat s^0,\hat \mu_0)$ as in \eqref{eq: initial point formula} \;
     $\kappa\gets 1$\;
    \While{$\kappa\leq T$}{
        $\targetmu\gets 1/(16384 n^3\Delta^2(\kappa+1)^2)$ as in \Cref{lem:terminate-kappa}\;
        $R\gets$ \FixedKappaNewton{$\what M,\what q,\hat x^0,\hat s^0,\hat \mu_0,\targetmu,\kappa$}\;
        \lIf{$R=(\what M\notin\cPs)$}{\Return{$M\notin\cPs$}}
        \If{$R=((x,\tilde x),(s,\tilde s),\mu)$}{
        $I\gets\{i\in [n]\,:\, x_i\le \sqrt{2n\mu}\}$ ; $J\gets\{i\in [n]\,:\, s_i\le \sqrt{2n\mu}\}$\label{line:I-J}\;
        $U\gets$ \LPSolve{$Mx+q=s\,,\, x,s\ge \0_n\,,\,x_I=\0_I\,,\,s_J=\0_J$}\;
        \lIf{$U=(x',s')$}{\Return{$(x',s')$}}
        }
        $\kappa\gets 2\kappa$\label{line:update-kappa}\;
    }
    \Return{$M\notin\cPs$}
\end{algorithm}
Now we can prove our main theorem \Cref{thm:main}.
\main*
\begin{proof}
In \Cref{alg:rescale-ipm-full}, we first check the feasibility of \ref{eq: LCP} by solving the linear relaxation $Mx+q=s$, $x,s\ge\0_n$. 
We terminate if this is infeasible. Otherwise, either \ref{eq: LCP} is feasible, or $M\notin\cPs$ must hold, according to  \Cref{thm:feasibile}. 

Let us start with analyzing the subroutine \textsc{FixedKappaNewton}. This returns $M\notin\cPs$ whenever \eqref{eq: Newton-step for basic IPM} has no solution or non-unique solution. The correctness of this conclusion is justified by \Cref{cor:invertiblity}. By the choice of $\beta$ and $\delta$ as in \eqref{eq:greek}, \Cref{thm:Newton-simple} and \Cref{lem:two-step-progress-rump} guarantee that as long as $\hk(M)\leq\kappa$, the algorithm either terminates in 
\begin{equation}\label{eq:kappa-running}
O\left(\min\left\{\sqrt{n} (1+\kappa),n^{9/8}\sqrt{1+\kappa}\right\}\cdot \log\left(\frac{\mu_0}{\mu_1}\right)\right)
\end{equation}
iterations. Note that \Cref{thm:Newton-simple} analyzes \Cref{alg:rescale-ipm-simple} that only does a single Newton step whereas here we perform two Newton steps; but this does not increase the asymptotic running time.

Consider now the cases when the subroutine returns $\hk(M)>\kappa$. This happens when either the iterate $(x', s')$ after the first Newton step falls outside of $\Feasp$, or the iterate $(x_g, s_g, \mu_g)$ after the second Newton step falls outside of $\cN^2(\beta)$. Both analyses  guarantee that $(x', s')\in\Feasp$ and $(x_g, s_g, \mu_g)\in\cN^2(\beta)$ respectively, under the assumption that $M\in\cPs(\kappa)$; in the case of \Cref{lem:two-step-progress-rump}, we recall that $\cNc(\beta)\subseteq \cN(\beta)$ by \eqref{eq:rump-vs-narrow}. From the failure of the analysis, it follows that the assumption $M\in\cPs(\kappa)$ does not hold, and hence we correctly infer $\hk(M)>\kappa$. Note that in the parameter choice $\kappa^4>n^5$ the iterates could leave $\cNc(\beta)$ but we do not notice as we only check containment in $\cN(\beta)$. Hence, the algorithm may be able to proceed for longer, until the broader neighborhood also get violated; or it might even succeed in reaching $\targetmu$ despite $M\notin\cPs(\kappa)$.

Consider now the embedded system \ref{eq: extended LCP} and $(\hat x^0,\hat s^0,\hat \mu_0)$ as in \Cref{lem: initial point} and the choice of $\targetmu$ as in \Cref{lem:terminate-kappa}. By 
 \Cref{lem: initial point}, $\hk(\what M)=\hk(M)$. Hence, if \textsc{FixedKappaNewton}$(\what M,\what q,\hat x^0,\hat s^0,\hat \mu_0,\targetmu,\kappa)$ terminates with $\what M\notin\cPs$ or $\hk(M)>\kappa$, these conclusions are also valid for $M$. Further, notice that $\log(\mu_0/\mu_1)=O(L+\log(2+ \kappa))$.

 In case that \textsc{FixedKappaNewton}$(\what M,\what q,\hat x^0,\hat s^0,\hat \mu_0,\targetmu,\kappa)$ returns a feasible $((x,\tilde x),(s,\tilde s),\mu)\in \cN^2(\beta_{2n}(\kappa))$ with $\mu\le\mu_1$, then \Cref{lem:terminate-kappa} guarantees that either $\hk(M)>\kappa$, or there is a complementary feasible solution $(x',s')$ to \ref{eq: LCP} with $x'_I=\0_I$ and $s'_J=\0_J$ for $I$ and $J$ defined as in line~\ref{line:I-J}. Hence, in case no feasible solution is returned, increasing $\kappa$ is again justified.

 Since $\kappa$ increases by a factor 2 between any two iterations, and the running time bound of each iteration is \eqref{eq:kappa-running} with $\log(\mu_0/\targetmu)=O(L+\log(2+\kappa))$, where this bound also applies to unsuccessful guess, we can conclude that the running time of the final iteration dominates all preceding iterations. Since $\kappa\le 2\hk(M)+1$, the overall running time bound follows, noting also that there are at most $O(\log(2+\hk(M)))$ calls to the LP solver.
 \end{proof}

\section{Deferred Proofs}\label{sec:proofs}

\subsection{Proof of Theorem~\ref{thm:equivalent_pstark_char}}\label{sec:proof-of-equikappathm}

First we restate the theorem.

\equikappathm*

The $\le$ direction is proved in  Corollary~\ref{cor:inverse-norm}, and the $\ge$ direction in Lemma~\ref{lem:pstark_reciprocal}. 

\begin{proposition}%
\label{prop:scalar-ineq} Let $\kappa\ge 0$, and let $P,N \geq 0$ be scalars. Then,  $\left(1+4\kappa\right)P-N\geq0$ if and only if $P-N\geq-\frac{2\kappa}{1+2\kappa}\left(P+N\right)$.
\end{proposition}

\begin{proof} If $\kappa = 0$, then the statement is trivial, hence assume $\kappa>0$.
Let $\alpha=\frac{2\kappa}{1+2\kappa}$, then $\frac{1+\alpha}{1-\alpha}=1+\frac{2\alpha}{1-\alpha}=1+\frac{2\cdot2 \kappa/\left(1+2\kappa\right)}{1-2\kappa/\left(1+2\kappa\right)}=1+4\kappa$.
Hence, $\left(1+4\kappa\right)P-N\geq0$ is equivalent to
\[
\frac{1+\alpha}{1-\alpha}P-N\geq0
\,.\]
By re-arranging terms,
$P-N\geq-\alpha\left(P+N\right)$,
which is further equivalent to $P-N\geq-\frac{2\kappa}{1+2\kappa}\left(P+N\right)$.
\end{proof}
\begin{lem}
\label{lem:lower bound on M+M^T} If $M\in\cPs(\kappa)\cap \R^{n \times n}$, then 
 \[
 M+M^\top%
 \succeq -\left(\frac{2\kappa}{1+2\kappa}\right)^2 \cdot I_n -  M^{\top}M\,.
  \]
\end{lem}

\begin{proof}
If $\kappa=0$, then the lemma is trivial, hence assume $\kappa > 0$.
    Let $x\in \R^n$, and $P=\sum_{i:x_{i}\left(Mx\right)_{i}\geq0}x_{i}\left(Mx\right)_{i}$,
and $N=-\sum_{i:x_{i}\left(Mx\right)_{i}<0}x_{i}\left(Mx\right)_{i}$. Since $(1+4\kappa)P-N\geq 0$, by Proposition \ref{prop:scalar-ineq}, we have  
\begin{align*}
P-N & \geq-\frac{2\kappa}{1+2\kappa}\left(P+N\right)\, ,
\end{align*}
where $x^\top(M+M^\top)x=2x^\top Mx=2(P-N)$. We can upper bound
\begin{align*}
P+N&=\sum_{i=1}^n\left|x\right|_{i}\left|Mx\right|_{i}\leq\sum_{i=1}^n\frac{1}{2}\left( \frac{2\kappa}{1+2\kappa}\cdot  x_{i}^{2}+\frac{1+2\kappa}{2\kappa}\cdot \left(Mx\right)_{i}^{2}\right)
\\
&=\frac{1}{2}\left( \frac{2\kappa}{1+2\kappa}\cdot  \left\Vert x\right\Vert^{2}+\frac{1+2\kappa}{2\kappa}\cdot\left\Vert Mx\right\Vert^{2}\right)\,.
\end{align*}
Combining the two inequalities yields 
\[
x^\top(M+M^\top)x\ge -\left(\frac{2\kappa}{1+2\kappa}\right)^2 x^\top  x-x^\top M^\top Mx\, \quad\forall x\in\R^n\, ,
\]
proving the claim.%
\end{proof}

\begin{cor}
\label{cor:inverse-norm} If $M\in \cPs(\kappa)\cap\R^{n\times n}$,
then for any $D\in\cDpp$, 
\[
\left\Vert \left(I_n+DMD\right)^{-1}\right\Vert_2 \leq\sqrt{  1+\frac{\kappa^2}{1/4+\kappa }}\,.
\]
\end{cor}
\begin{proof}
Note that by \Cref{lem:DMD}, $M\in\cPs(\kappa)$ if and only if $DMD\in\cPs(\kappa)$.
Hence, it suffices to show the statement for $D=I_n$, that is, $\left\Vert \left(I_n+M\right)^{-1}\right\Vert_2 \leq
\sqrt{  1+\frac{\kappa^2}{1/4+\kappa }}$.
From Lemma \ref{lem:lower bound on M+M^T} we can lower bound
\[
\left(I_n+M\right)^{\top}\left(I_n+M\right)=I_n+M^{\top}M+M+M^{\top}
\succeq I_n+M^\top M  -\left(\frac{2\kappa}{1+2\kappa}\right)^2 \cdot I_n -  M^{\top}M
=  \frac{1+4\kappa}{(1+2\kappa)^2}  \cdot I_n\,.
\]
Equivalently we have
\[
\left\Vert \left(I_n+M\right)^{-1}\right\Vert_2 \leq
\sqrt{  1+\frac{\kappa^2}{1/4+\kappa }}\,.
\]
which is what we needed.
\end{proof}
We now show the converse direction, completing the proof of \Cref{thm:equivalent_pstark_char}.
\begin{lem}\label{lem:pstark_reciprocal}
Suppose that for $M\in\R^{n\times n}$ there exists a scalar $\kappa  \geq 0 $ such that for any
$D\in\cDpp$, 
\[
\left\Vert \left(I_n+DMD\right)^{-1}\right\Vert_2 
\leq 
\sqrt{  1+\frac{\kappa^2}{1/4+\kappa }}
\,,
\]
then $M$ is a $\cPs(\kappa)$ %
matrix, that is, $\hk(M)\le \kappa$.
\end{lem}

\begin{proof}
Define $\varepsilon = 1/\left( 1+\frac{ \kappa^2}{1/4+\kappa }\right) = \frac{1+4\kappa}{1+4\kappa+4\kappa^2}$, so that $\left\Vert \left(I_n+DMD\right)^{-1}\right\Vert_2 \leq \sqrt{1/\varepsilon}$.
Equivalently we can write the condition as 
\[
\left(1-\varepsilon\right)I_n+D\left(M+M^{\top}\right)D+DM^{\top}D^{2}MD\succeq0
\]
and again equivalently as
\[
\left(1-\varepsilon\right)D^{-2}+M^{\top}D^{2}M+M+M^{\top}\succeq0\,.
\]
Thus for any $x\in\R^n$ and $D\in\cDpp$, \begin{equation}\label{eq: to derive kappa by spectrual norm}
\sum\limits_{i=1}^n\left(\left(1-\varepsilon\right)\frac{x_{i}^{2}}{d_{i}^{2}}+\left(Mx\right)_{i}^{2}\cdot d_{i}^{2}\right)+2x^{\top}Mx\geq0
\,.  
\end{equation}
Now we take the infimum over all strictly positive $d_i$. For any  $a,b \geq 0$, we have $\inf_{t>0} \left( \frac{a}{t} + bt \right) = 2\sqrt{ab}$. Indeed, if both $a,b>0$, the infimum is attained at $t=\sqrt{a/b}$. If $a = 0$, the same identity follows by taking $t\rightarrow 0$. If $b=0$, we take $t \rightarrow \infty$. If $a=b=0$, any positive $t$ attains the infimum. Note that in the cases where the infimum is not attained, it approaches $0$ in the limit.

Therefore, applying this coordinatewise with $t= d_i^2$, $a=(1-\varepsilon)x_i^2$, $b= (Mx)_i^2$ in \eqref{eq: to derive kappa by spectrual norm} we obtain that
\[
\sqrt{1-\varepsilon}\sum\limits_{i=1}^n\left|x_{i}\right|\left|Mx\right|_{i}+x^{\top}Mx\geq0
\]
Again letting $P=\sum_{i:x_{i}\left(Mx\right)_{i}\geq0}x_{i}\left(Mx\right)_{i}$,
and $N=-\sum_{i:x_{i}\left(Mx\right)_{i}<0}x_{i}\left(Mx\right)_{i}$, this equivalently
gives us
\[
\left(P-N\right)\geq-\sqrt{1-\varepsilon}\left(P+N\right)\, .
\]
By \Cref{prop:scalar-ineq}, we  equivalently have
that $\left(1+4\tilde{\kappa}\right)P-N\geq0$ for $\tilde{\kappa}$ such that $\sqrt{1-\varepsilon}=\frac{2\tilde{\kappa}}{1+2\tilde{\kappa}}=1-\frac{1}{1+2\tilde{\kappa}}$,
thus $2\tilde{\kappa}=\frac{1}{1-\sqrt{1-\varepsilon}}-1$. Substituting $\varepsilon$ into this expression yields
\[
2\tilde{\kappa} = \frac{1}{1-\sqrt{ \frac{4 \kappa^2}{(1+2\kappa)^2} }}-1 = 2\kappa\, ,\]
which shows that $M$ is in $\cP_\star(\kappa)$, as claimed.
\end{proof}

\subsection{Proof of Theorem~\ref{lem:two-step-progress-rump}}\label{sec:two-step-no-recentering}
\thmimprovedanalysis*
\begin{proof}
We keep the same notation as in the proof of \Cref{lem:two-step-progress}, and recall  $\|Q_{H_{x,\mu}}\err(x, s, \mu)\|_2\le\beta$. Let $\Delta \bx \coloneqq p + r$, where $p\coloneqq-\delta H\1_n$ and $r\coloneqq(1+\delta)H\err(x,s,\mu)$.

Applying \Cref{lem:main-stability} to the change from $\mu$ to $\mu'$ gives $q_H\le\exp(8\Kp n\delta)q_{H_{x,\mu}}$, so $\|(1+\delta)Q_H\err(x, s, \mu)\|_2\le\exp(65/1024)\beta\le3\beta/2$. Since $p=-\delta H\1_n$ and $\|\1_n\|_2=\sqrt n$, the arguments used for \eqref{eq:two-step-compact-first} and \eqref{eq:two-step-first-matrix} give
\[
\begin{aligned}
    \|p\|_2&\le\delta\sqrt{n(1+\kappa)},\\
    \|Q_Hp\|_\infty&\le32\delta\sqrt{n(1+\kappa)},\\
    \|H\diag{p}\|_2&\le2\delta n\sqrt{1+\kappa}.
\end{aligned}
\]
For $r$, we instead use $\|Q_HHQ_H^{-1}\|_{\max}\le\Kp=8$ from \Cref{lem:weights-well-defined} and $q_H\ge\1$ to obtain bounds independent of $\kappa$:
\[
\begin{aligned}
    \|Q_Hr\|_\infty&\le12\sqrt n\,\beta,\\
    \|r\|_2 &\le12n\beta,\\
    \|H\diag{r}\|_2&\le\|H\diag{r}\|_F
    \le8\sqrt n\,\|Q_Hr\|_2
    \le96n^{3/2}\beta.
\end{aligned}
\]
These replace the use of \eqref{eq:two-step-smallness}. In particular, with $U=\diag{\1_n+\Delta\bx}$ as before, using the bounds on $\delta$ and $\beta$,
\[
\begin{aligned}
    \|\Delta\bx\|_1&\le\delta n\sqrt{1+\kappa}+12n^{3/2}\beta\le\frac1{256},\\
    \|H(U^{-2}-I_n)\|_2&\le4\|H\diag{\Delta\bx}\|_2
    \le8\delta n\sqrt{1+\kappa}+384n^{3/2}\beta\le\frac{13}{128}<\frac12.
\end{aligned}
\]
Thus $x',s'>0$ by \eqref{eq:two-step-product-identity}, and the argument leading to \eqref{eq:two-step-matrix-transfer} still gives $\|H'\diag{p}\|_2\le4\|H\diag{p}\|_2$. To compare the weights, use $|\log(1+a)|\le\frac43|a|$ for $|a|\le\frac14$ in \Cref{lem:main-stability}. This gives $\|\log q_{H'}-\log q_H\|_\infty\le\frac{64}{3}\Kp\|\Delta\bx\|_1\le\frac23<\log2$, so again $q_{H'}\le2q_H$.

For the second step, write $\Delta\bx'=z+w$, where $z\coloneqq H'p^{2}$ and $w\coloneqq H'(2p r+r^{2})$. The bounds for $z$ are the same as in \eqref{eq:two-step-compact-second} and \eqref{eq:two-step-second-infty}:
\[
\begin{aligned}
    \|z\|_2&\le4\|H\diag{p}\|_2\|p\|_2
    \le8\delta^2n^{3/2}(1+\kappa),\\
    \|Q_{H'}z\|_\infty&\le2\Kp\|Q_Hp\|_\infty\|p\|_1
    \le512\delta^2n^{3/2}(1+\kappa).
\end{aligned}
\]
To bound the terms involving $r$, we use the entry bound for $Q_{H'}H'Q_{H'}^{-1}$ and $q_{H'}\le2q_H$ to get 
\[
\begin{aligned}
    \|Q_{H'}w\|_\infty
    &\le2\Kp\|Q_Hr\|_\infty(2\|p\|_1+\|r\|_1)\\
    &\le384\delta n^{3/2}\sqrt{1+\kappa}\,\beta+2304n^2\beta^2
    \le\left(\frac38+\frac9{16}\right)\sqrt n\,\beta
    \le\sqrt n\,\beta.
\end{aligned}
\]
Consequently,
\[
    \|\Delta\bx'\|_1
    \le\sqrt n\,\|z\|_2+n\|Q_{H'}w\|_\infty
    \le8\delta^2n^2(1+\kappa)+n^{3/2}\beta
    \le\frac{33}{2^{17}}<\frac1{512}.
\]
As in the proof of \Cref{lem:two-step-progress}, this ensures $x'',s''>0$ and $q_{H''}\le2q_{H'}$. Finally, applying $(z+w)^{2}\le2z^{2}+2w^{2}$ gives the following strengthening of  \eqref{eq:two-step-final-error}:
\begin{equation}\label{eq:two-step-rump-initial-error}
\begin{aligned}
    \|Q_{H''}\err(x'',s'',\mu')\|_2
    &\le4\left(\|Q_{H'}z\|_\infty\|z\|_2
        +\sqrt n\,\|Q_{H'}w\|_\infty^2\right)\\
    &\le2^{14}\delta^4n^3(1+\kappa)^2+4n^{3/2}\beta^2\\
    &\le\frac{\beta}{2^{14}}+\frac{\beta}{1024}<\beta,
\end{aligned}
\end{equation}
where in the last inequality we use the intervals of $\delta$ and $\beta$, in particular $\delta\leq \frac{\beta^{1/4}}{128n^{3/4}\sqrt{1+\kappa}}$.
Both steps preserve linear feasibility, so $(x'',s'',\mu')\in\cNc(\beta)$.
\end{proof}

\subsection{Deferred proofs from Section~\ref{sec:hardness}}\label{sec:deferred_hardness}
In this section, we prove \Cref{lem: rational gap lemma}, \Cref{thm: closed Rump's theorem}, \Cref{prop:change-of-variable-handicap-number} and \Cref{prop:identity-for-max-cut}. 
\rationalgap*
\begin{proof}
 Let $S\in\cS$, $(\lambda_j)_{j=1}^n$ denote all the eigenvalues of $SM$, and $\lambda_i$ with $i\in [n]$ be any real eigenvalue of $SM$. Since $\rhos(M)<1$, we have $|\lambda_i|<1$ and $\det(\pm I_n-SM)\neq 0$. We first suppose that $\lambda_i\geq 0$. Since $\det(I_n- SM)$ is an integer multiple of $1/C^n$, $\prod_{j=1}^n |1-\lambda_j|=|\det(I_n-SM)|\geq 1/C^n$. For all $j\in [n]$, $|\lambda_j|\leq\|SM\|\leq \|M\|_F\leq n\|M\|_{\max}$. Then, \[1-\lambda_i\geq \frac{1}{C^n\prod\limits_{j\neq i}|1-\lambda_j|}\geq \frac{1}{C^n\prod\limits_{j\neq i}(1+|\lambda_j|)}\geq \frac{1}{C^n(1+n\|M\|_{\max})^{n-1}}\] giving $|\lambda_i|\leq 1-\frac{1}{C^n(1+n\|M\|_{\max})^{n-1}}=1-\frac{1}{R}$. If $\lambda_i<0$ we instead use $\det(-I_n-SM)$ being an integer multiple of $1/C^n$, which gives $\prod_{j=1}^n |1+\lambda_j|\geq 1/C^n$. Then, \[1+\lambda_i\geq \frac{1}{C^n\prod\limits_{j\neq i}|1+\lambda_j|}\geq \frac{1}{C^n\prod\limits_{j\neq i}1+|\lambda_j|}\geq \frac{1}{C^n(1+n\|M\|_{\max})^{n-1}}\, , \] 
which also gives $|\lambda_i|\leq 1-\frac{1}{C^n(1+n\|M\|_{\max})^{n-1}}=1-\frac{1}{R}$. We have \[\log R=n\log C+(n-1)\log(1+n\|M\|_{\max})=O(\poly(n, \length(M)))\,.\] 

\end{proof}
\closedRump*
\begin{proof}
    Suppose that $\rhos(M)\leq 1$. Let $(\alpha_k)_{k\in\N}=1-\frac{1}{k+1}$, we have $\rhos(\alpha_kM)=\alpha_k\rhos(M)\leq \alpha_k<1$, and also, $I_n-\alpha_kM$ is nonsingular because if otherwise, $M$ has an eigenvalue of $1/\alpha_k>1$ which contradicts $\rhos(M)\leq 1$. Then using \Cref{thm: Rump's theorem}, $\cC(\alpha_k M)=(I_n-\alpha_k M)^{-1}(I_n+\alpha_k M)\in \cP$. Thus, $\lim_{k\rightarrow\infty}\cC(\alpha_kM)\in\cP_0$ by \Cref{lem: matrix class relation}, where $\lim_{k\rightarrow\infty}\cC(\alpha_k M)=\lim_{k\rightarrow \infty}(I_n-\alpha_k M)^{-1}(I_n+\alpha_k M)=(I_n-M)^{-1}(I_n+M)=\cC(M)$.

    Suppose that $\cC(M)\in\cP_0$. We have $\cC(M)+\varepsilon I_n\in\cP$ for any $\varepsilon>0$ by \Cref{lem: matrix class relation} and $(1+\varepsilon)I_n+\cC(M)$ is nonsingular by \Cref{cor:invertiblity}, so $N_\varepsilon:=\cC^{-1}(\cC(M)+\varepsilon I_n)$ exists. As a result of \Cref{prop: inverse Cayley transform} and \Cref{thm: Rump's theorem}, $\rhos(N_\varepsilon)<1$. We also have $\lim_{\varepsilon\rightarrow 0}\cC(M)+\varepsilon I_n=\cC(M)$ and \[\lim_{\varepsilon\rightarrow 0}N_\varepsilon=\lim_{\varepsilon\rightarrow 0} (\cC(M)+\varepsilon I_n-I_n)(\cC(M)+\varepsilon I_n+I_n)^{-1}=(\cC(M)-I_n)(\cC(M)+I_n)^{-1}=\cC^{-1}(\cC(M))=M\,.\]
    Therefore, it follows from the continuity of $\rhos(X)$ in $X$ entry-wise \cite[Corollary 2.5]{rump1997theorems} that \[\rhos(M)=\rhos\left(\lim_{\varepsilon\rightarrow 0} N_\varepsilon\right)=\lim_{\varepsilon\rightarrow 0} \rhos(N_\varepsilon)\leq 1\,.\]
\end{proof}
\changeofvariable*
\begin{proof}
    By \Cref{lem: matrix class relation}, $I_n+M$ is invertible so $y=\frac{1}{2}(I_n+M)x$ is a bijective linear map. We have \begin{equation*}
        \begin{aligned}
         Mx&=2M(I_n+M)^{-1}y=2\left(I_n-(I_n+M)^{-1}\right)y\\
         &=2\left(I_n-\frac{1}{2}(I_n+N)^{-1}\right)y%
    =\left(I_n+N(I_n+N)^{-1}\right)y\, ,  
        \end{aligned}
    \end{equation*}
    and $x=2(I_n+M)^{-1}y=(I_n+N)^{-1}y=(I_n-N(I_n+N)^{-1})y$. Using $y^\top(I_n-A)^\top(I_n+A)y=\|y\|^2-\|A y\|^2$ for any $A\in\R^{n\times n}$, we have \begin{equation*}
        x^\top Mx=y^\top\left(I_n-N(I_n+N)^{-1}\right)^\top\left(I_n+N(I_n+N)^{-1}\right)y=\|y\|^2-\|N(I_n+N)^{-1}y\|^2
    \end{equation*}
    When $N$ is nilpotent of degree $2$, i.e., $N^2=\0_{n\times n}$, $(I_n+N)(I_n-N)=I_n$, so $(I_n+N)^{-1}=I_n-N$, thus we have \[x^\top Mx=\|y\|^2-\|N(I_n-N)y\|^2=\|y\|^2-\|Ny\|^2=y^\top(I_n-N^\top N)y\,.\]
\end{proof}
\maxcutidentity*
\begin{proof} If $v=\0$ then the statement is trivial, hence assume $v\neq \0$.
 Let $\hat v=v/\|v\|_\infty$, then \begin{equation*}
    v^\top (L_G-I_n)v=\|v\|_\infty^2 \cdot \hat v^\top(L_G-I_n)\hat v\leq \|v\|_\infty^2\cdot \max_{s\in [-1, 1]^n} s^\top (L_G-I_n)s\, .
\end{equation*}
The maximum of $s\in [-1,1]^n$ is taken at a vertex $s\in \{-1,1\}^n$, by \cite[Corollary 8.1]{punnen2022}. To see the proof of this, note that
 $(L_G - I_n)_{ii}\geq 0$ for all $i\in [n]$ because $G$ is connected (and thus every vertex has degree $\ge 1$). Hence, the function $f(s)=s^\top (L_G-I_n) s$ is convex in each coordinate direction $s_i$. This implies that there is a maximizer $s\in \{-1,1\}^n$.
Consequently,
\begin{equation*}
 \begin{aligned}
       v^\top (L_G-I_n)v\leq  \|v\|_\infty^2\cdot \max_{s\in \{-1, 1\}^n} s^\top (L_G-I_n)s=\|v\|_\infty^2\left(4\mc(G)-n\right)\, ,
\end{aligned}
 \end{equation*} where the last equality follows by \Cref{prop:maxcut-laplacian}.
\end{proof}
\zeroeigencor*
\begin{proof}
    Suppose $0$ is not a semi-simple eigenvalue of $M$. Then there exists $x, y\in\R^n$, $x\neq \0$ such that $Mx=\0$ and $My=x$. Since $M\in\cPs$, $I_n+M$ is nonsingular by \Cref{cor:invertiblity}. Let $\wt M=\cC^{-1}(M)=(M-I)(M+I)^{-1}$, then $\wt Mx=-x $ and $\wt My=(M-I)\left(y-(I+M)^{-1}x\right)=(M-I)(y-x)=2x-y$. Since $|\wt Mx|=|x|$ and \[\wt Mx \circ \wt M(-y)-x\circ (-y)=(-x)\circ (y-2x)+x\circ y=2x^2\geq \0\,,\]
    \Cref{thm:sufficient-Rump-theorem} gives $\wt Mx\circ \wt My=x\circ y$, which implies that $x=\0$, a contradiction. As a result, $M^p$ has the same null space and rank as $M$, for all $p\in\mathbb{N}$.

  Let $x\in\ker(M)$, so $\wt Mx=-x$. \Cref{thm:sufficient-Rump-theorem} ensures that if $(-x)\circ \wt My\geq x\circ y $ then $(-x)\circ \wt My=x\circ y$, for all $y\in\R^n$
. If we set $A_x:=-\diag{x}(I+\wt M)$, this is equivalent to having for all $y\in\R^n$, if $A_xy\geq\0$ then $A_x y=\0$, that is, there exists no $y\in\R^n$ such that $A_xy\geq \0$ and $A_xy\neq \0$. By Stiemke's alternative, there exists $w>\0$ such that $A_x^\top w=\0$. Let $\tilde x:=\diag{w}x$, then $(I+\wt M^\top)\tilde x=(I+\wt M^\top)\diag{x}w=-A_x^\top w=\0$ which gives $\tilde x\in \ker(I+\wt M^\top)$. Since \[I+\wt M^\top=I+(M^\top+I)^{-1}(M^\top-I)=2(M^\top+I)^{-1}M^\top\,,\] $\ker(I+\wt M^\top)=\ker(M^\top)$. Thus, we have $\tilde x\in\ker(M^\top)$, and $\sgn(\tilde x)=\sgn(x)$ because $\tilde x=\diag{w} x$ where $w>\0$. This shows that $\sgn(\ker(M))\subseteq \sgn(\ker(M^\top))$. We then have $\sgn(\ker(M))=\sgn(\ker(M^\top))$ by applying the same argument to $M^\top$ to show $\sgn(\ker(M^\top))\subseteq\sgn(\ker(M))$, because $M^\top\in\cPs$ by \Cref{thm: handicap number of matrix transpose}. \end{proof}

\subsection*{AI Disclosure} 
The authors uploaded to ChatGPT a writeup of their own simple IPM analysis from Section 4 and characterization from \Cref{thm:equivalent_pstark_char}, seeking to avoid the linear dependence on the handicap number present in this analysis and in earlier work. Prompted to test a range of approaches, ChatGPT suggested measuring the Newton step $H \err$ in a
norm induced by an M-matrix inverse $(I-|H|/\alpha)^{-1}$ for a suitable scalar $\alpha$. Further refinement with the authors
led ChatGPT to suggest using Rump's theorem to provide an analysis of Newton's method that uses the local rescaling rather than the global worst-case guarantee on the norm of $H$. The authors developed the final analysis with LLM assistance, particularly the stability proofs from Sections~\ref{sec:bounding-Rump-weights} and~\ref{sec:stability-Rump-weights}. 

DeepSeek assisted with the hardness results, and independently suggested using Rump's sign spectral radius to
prove coNP-hardness of testing sufficiency. The authors refined all ideas and strengthened the arguments. 

We view the connection to Rump's work, which appears to have been overlooked in the IPM literature, as the most significant LLM contribution. LLM use in writing was limited to proofreading and helping to improve the flow of the text.

\bibliographystyle{alphaabbrv}

\bibliography{2VPI}
\end{document}